\documentclass{article}

\usepackage{arxiv}

\usepackage[utf8]{inputenc} 
\usepackage[T1]{fontenc}    
\usepackage{hyperref}       
\usepackage{url}            
\usepackage{booktabs}       
\usepackage{amsfonts}       
\usepackage{nicefrac}       
\usepackage{microtype}      
\usepackage{lipsum}		
\usepackage{graphicx}
\usepackage{natbib}
\usepackage{doi}

\usepackage{enumitem}
\usepackage{tikz}
\usetikzlibrary{calc,arrows.meta,positioning}

\usetikzlibrary{shapes,decorations,arrows,calc,arrows.meta,fit,positioning}
\tikzset{
    -Latex,auto,node distance =1 cm and 1 cm,semithick,
    state/.style ={circle, draw, minimum width = 0.5 cm},
    point/.style = {circle, draw, inner sep=0.04cm,fill,node contents={}},
    bidirected/.style={Latex-Latex,dashed},
    el/.style = {inner sep=2pt, align=left, sloped}
}

\makeatletter
\newcommand{\xleftrightarrow}[2][]{\ext@arrow 3359\leftrightarrowfill@{#1}{#2}}
\makeatother

\newcommand{\xdasharrow}[2][->]{
\tikz[baseline=-\the\dimexpr\fontdimen22\textfont2\relax]{
\node[anchor=south,font=\scriptsize, inner ysep=1.5pt,outer xsep=2.2pt](x){#2};
\draw[shorten <=3.4pt,shorten >=3.4pt,dashed,#1](x.south west)--(x.south east);
}
}

\usepackage{longtable}
\usepackage{soul}

\usepackage{amsmath}
\usepackage{amsthm}
\usepackage{amsfonts}
\usepackage{verbatim}

\usepackage{times}
\usepackage{latexsym}
\usepackage{graphicx}
\usepackage{color}
\usepackage{dsfont}

\usepackage{multirow}
\usepackage{float}

\usepackage{enumitem}

\usepackage{booktabs} 

\usepackage{subfig}

\newtheorem{thm}{Theorem}[section]

\newtheorem{prop}[thm]{Proposition}
\newtheorem{ex}[thm]{Example}
\newtheorem{defn}[thm]{Definition}
\newtheorem{rem}[thm]{Remark}

\newcommand{\argmin}{\mathrm{argmin}}

\newcommand{\field}[1]{\mathbb{#1}}
\DeclareMathOperator{\E}{\field{E}}      

\def\A{\field{A}}        
\def\X{\field{X}}        

\newcommand{\ls}[1]
  {\dimen0=\fontdimen6\the\font \lineskip=#1\dimen0
  \advance\lineskip.5\fontdimen5\the\font \advance\lineskip-\dimen0
  \lineskiplimit=.9\lineskip \baselineskip=\lineskip
  \advance\baselineskip\dimen0 \normallineskip\lineskip
  \normallineskiplimit\lineskiplimit \normalbaselineskip\baselineskip
  \ignorespaces }

\graphicspath{{./}{Figs/}}

\title{Optimizing Patient Transitions to \\ Skilled Nursing Facilities}

\author{ {Oguzhan Alagoz}\\
	Department of Industrial and Systems Engineering \\
	University of Wisconsin-Madison\\
	Madison, WI \\
	\texttt{alagoz@engr.wisc.edu} \\
	\And
	{Sebastian A. Alvarez Avendaño} \\
	Department of Industrial and Systems Engineering\\
	University of Wisconsin-Madison\\
	Madison, WI \\
	\texttt{alvarezavend@wisc.edu} \\
	\And
	{John Oruongo} \\
	Sykes College of Business\\
	University of Tampa\\
	Tampa, FL \\
	\texttt{joruongo@ut.edu} \\
	\And
	{Gabriel Zayas-Cab\'{a}n}\\ 
	Department of Industrial and Systems Engineering\\
	University of Wisconsin-Madison\\
	Madison, WI \\
	\texttt{zayascaban@wisc.edu} \\
}

\renewcommand{\headeright}{Preprint}
\renewcommand{\undertitle}{Preprint}
\renewcommand{\shorttitle}{Optimizing Patient Transitions to Skilled Nursing Facilities}

\hypersetup{
pdftitle={Optimizing Patient Transitions to Skilled Nursing Facilities},
pdfauthor={Sebastian A. Alvarez Avendaño, John Oruongo, Oguzhan Alagoz, Gabriel Zayas-Cab\'{a}n},
pdfkeywords={care transitions; hospitals; skilled nursing facilities; readmissions; stochastic modeling applications},
}

\begin{document}
\maketitle

\begin{abstract}
Hospital inpatient care costs is the largest component of health care expenditures in the US. At the same time, the number of non-hospital rehabilitative settings, such as skilled nursing facilities (SNFs), has increased. Lower costs and increased availability have made SNFs and other non-hospital rehabilitation settings a promising care alternative to hospitalization.  To maximize their benefits, transitions to SNFs require special attention, since poorly coordinated transitions can lead to worse outcomes and higher costs via unnecessary hospital readmissions. This study presents a framework to improve care transitions based on the premise that certain SNFs may provide better care for some patients. We estimate readmission rates by SNF and patient types using observational data from a tertiary teaching hospital and nearby SNFs. We then analyze and solve a stochastic model optimizing patient transfer decisions to minimize readmissions. Our model accounts for patient discharge patterns and SNF capacity availability. We provide conditions for when an easy-to-use myopic policy, which assigns discharged patients to the SNF with the lowest readmission rate that is available, is optimal. We also show when an optimal policy has a threshold-like structure. Using estimated readmission rates, we compare the performance of the myopic policy and a proposed policy that depends on the discharge process, readmission rates, and future SNF availability. We evaluate when the myopic policy may be beneficial to use and when the proposed transfer heuristic provides a better alternative. Otherwise, we contend that using a stochastic optimization model for guiding transfer decisions may help reduce readmissions. 
\end{abstract}


\section{Introduction} 

Inadequate care coordination was associated with \$25 to \$45 billion in wasteful spending in 2011 through avoidable complications and with unnecessary hospital readmissions \citep{burton2012health}.  This study is on a specific and important component of coordination: transitions of care.  Care transitions refer to moving patients from one health care facility (e.g., inpatient hospital unit) to another (e.g., skilled nursing facility [SNF]) \citep{burton2012health,li2014optimizing}.  As the above numbers suggest, poorly managed transitions have been associated with increased costs and worse health outcomes due to high readmission rates \citep{burton2012health}. 

Improving hospital post-discharge management through more effective use of post-acute care is a national priority \citep{burke2019transitional,ryskina2020hospitalist}. Under traditional fee-for-service reimbursement, hospitals had few incentives to invest in reducing readmission rates \citep{zohrabian2018economic}. Since the passage of the Patient Protection and Affordable Care Act, hospitals are increasingly being held accountable for patient outcomes after discharge \citep{andritsos2018incentive,ryskina2019assessing,ryskina2020hospitalist}. One way this has occurred is through a shift from fee-for-service to value-based purchasing, bundled payments, and Accountable Care Organizations. Another way is through financial penalties on hospitals with excess readmission rates as part the Medicare's Hospital Readmissions Reduction Program \citep{zuckerman2016readmissions,werner2018trends,ryskina2019assessing}. As a result, hospitals are increasingly relying on post-acute care to reduce readmissions.

SNFs represent the most common setting for post-acute care in the US, but still result in high rates of readmissions. SNFs provide short-term skilled nursing or rehabilitation services after hospital discharge \citep{ryskina2020association}. It is estimated that more than 40\% of patients who require post-acute care after hospitalization are discharged to SNFs, and one in five Medicare beneficiaries are transferred to a SNF after hospitalization \citep{medicare2019report,ryskina2020association}. Rates of readmission from SNFs are high and quite variable between SNFs \citep{neuman2014association,ouslander2016hospital,feder2018they,ryskina2019postacute,burke2020variability}. One in four patients discharged to a SNF is readmitted within 30 days, and two-thirds of these readmissions may be preventable \citep{neuman2014association,vasilevskis2017potentially}.  

One way to reduce readmissions is to transfer patients to the most appropriate SNF. This may be possible, because some SNFs may be better equipped to handle certain patients (e.g., cardiac rehabilitation patients) than others \citep{rahman2017contribution,ryskina2020association}. While many studies focused on reducing readmissions from SNFs \citep{voss2011care,chen2012skilled,king2013consequences,gardner2014implementation,neuman2014association,kripalani2014reducing,li2014optimizing,ouslander2016hospital,davidson2017improving,gilmore2017transitions,mileski2017investigation,moore2017improving,rothberg2017impact,gadbois2019lost,ryskina2019assessing,ryskina2019postacute,ryskina2020association,ryskina2020hospitalist}, 
few have considered data driven optimization approaches for supporting these decisions.  

Our goal is to support the management of care transitions by optimizing the transfer of discharged patients to SNFs. Our objective is to reduce readmission rates, considering that readmission rates are a prominent target in the literature described above. An additional benefit of using readmission rates, over other possible objectives, is that it is both a measure for quality of care and a proxy for cost. Quality of care is an important expectation for patients and/or their families and may take other forms than readmissions. Other measures of quality of care by SNFs include falls and medication error rates. While we focus on readmission rates, these alternative measures of quality are highly correlated with the readmission rates as well. For instance, surgery patients going to a SNF that has a high fall rate may lead to complications that require re-hospitalization. Thus, readmission rates capture elements of the quality of care offered at SNFs. Insurers, on the other hand, are interested in cost savings. Preventable readmissions results in increased health costs, and identifying opportunities to minimize readmissions result in cost savings. The Centers for Medicare and Medicaid Services, or CMS, for example, uses readmissions as a measure for cost. With reducing readmissions as our objective, optimizing transfers to SNFs requires confronting two issues.  First, readmission rates may vary by patient treatment requirements \citep{coleman2003falling,coleman2004lost,coleman2006care}. Second, SNF availability is unpredictable at the time of discharge.  It may depend on patient demographics, clinical complexity, hospital treatment, demand from other hospitals, and other factors.  

In this study, we develop an analytical approach for optimizing care transitions.  It consists of two complementary frameworks that address the two aforementioned issues. The first is estimating readmission rates by patient characteristics using observational data from a large tertiary teaching hospital.  This analysis suggests that different SNFs yield significantly different readmission rates for different patient types depending on demographic factors and diagnosis-related groups (DRGs), which are a patient classification scheme providing a means of relating the type of patients a hospital treats (i.e., its case mix) to the costs incurred by the hospital and treatment \citep{fetter1980case}.  Hospitals may thus reduce readmission rates by choosing the ``right" SNF for a patient. Addressing this latter point, the second framework seeks to optimize transfer decisions using estimated readmission rates in a sequential decision-making model.  For this, we develop and analyze a Markov decision process process (MDP) formulation for optimizing transfer decisions based on the SNF discharge process at our partner hospital that minimizes estimated readmission rates under discounted expected and long-run average cost criteria. Our modeling approach is based on the transfer process utilized by our partner hospital, which is used by many hospitals across the US.

\textbf{Transfer process from hospital to SNFs at our partner hospital.}  Our partner hospital is a large tertiary teaching hospital and health system and transfers patients to SNFs in the surrounding region.  Typically, a care team determines whether the patient is ready for discharge to a SNF. Once ready, the case is assigned to a social worker who then contacts potential SNFs to determine availability for that patient.  Importantly, SNF congestion information is not available to the social worker at the time of SNF transfer decision. Instead, the social worker has information on only whether the SNF is willing (or not) to `accept' the patient needing the SNF at the time of discharge.  The patient is then provided a list of SNFs that have available capacity from which to choose.  The patient chooses one of the SNFs in the list, after which they are discharged to that particular SNF.  Patient preferences play an important factor in choosing the SNF.  For example, patients often prefer SNFs that are closest to their home or the home of their caregivers.  However, patients and their caregivers typically follow the recommendations of their assigned social workers in choosing a SNF. 

There are two common features in post-discharge management.  First, the transfer process in many US hospitals is decentralized. Different medical units such as transplantation and cardiothoracic surgery have their own social workers arranging transfers to SNFs from their unit without coordinating with each another. This has led to efforts to centralize how these decisions are made.  For example, at our partner hospital, administrators are working to centralize the transfer process by coordinating discharges among various units in the hospital. They have also recently established a collaborative working group to help identify the best SNFs caring for patients utilizing historical data on readmissions.  Secondly, SNF availability is not directly controlled by the hospitals. Different types of patients with different care needs from other health systems may also use available capacity and the discharging hospital does not have full knowledge of available SNF at all times. The discharging hospital knows only whether a particular SNF has available capacity at the time of the discharge.  Based on these observations, our model reflects the perspective of a single hospital which uses a centralized controller to make transfer decisions to SNFs after discharge.

\textbf{Contributions to care transitions.}  To the best of our knowledge, our study is the first to propose, analyze, and solve an optimization model for representing transitions from hospitals to the SNFs. We use real-life observational data at the patient-level from a large tertiary teaching hospital to inform our models. In addition, our work makes the following contributions.

We provide conditions for when a particular easy-to-use policy, which we refer to as myopic, is optimal.  The myopic policy assigns discharged patients to the available SNF with the lowest readmission rate for the discharged patient type.  Second, we provide conditions under which an optimal policy has a threshold structure: if it is optimal to assign a patient type to an available SNF, then it is again optimal to send the same patient type to the same SNF when there are fewer available SNFs. The structure of the optimal policy is difficult to discern outside these cases.  This is because the transfer decision depends on readmission rates, which themselves depend on patient type and SNF, and on the discharge process for each patient type and SNF availability, which also depends on the SNF and the transfer decision.  We present examples showing that when the aforementioned conditions do not hold, the myopic policy is not optimal.  Lastly, our analysis highlights that when transfer decisions have little to no impact on future SNF availability, the myopic policy performs close to optimal.  Otherwise, when the transfer decision significantly impacts SNF availability, it performs poorly.  We formally examine this tradeoff in a numerical study.  We compare the performance of the myopic policy and a proposed state-dependent policy that minimizes the sum of immediate readmission risks and one-step expected readmission risks when using the myopic policy and thus depends on the discharge process, readmission rates, and future SNF availability, with the optimal MDP policy.  We compare the two policies in three scenarios,  where each scenario considers increasing dependency between transfer decision and overall SNF availability. Our numerical study provides several useful insights for hospitals as well. For example, we find that when the availabilities of SNFs change quickly and when all SNFs have similar availability in the long term, the optimal policy obtained by the MDP model provides significant benefits over the two heuristic policies in terms of long-run average readmission rates.

\textbf{Organization of the paper.}  The remainder of this paper is organized as follows.  Section~\ref{sec:lit} summarizes the literature on readmission risk estimation and related operations research (OR) and management science (MS) studies.   Section~\ref{sec:model} describes our modeling approach.  We present our MDP formulation and analysis of the optimal policy.  Section~\ref{sec:paramest} presents how we estimate readmission rates for various SNFs based on real-life data from a large, tertiary teaching hospital and health system. Section~\ref{sec:numstudy} presents our numerical study.  We conclude in Section~\ref{sec:conclusion}. 

\section{Literature Review} \label{sec:lit}

This study is a synthesis of three main areas of the literature: the clinical literature on factors associated with SNF readmissions, the OR/MS literature on care transitions, and sequential decision-making under uncertainty.  Here, we summarize the studies most closely related to the current work.

\subsection{Estimating SNF readmissions} 

We use observational data at the patient-level from a large, tertiary teaching hospital to estimate readmission rates for SNFs using logistic regression.  There are several related studies that also analyze factors associated with SNF readmissions.
We classify this large body of literature into two broad categories.

The first consists of studies that analyze associations between specific factors (e.g., state-level policies, specific interventions, demographics, length of stay [LOS], rural vs. urban SNF) and SNF readmissions using a similar approach to ours (e.g., logistic regression).  In~\cite{grabowski2010medicaid}, for example, the authors analyze the impact of Medicaid bed-hold policies, which are state-level policies that pay nursing homes to reserve beds of acutely hospitalized Medicaid residents, on rehospitalization for Medicare post-acute SNF residents. Links between racial disparities and rehospitalization rates among non-Hispanic White and Black Medicare beneficiaries admitted to SNFs for postacute care are analyzed by~\cite{li2011racial}. We do not analyze these specific relationships, but focus instead on estimating readmission rates by patients with different care needs across different SNFs.  This is done to capture variability in readmission rates that can be used for improving how transfer decisions can be made in an optimization model. 

The second body of work focuses on  variability in SNF readmissions, which is relevant to our work since our optimization model utilizes variability in SNF readmissions  \citep{chen2018readmission,dombrowski2012factors,kim2019derivation,rahman2016skilled}. For example, \cite{dombrowski2012factors} examine risk factors for rehospitalization, but unlike our study, focus only on rehabilitation patients.   \cite{ogunneye2015association} assess associations between quality of care and 30-day risk-adjusted readmission rate only on patients with acute decompensated heart failure. \cite{fernandes2018thirty} focus on vascular surgery patients discharged to SNFs. Although our study considers all patient types using broader categories and uses risk-adjusted readmission rate estimates, our goal differs from the aforementioned studies in that we use these estimates to guide the transfer of the patients to SNFs in a sequential decision-making model that also incorporates stochastic models for the patient discharge process and capacity availability. 

\subsection{Care transitions using OR/MS tools and concepts} 

Focusing on performance analysis, \cite{mo2014modeling} considers patients transitioning from inpatient care setting to home or community based settings and patients transitioning between primary, secondary, and tertiary rehabilitation facilities. They use stochastic models and simulation to analyze the impact of LOS that is covered by insurance on reducing total care spending versus improving the quality of care to individual patients.  \cite{li2016capacity} consider transitions between long-term care services, such as nursing homes, home, and community-based services.  They also use stochastic modeling to capture transition dynamics between these services and formulate capacity planning needs as a newsvendor model to identify factors (e.g., cost for capacity shortage, transition rates between different services) for making capacity decisions.  \cite{mo2016simulation} consider transitions between inpatient acute, outpatient sub-acute, and general residential care for patients with traumatic brain injuries.  They use prediction models to estimate readmission rates and embed this in a simulation model to examine how coverage duration of publicly funded rehabilitation impacts sub-acute rehabilitation readmission and total rehabilitation spending.  Lastly, \cite{li2017evaluation,li2018threshold} consider transitions between upper level hospitals (e.g., comprehensive hospitals) to lower level hospitals (e.g., community based care) as a result of reverse referrals.  In~\cite{li2017evaluation}, queueing models are used to capture patient flow (e.g., waiting times, service, and transitions) and to evaluate performance based on profitability of each hospital and its willingness to participate in these referrals.  In~\cite{li2018threshold}, the authors propose and analyze the impact of reverse referral practices, defined as referring patients from a comprehensive hospital to community healthcare centres, which is a common practice in China, on system performance measures and profit using stochastic modeling techniques.

The aforementioned studies assess performance of specific care transitions assuming a given transition policy whereas we are interested in controlling transitions by dynamically varying how patients are routed from hospital to SNFs immediately after discharge in order to minimize readmission rates.  A related study that analyze discharge timing strategies include~\cite{chan2012optimizing} and~ \cite{crawford2014analyzing}. \cite{chan2012optimizing} consider discharge timing strategies for patients receiving care in an intensive care unit.  They use a sequential decision-making framework to analyze strategies that prioritize patients based on some measure of criticality assuming a model of readmission risk.  \cite{crawford2014analyzing} use simulation of patient pathways through an acute care hospital that comprises an emergency department (ED) and several inpatient units to analyze the impact of discharge timing on ED-related measures and readmissions.  By contrast, our study focuses on transfer decisions after discharge (and not on the timing of discharge).  We thus require a model of both the discharge process and of SNF availability following each transfer decision.  To our knowledge, no study in the OR/MS literature has considered optimizing transitions from hospitals to SNFs.  

\subsection{Stochastic optimization models} 

We propose and analyze an MDP model for allocating patients to SNFs following discharge.  Closely related stochastic optimization models are those that focus on routing randomly arriving jobs to multiple service facilities \citep{ephremides1980simple,hajek1984optimal,altman2000open,bhulai2009dynamic} and centralized dynamic stochastic matching models \citep{gurvich2015dynamic,ozkan2020dynamic}. Among these, the model proposed and analyzed by \cite{altman2000open} is closest to the one proposed and analyzed here.  \cite{altman2000open} study the assignment of jobs to a finite number of parallel, heterogeneous servers with no buffers.  The states of each of the servers are unknown and decisions are based on the number of time slots from the current decision epoch that jobs have been sent to the different servers. They use an MDP formulation to minimize the expected average cost.  

We also consider an arrival process of different job classes corresponding to the discharge process for each patient type, which is independent of the service capacity, and where each discharged patient can be assigned to one of a finite number of service facilities corresponding to each SNF.  In contrast to ~\cite{altman2000open} and traditional routing problems in the literature, immediate costs, corresponding to readmission rates, are incurred at each decision epoch that depend on both the patient class and SNF destination.  Additionally, we do not have full information about SNF capacity, instead we only observe whether a SNF is available or not at the time of discharge. Furthermore, SNF availability transitions according to a probability matrix that depends on the current availability in the SNF and transfer decision, regardless of whether a patient is assigned to that particular SNF.  This generalization comes at the cost of tractability.  To our knowledge, the proposed model has not been considered by the stochastic optimization literature before.  

\section{Modeling Framework} \label{sec:model}

\subsection{Rationale for a sequential decision making approach} \label{sec:rationale}

As alluded to, at the time of discharge, the case/social worker at the hospital works with the patient or their family to find a facility. This involves reaching out to possibly different facilities to learn whether there is capacity to admit the patient or not. The information shared with the case/social worker is binary – whether there is room for the patient or not. This decision is made for the current patient who needs to be admitted rather than several, or batch, admissions.  Now, assuming a centralized decision maker (i.e., the case/social worker) who is given estimates for readmission rates by patient type (See Section~\ref{sec:rates} below), which are our main proxy for quality of service and costs, patients could be transferred myopically: each time a patient is discharged, allocate the patient to the SNF for which the corresponding estimated readmission rate is the lowest among the available SNFs at the time of discharge. We refer to this policy as the \emph{myopic policy} (see Definition~\ref{defn:myopic} below). The patient discharge process, however, depends on the patient type and is independent of SNF availability.  Moreover, transfer decisions impact readmission rates and SNF availability. Therefore, it is not obvious that using a myopic policy is always optimal (see Examples~\ref{ex1} and~\ref{ex2} below).

Indeed, if a particular patient type is transferred to the ``best" available SNF (i.e., the one with the lowest readmission rate that is available), this SNF may not be available for a more complicated patient at the next discharge period. If it turns out that the latter SNF is also the best SNF for the more complicated patient, then it may be better to save that capacity for the more complicated patients later on.

To examine this tradeoff formally, the following sections propose and analyze an MDP formulation.

\subsection{MDP formulation} \label{sec:mdp}

Let $i$ denote patient type $i=0,1,\ldots,k$; $j$ denote SNF $j=0, 1,\ldots,\ell$; $s_j \in \{0,1\}$ the state of SNF $j=0,1,\ldots,\ell$, where  $1$ ($0$) represents the case that a SNF is (is not) available.  The state space is given by the set of $\ell+1$-tuples $\X:=\{(i,s_0,s_1,\ldots,s_\ell)| i \in\{0,1,\ldots,k\}, s_j \in \{0,1\}, j=0,1,\ldots,\ell\}$.  For patient type $i$, let $A(i) \subset \{0,1,\ldots,\ell\}$ denote the set of feasible SNFs that patient type $i$ can transition to.  For $x=(i,s_0,s_1,\ldots,s_{\ell})$, the set of feasible transfer options/actions is given by  $A(x)=\{j \in A_i | s_j = 1\}$.  Finally, let $r^j_i$ denote the readmission risk of patient type $i$ when the patient is transferred to SNF $j$.

We assume that at most one patient is discharged in successive time periods and that patient type $i$ is discharged in a given time period with probability $\lambda_i$ and that no patient is discharged in a time period with probability $\lambda_0=1-\sum_{i=1}^k \lambda_i$ with $0<\lambda_i<1$ for all $i=1,\ldots,k$.  We also assume that the patient arrival process is independent of the SNF availability process.  Lastly, for a fixed SNF $j \in \{0,1,\ldots,\ell\}$, the state of SNF $j$ representing whether the SNF is available or not transitions from $s_j \in \{0,1\}$ to $s'_j \in \{0,1\}$ under transfer decision $a$ according to a probability $p^{a,j}_{s_j,s'_j}$.  The transition probabilities capture SNF availability changes that happen throughout the day, which includes the effect of a decision taken to use or not use available SNFs.  In Remarks~\ref{rem:1} and~\ref{rem:2} below, we highlight how the same framework can be used to include more complicated transition probabilities that may accommodate more complicated history-dependence.

The decision-making scenario is as follows. At the beginning of each time period $t$, the decision-maker views the state corresponding to patient discharges.  If the state corresponding to patient discharges is $i=0$, then there is no transfer decision to be made and no costs (i.e.\ readmission rate) are incurred.  If, however, there is a discharged patient and  all SNFs are occupied (i.e., $s_j=0$ for all $j\in \{1,\ldots,\ell\}$, then the current patient is `lost' and a `loss' penalty of $K >> \max_{\{i,j\}}\{r^j_i\}$ is incurred.  This is done by including a SNF, denoted by $0$, that is always available (i.e. $p^{a,0}_{s_0,1}=1$ for all $a$ and $s_0$) and so that $r^0_0=0$ and $r^0_i=K$ for $i \neq 0$.  Otherwise, the decision-maker makes a transfer decision to an actual SNF and an immediate cost of $r^j_i$ is incurred when patient type $i$ is transferred to SNF $j$. In all cases, after system costs are incurred, each SNF availability transitions to a potentially new availability at the beginning of the next time period according to the given transition probabilities. We seek a policy that after each patient arrival, describes where to allocate the patient based on the current state and history of states and actions (i.e., a non-anticipating policy).  

The finite horizon, discounted expected cost for a non-anticipating policy $\pi$ is $v_{t,\alpha}^{\pi}(x) \equiv \E_{x} \sum_0^t \Big[ \alpha^s r^{A^{\pi}(s)}_{X^{\pi}(s)}\Big] ds$, where $A^{\pi}(s)$ and $X^{\pi}(s)$ represent, respectively, the action and patient state process at time $s$.  For $1 > \alpha > 0$, the infinite-horizon discounted expected cost under policy $\pi$ is $v_{\alpha}^{\pi}(x) \equiv \lim_{t \rightarrow \infty} v_{t,\alpha}^{\pi}(x)$.  The long-run average cost rate is $\rho^{\pi}(x) \equiv \limsup_{t \rightarrow \infty} v_{t,0}^{\pi}/t$.  Under either optimality criterion we seek a policy $\pi^\ast$ such that $w^{\pi^\ast}(x) = \inf_{\pi \in \Pi} w^{\pi}(x)$ where $\Pi$ is the set of all non-anticipating policies and $w = v_{\alpha}$ or $\rho$.

To analyze the structure of optimal policies and to approximate their values, we rely on the optimality equations.  Before we can use the optimality equations, we must first provide conditions that guarantee that the optimality equations under each criterion has a solution. The proofs can be found in Section~\ref{appendix:0} of the Appendix.

For any nonnegative function $v$ on $\X$ define the mapping $T$ such that
\begin{align*}
 T v(x) = & \min_{a \in A(x) } \{ r^a_i + \alpha \sum_{x'=(i',1,s'_1,\ldots,s'_{\ell}) \in \X} \lambda_{i'} \Pi_{j=1}^\ell p^{a,j}_{s_j,s'_j} v(x') \}
\end{align*}
for $x = (i,s_0,s_1,\ldots,s_{\ell})\in \X$. 

\begin{prop} \label{prop:dcoe}
Suppose $1>\alpha > 0$.  The following hold:
\begin{enumerate}
 \item The function $v_{\alpha}$ satisfies the discounted-cost optimality equations, i.e.,
 \begin{align*}
   v_{\alpha} = T v_{\alpha}.
 \end{align*}  
 \item There exists a stationary, deterministic policy $f_{\alpha}$ (depending on the discount factor $\alpha$) that attains the minimum in the right hand side of the optimality equations, and hence, $f_{\alpha}$ is discounted cost optimal.
\end{enumerate}
\end{prop}

\begin{prop} \label{prop:acoe}
Suppose $0 <p^{a,j}_{s,s'} < 1$ for all $a$,$j \geq 1$,$s$, and $s'$.  The following hold:
\begin{enumerate}
 \item There exists a constant $g$ and function $h$ on the state space such that $(g,h)$ satisfies the average cost optimality equations, \[g \vec{1} + h = T h,\] where $\vec{1}$ is the vector of ones, $g$ is the optimal average cost, and hence, is unique, and $h$, known as the relative value function, is unique up to additive constants.
 \item A deterministic stationary policy $f$ is average cost optimal if and only if $f$ satisfies the minimum in the average cost optimality equations.
\end{enumerate}
\end{prop}

\begin{rem}[Non-stationary formulations] \label{rem:1}
In practice, discharge rates and/or transition probabilities may depend on the discharge period (i.e., they are non-stationary).  We remark that the model above can easily be modified to a finite horizon model with non-stationary discharge rates and transition probabilities, or to an infinite horizon model with nonstationary, but periodic, discharge rate and transition probabilities.  For both variations, solutions to optimality equations can be shown (proofs omitted).  Further, Examples~\ref{ex1}-~\ref{ex2} hold, and natural extensions to Propositions~\ref{prop:1}-~\ref{prop:3} and Heuristic $r + pr$ (see Section~\ref{sec:numstudy}) remain valid for those variations, and as such, we focus on the infinite-horizon, stationary model described above.
\end{rem}

\begin{rem}[History-dependent process] \label{rem:2}
In practice, it may also be the case that SNF availability depends not only on the current state and allocation decision, but on previous availability and allocation decisions.  To this end, within the same framework, we can redefine the state process to be $\hat{X}^j(\ell)_t=(s^j_{t-\ell}, a^j_{t-\ell}, s^j_{t+1-\ell},\ldots,a^j_{t-1},s_t)$ for SNF $j \in \{0,\ldots,k\}$, and with $s_k$ and $a_k$ ($k=t-\ell,\ldots,t)$, respectively, denoting availability and allocation in time period $k$, and where $\ell$ is some prespecified lag, which can be estimated from data.  The caveat is a larger state space.
\end{rem}

\subsection{The Optimality of Myopic Policy}
In this section, we analyze the structure of the optimal policy.  We first provide conditions under which the myopic policy is optimal in Proposition~\ref{prop:1}. This result is useful for hospitals that consider implementing the  easy-to-use myopic policy instead of using an optimization model. We then provide conditions under which the optimal policy has a threshold-type structure in Proposition~\ref{prop:3}.  We first formally define the myopic policy, where for ease of exposition, we assume that $r^j_i \neq r^{j'}_i$ when $j \neq j'$ for all $i \geq 1$.  

\begin{defn} \label{defn:myopic}
A stationary deterministic policy $f$ is \textbf{myopic} if, for every state $x=(i,s_0,s_1,\ldots,s_{\ell})$, $f(x) = \argmin_{a \in \{1,\ldots,\ell\}|s_a \neq 0}\{r^a_i\}$.  
\end{defn}

One might conjecture that the myopic policy given by Definition ~\ref{defn:myopic} is always optimal, but as conjectured in the previous section, the following example confirms that not only is this not the case, but that the myopic policy may yield long-term average costs (i.e., readmission rates) that are significantly higher than those from the  optimal policy.
\begin{ex} \label{ex1}
Suppose we have the following inputs for the model: $k=\ell=2$; $A(0) = \{0\}$, $A(1) = A(2) = \{1,2\}$; $\lambda_1=\lambda_2=0.4$; transition probabilities associated with each transfer decision and SNF are
\begin{align*}
& \mathbf{P}^{0,0}=
\begin{bmatrix}
0&1\\
0&1
\end{bmatrix}, \mathbf{P}^{1,0}= \begin{bmatrix}
0&1\\
0&1
\end{bmatrix}, \mathbf{P}^{2,0}=\begin{bmatrix}
0&1\\
0&1
\end{bmatrix},\\
& \mathbf{P}^{0,1} = \begin{bmatrix}
0&1\\
0&1
\end{bmatrix}, \mathbf{P}^{1,1} = \begin{bmatrix}
0.49&0.51\\
0.99&0.01
\end{bmatrix}, \mathbf{P}^{2,1} = \begin{bmatrix}
0.05&0.95\\
0.51&0.49
\end{bmatrix},\\
& \mathbf{P}^{0,2} = \begin{bmatrix}
0&1\\
0&1
\end{bmatrix}, \mathbf{P}^{1,2} = \begin{bmatrix}
0.01&0.99\\
0.05&0.95
\end{bmatrix}, \mathbf{P}^{2,2} = \begin{bmatrix}
0.5&0.5\\
0.95&0.05
\end{bmatrix};
\end{align*}
costs given by
\begin{align*}
& \mathbf{r}^{0}=
\begin{bmatrix}
0\\
10\\
10
\end{bmatrix},\mathbf{r}^{1}= \begin{bmatrix}
n/a\\
0.5\\
1.3
\end{bmatrix},\mathbf{r}^{2}= \begin{bmatrix}
n/a\\
0.55\\
1.2
\end{bmatrix}.
\end{align*}
For this example, the optimal average cost is approximately $2.9$ and the cost of the myopic policy is approximately $6.34$: the myopic leads to an average cost that is 54\% higher than the optimal average cost.  The optimal policy differs from the myopic policy in state $(2,1,1)$ corresponding to discharging patient type 2 when both SNFs are available (see Table~\ref{tab:ex1} in the Appendix~\ref{appendix:0.5}).  In this case, the optimal policy chooses to send patient type $2$ to SNF $1$ and the myopic policy sends patient type $2$ to SNF $2$.  Sending patient type $2$ to SNF $2$, as the myopic policy does, will make SNF $2$ unavailable at the next stage with probability $0.95$ whereas sending patient type $2$ to SNF $1$ will make SNF $2$ unavailable with probability $0.51$.  At the same time, the difference between $r^1_2=1.3$ and $r^2_2=1.2$ is relatively small.  As a result, it is better to send patient type 2 to SNF 1.
\end{ex}

A natural question then becomes, are there conditions under which the myopic policy is optimal?  We begin our analysis with the following intuitive result, which will form the basis of our numerical study (Section~\ref{sec:numstudy}).  The result says that if the transfer decision has no impact on subsequent SNF availability, then it is optimal to use the myopic policy as given by  Definition~\ref{defn:myopic}.  The proof can be found in Section~\ref{appendix:1} of the Appendix.
\begin{prop} \label{prop:1}
 Under both the $\alpha-$discounted cost and average cost criteria, if $p^{a,j}_{s,s'}=p^j_{s,s'}$ for all $a,j,s,s'$, then it is optimal to allocate patients myopically.
\end{prop}

One may also conjecture that if SNF availability only depends on the transfer decision and not the SNF (i.e., $p^{a,j}_{s_j,s'_j}=p^a_{s,s'}$ for all $a$, $j \geq 1$, $s_j$, and $s'_j$), then the myopic policy is once again optimal.  The following example illustrates that this may not be the case.  
\begin{ex} \label{ex2}
Suppose we have the following inputs for the model: $k=\ell=2$; $A(0) = \{0\}$, $A(1) = \{1,2\}$, $A(2) = \{1,2\}$; $\lambda_1=\lambda_2=0.4$; transition probabilities associated with each transfer decision and SNF are given by
\begin{align*}
& \mathbf{P}^{0,0}=
\begin{bmatrix}
0&1\\
0&1
\end{bmatrix}, \mathbf{P}^{1,0}= \begin{bmatrix}
0&1\\
0&1
\end{bmatrix}, \mathbf{P}^{2,0}=\begin{bmatrix}
0&1\\
0&1
\end{bmatrix},\\
& \mathbf{P}^{0,1} = \begin{bmatrix}
0&1\\
0&1
\end{bmatrix}, \mathbf{P}^{1,1} = \begin{bmatrix}
0.02&0.98\\
0.82&0.18
\end{bmatrix}, \mathbf{P}^{2,1} = \begin{bmatrix}
0.3&0.7\\
0.08&0.92
\end{bmatrix},\\
& \mathbf{P}^{0,2} = \begin{bmatrix}
0&1\\
0&1
\end{bmatrix}, \mathbf{P}^{1,2} = \begin{bmatrix}
0.02&0.98\\
0.82&0.18
\end{bmatrix}, \mathbf{P}^{2,2} = \begin{bmatrix}
0.3&0.7\\
0.08&0.92
\end{bmatrix};
\end{align*} costs given by
\begin{align*}
& \mathbf{r}^{0}=
\begin{bmatrix}
0\\
10\\
10
\end{bmatrix},\mathbf{r}^{1}= \begin{bmatrix}
n/a\\
0.5\\
1.3
\end{bmatrix},\mathbf{r}^{2}= \begin{bmatrix}
n/a\\
0.55\\
1.2
\end{bmatrix}
\end{align*}
where $\mathbf{r}^a$ denotes the readmission rates corresponding to transfer decision $a$.  For this example, the average costs for the optimal and myopic policy are $1.3$ and $14.7$, respectively; the average cost of the myopic policy is 91\% higher than the optimal average cost.  The optimal policy differs from the myopic policy only in state $(1,1,1)$ corresponding to discharging patient type $1$ when both SNFs are available (see Table~\ref{tab:ex2} of the Appendix \ref{appendix:0.5}).  In this case, the optimal policy chooses to send patient type $1$ to SNF $2$ and the myopic policy sends patient type $1$ to SNF $1$.  Sending patient type $1$ to SNF $1$, as the myopic policy does, will make SNF $1$ unavailable at the next stage with probability $0.82$ whereas sending patient type $1$ to SNF $2$ will make SNF $2$ unavailable with probability $0.08$.  At the same time, the difference between $r^1_1=0.5$ and $r^2_1=0.55$ is relatively small.  As a result, it is better to send patient type $1$ to SNF $2$.
\end{ex}

We make the following remarks regarding Proposition~\ref{prop:1}.  First, Proposition~\ref{prop:1} suggests that when SNF availability is robust to transfer decisions in that they do not severely impact their availability, it may be advantageous to use a myopic policy.  Second, we have assumed that the probability of transitioning from availability $x=(1,s_1,\ldots,s_{\ell})$ to availability $x'= (1,s'_1,\ldots,s'_{\ell})$ has the multiplicative form $ \Pi_{j=1}^\ell p^{a,j}_{s_j,s'_j}$ but Proposition~\ref{prop:1} holds for any transition probability governing SNF availability so long as it does not depend on the transfer decision $a$.  Third, once transfer decisions impact SNF availability, myopic policy may no longer be optimal, as illustrated in Examples~\ref{ex1} and ~\ref{ex2}.  Further, Proposition~\ref{prop:1} and Examples~\ref{ex1}-~\ref{ex2} suggest a tradeoff between using a myopic policy when transfer decisions have little to no impact on SNF availability and when they do have a significant impact on SNF availability.  This tradeoff is examined in detail in Section~\ref{sec:numstudy} where we also propose a heuristic, which we label $r + p r$, that attempts to circumvent the short-sightedness of the myopic policy by making allocation decisions that minimize the sum of immediate readmission risk and one-step expected readmission risks.


Unfortunately, additional conditions that guarantee that the optimality of the myopic policy and other index-like policies remain elusive, even under seemingly more tractable assumptions as highlighted in Examples~\ref{ex1}-~\ref{ex2}.  For example, when the myopic policy yields exactly the same policy for each patient type and if it is always more likely that each SNF is available at the next decision epoch in general (i.e., $p^{a,j}_{s_j,1} \geq p^{a,j}_{s_j,0}$ for all $a$, $j$, and $s_j$) \emph{and} when the myopic policy is used, the myopic policy was still not optimal.  A natural question is then whether additional properties of the optimal policy can be discerned.  We conclude this section with one such additional property.  It is a \emph{second order} result that provides a set of conditions on the transition probabilities under which a threshold-type structure is optimal.  The result shows that if the rate of change in SNF availability under two different transfer decisions is nondecreasing as more SNFs become available, then, if it is optimal to send a patient type to a specific SNF, then it is also optimal to send the same patient type to the same SNF (if available) when overall SNF availability decreases.  The proof can be found in Section~\ref{appendix:1} of the Appendix.
\begin{prop} \label{prop:3}
Under both the $\alpha-$discounted cost and average cost criteria, if it is optimal to allocate patient type $i$ to SNF $a^*$ in state $x=(i,s_1,s_2,\ldots,s_{\ell})$ and
\begin{align}
 \Pi_{j=1}^\ell p^{a^*,j}_{s_j',s''_j} - \Pi_{j=1}^\ell p^{a',j}_{s_j',s''_j} \leq \Pi_{j=1}^\ell p^{a^*,j}_{s_j,s''_j} - \Pi_{j=1}^\ell p^{a',j}_{s_j,s''_j} \label{ass:prop3}
\end{align}
holds for $x'=(i,s_1',s_2',\ldots,s_n')$ with $s_{a^*}=s'_{a^*}$ and $s_j \geq s_j'$, all feasible $a' \in A(x')$, and all $x''=(i'',s''_1,\ldots,s''_{\ell}) \in \X$, then it is also optimal to allocate patient type $i$ to SNF $a^*$ in state $x'$. 
\end{prop}

\section{Logistic Regression Model of Readmission Rates by SNF and Patient Type} \label{sec:paramest}
We test the performance of transfer policies obtained by our MDP model (Section~\ref{sec:mdp}) using real-life data from our partner hospital. For this purpose, we develop a logistic regression model to estimate the readmission rates $r^j_i$, which provide a key input to our MDP model. We note that the same data has been used to estimate readmission rates in \citet{oruongo2020skilled}. We first describe our data sources and then present our logistic regression model. 

\subsection{Data sources}
We estimated the readmission rates using CMS claims data and identified all discharges for patients assigned to our partner health system in the Medicare Shared Savings Program. For each patient discharge, we recovered \ul{demographic information} (age, sex);  \ul{the LAC index} (a readmission risk score derived from \ul{L}ength of stay, \ul{A}cuity of admission, and \ul{C}omorbidity which ranges from 1 to 19 \citep{van2010derivation}); \ul{CMS Hierarchical Condition Category} (HCC; a risk-adjusted score designed to estimate future health care costs for patients \cite{pope2000diagnostic}); and \ul{number of Elixhauser comorbidities} \citep{elixhauser1998comorbidity}.  Additional binary variables are recovered on prior medical history specifying whether each patient had prior hospitalization, congestive heart failure (CHF), chronic obstructive pulmonary disease (COPD), uncomplicated diabetes (diabetes without other chronic complications), complicated diabetes (diabetes with other chronic complications), prior stroke, obesity, depression, dementia, hypertension, hypothyroidism, metastic solid tumor, valvular disease, peripheral vascular disease (PVD), and deficiency anemias.  We restricted attention to discharges from an inpatient stay to a SNF. During an inpatient stay, a patient was grouped based on their DRG into one of four disjoint types: \ul{uncomplicated medical}, \ul{complicated medical}, \ul{complicated surgery}, and \ul{joint surgery}. These patient types were derived by grouping inpatient DRGs based on the joint likelihood of discharge to a SNF and readmission. Complicated medical and surgery groups, for example, had major DRG clinical complications designation during their inpatient stay, required ventilation support, or needed skin graft as part of care received.  We focused on the top five facilities, labeled \textbf{A}, \textbf{B}, \textbf{C}, \textbf{D}, and \textbf{E}, that had the largest number of discharged patients and that are relatively close to each other so as to minimize the impact of distance on patient preferences for specific SNFs, as patients may choose the closest SNF to their homes or the residency of their caregivers. The discharge dates in our sample ranged from January 2012 to April 2018.  In total, there were 4,116 discharges to these five SNFs during the study period.

\subsection{Summary statistics}

Summary statistics were recovered (Tables~\ref{tbl1:snf} and~\ref{tbl1:drg}). For continuous variables, we performed one-way ANOVA to check for differences in means between SNFs or patient types.  For binary variables, we used $\chi^2$ tests to check for differences in proportions.  Significance was considered an $\alpha$ level of 0.05. 

Each SNF treats significantly-different groups of patients, which may impact readmission rates (Table~\ref{tbl1:snf}). For example, patients in SNF C have an average age of 83.1 years compared to 78.7 years in SNF A. Variables that did not differ significantly by SNF were prior hospitalization, hypothyroidism, metastatic solid tumor, or having at least three Elixhauser chronic conditions, prior stroke and proportion Caucasian.  

\begin{table}[ht!]
\footnotesize
\centering
\begin{tabular*}{\textwidth}{l@{\extracolsep{\fill} }ccccccr}
\toprule
& & \multicolumn{5}{c}{SNF}  \\
\cline{3-7}
& Overall & A & B & C & D & E & \multicolumn{1}{c}{\emph{p}} \\ 
\midrule
Mean Age (SD) & 80.3 (9.5) & 78.7 (9.8) & 80.2 (9.3) & 83.1 (8.6) & 80.2 (9.2) & 79.7 (9.8) & $<$0.01 \\ 
Mean LAC (SD) & 5.1 (2.0) & 5.1 (1.9) & 5.4 (2.0) & 5.0 (1.9) & 4.8 (1.9) & 5.2 (2.0) & $<$0.01 \\ 
Mean HCC (SD) & 2.2 (1.1) & 2.3 (1.2) & 2.4 (1.1) & 2.1 (1.1) & 2.0 (1.0) & 2.3 (1.2) & $<$0.01 \\ 
Mean Elix Count (SD) & 4.1 (2.9) & 4.1 (2.9) & 4.4 (3.2) & 4.0 (2.9) & 3.8 (2.6) & 4.2 (2.9) & $<$0.01 \\ 
Caucasian (\%) & 96 & 95 & 96 & 98 & 97 & 96 & 0.06 \\ 
Female (\%) & 69 & 65 & 67 & 70 & 75 & 69 & $<$0.01 \\ 
\midrule
\multicolumn{2}{l}{Prior medical history (\%)} \\
\enskip Prior hospitalization & 59 & 59 & 61 & 59 & 57 & 61 & 0.45 \\ 
\enskip CHF & 31 & 31 & 35 & 29 & 27 & 34 & 0.02 \\ 
\enskip COPD & 29 & 33 & 32 & 24 & 27 & 27 & $<$0.01 \\ 
\enskip Uncomplicated diabetes & 7 & 7 & 7 & 5 & 10 & 8 & 0.01 \\ 
\enskip Complicated diabetes & 14 & 12 & 21 & 11 & 12 & 14 & $<$0.01 \\ 
\enskip Prior stroke & 10 & 9 & 12 & 10 & 10 & 12 & 0.08 \\ 
\enskip Obesity & 13 & 12 & 15 & 10 & 15 & 16 & $<$0.01 \\ 
\enskip Depression & 26 & 30 & 25 & 24 & 24 & 24 & 0.01 \\ 
\enskip Dementia & 15 & 16 & 14 & 17 & 11 & 16 & 0.01 \\ 
\enskip Hypertension & 62 & 57 & 63 & 62 & 66 & 65 & $<$0.01 \\ 
\enskip Hypothyroidism & 25 & 26 & 25 & 27 & 25 & 23 & 0.60 \\ 
\enskip Metastatic solid tumor & 8 & 8 & 6 & 8 & 8 & 7 & 0.52 \\ 
\enskip Valvular disease & 16 & 13 & 15 & 19 & 15 & 18 & $<$0.01 \\ \enskip peripheral vascular disease & 25 & 26 & 25 & 29 & 20 & 25 & $<$0.01 \\ 
\enskip Deficiency anemias & 29 & 28 & 38 & 25 & 24 & 30 & $<$0.01 \\ 
\enskip 3+ Elix conditions & 65 & 65 & 66 & 63 & 65 & 67 & 0.50 \\ 
		Readmission rate & 16 & 16 & 19 & 16 & 11 & 17 & $<$0.01 \\
\bottomrule
\end{tabular*}
\caption[Discharge characteristics by SNF]{Patient characteristics of dataset by SNF.}
\label{tbl1:snf}
\end{table}

Patient types also differed significantly in most patient variables (Table~\ref{tbl1:drg}). For example, medical patient types (i.e., uncomplicated medical and complicated medical)  were older on average than the surgical patients (i.e., complicated surgery and joint surgery). Medical patient types also had higher prevalences of prior stroke, prior hospitalizations, and at least 3 Elixhauser chronic conditions  when compared to surgical patients. Additionally, female patients were most prevalent in the joint surgery group, with more than three out of four joint surgery discharges being female. Most notably, complicated patient types had the highest unadjusted readmission rates. The joint surgery patient type had unadjusted readmission rate that were less than half of the next highest unadjusted readmission rate. We did not observe significant differences across patient types with either diabetes without complications or obesity.

\begin{table}[htb!]
\footnotesize
\centering
\begin{tabular}{l l l l l l r}
  \toprule
  && \multicolumn{4}{c}{Patient type} \\
  \cmidrule(lr){3-6}
  & Overall & UM & CM & CS & JS & \emph{p} \\ 
  \midrule
		Mean Age (SD) & 80.3 (9.5) & 82.8 (8.9) & 82.8 (9.2) & 78.9 (9.5) & 75.7 (8.7) & $<$0.01 \\ 
		Mean LAC (SD) & 5.1 (2.0) & 5.2 (1.8) & 6.1 (1.8) & 6.0 (2.0) & 3.7 (1.4) & $<$0.01 \\ 
		Mean HCC (SD) & 2.2 (1.1) & 2.0 (0.8) & 3.0 (1.2) & 3.0 (1.4) & 1.5 (0.6) & $<$0.01 \\ 
		Mean Elix Count (SD) & 4.1 (2.9) & 4.6 (2.8) & 5.1 (2.9) & 4.3 (3.0) & 2.6 (2.3) & $<$0.01 \\ 
		Caucasian (\%) & 96 & 97 & 95 & 95 & 97 & 0.02 \\ 
		Female (\%) & 69 & 69 & 63 & 62 & 76 & $<$0.01 \\ 
		\midrule
		\multicolumn{2}{l}{Prior medical history (\%)} \\
		\enskip Prior hospitalization & 59 & 68 & 73 & 60 & 37 & $<$0.01 \\ 
		\enskip CHF & 31 & 35 & 48 & 31 & 11 & $<$0.01 \\ 
		\enskip COPD & 29 & 28 & 38 & 32 & 20 & $<$0.01 \\ 
		\enskip Uncomplicated diabetes & 7 & 7 & 7 & 9 & 7 & 0.72 \\ 
		\enskip Complicated diabetes & 14 & 15 & 20 & 17 & 7 & $<$0.01 \\ 
		\enskip Prior stroke & 10 & 13 & 14 & 9 & 4 & $<$0.01 \\ 
		\enskip Obesity & 13 & 13 & 12 & 15 & 14 & 0.27 \\ 
		\enskip Depression & 26 & 29 & 28 & 26 & 20 & $<$0.01 \\ 
		\enskip Dementia & 15 & 22 & 19 & 13 & 6 & $<$0.01 \\ 
		\enskip Hypertension & 62 & 65 & 64 & 62 & 56 & $<$0.01 \\
		\enskip Hypothyroidism & 25 & 28 & 27 & 25 & 20 & $<$0.01 \\ 
		\enskip Metastatic solid tumor & 8 & 10 & 7 & 9 & 6 & 0.01 \\ 
		\enskip Valvular disease & 16 & 17 & 22 & 19 & 8 & $<$0.01 \\ 
		\enskip PVD & 25 & 30 & 34 & 23 & 13 & $<$0.01 \\ 
		\enskip Deficiency anemias & 29 & 32 & 38 & 34 & 14 & $<$0.01 \\ 
		\enskip 3+ Elix conditions & 65 & 74 & 78 & 66 & 42 & $<$0.01 \\ 
		Readmission rate & 16 & 16 & 23 & 20 & 7 & $<$0.01 \\ 
\bottomrule
\end{tabular}
\caption[Discharge characteristics by patient type]{Patient characteristics of dataset by patient type. \\ \emph{Note.} UM = Uncomplicated Medical, CM = Complicated Medical, CS = Complicated Surgery; JS = Joint Surgery.}
\label{tbl1:drg}
\end{table}

\subsection{Readmission rates} \label{sec:rates}

We developed a logistic regression model and estimated readmission rates for different patient types using data on discharges from the hospital to the SNFs.  This work showed that differences in readmission rates exist between patient types and SNFs, providing supporting evidence that readmission rates could be reduced using our optimization approach. Estimated readmission rates are used as an input to our MDP model in our numerical study (Section~\ref{sec:numstudy}).

We analyzed unadjusted readmission rates and number of discharges ($n$) by patient type (Tables~\ref{tbl:drg_rates}) and by SNF (Table~\ref{tbl:snf_rates}). These tables again highlight that there exist significant differences in readmission rates among patient types. For example, we found that patients discharged after a complicated surgery had higher unadjusted readmission rates compared to those discharged after uncomplicated medical care (Table~\ref{tbl:drg_rates}). Similarly, readmission rates for patients who had joint surgery had significantly lower unadjusted readmission rates compared to other patient types (Table~\ref{tbl:drg_rates}).  

\begin{table}[h]
\footnotesize
\centering
\begin{tabular}{ccccccccc}
\toprule
& \multicolumn{2}{c}{$\begin{matrix} \rm{Uncomplicated} \\ \rm{Medical}\end{matrix}$} & \multicolumn{2}{c}{$\begin{matrix} \rm{Complicated} \\ \rm{Medical}\end{matrix}$} & \multicolumn{2}{c}{$\begin{matrix} \rm{Complicated} \\ \rm{Surgery}\end{matrix}$} & \multicolumn{2}{c}{$\begin{matrix} \rm{Joint} \\ \rm{Surgery}\end{matrix}$} \\
Year & \% &  n & \% &  n  & \% &  n  & \% &  n  \\ 
\midrule
  2012 & 0.13 & 121 & 0.14 & 73 & 0.23 & 56 & 0.06 & 223 \\ 
  2013 & 0.15 & 169 & 0.17 & 132 & 0.23 & 77 & 0.07 & 221 \\ 
  2014 & 0.15 & 197 & 0.25 & 159 & 0.25 & 80 & 0.03 & 202 \\ 
  2015 & 0.16 & 260 & 0.24 & 222 & 0.20 & 87 & 0.07 & 210 \\ 
  2016 & 0.18 & 293 & 0.26 & 263 & 0.14 & 108 & 0.09 & 169 \\ 
  2017 & 0.17 & 234 & 0.19 & 208 & 0.18 & 83 & 0.11 & 143 \\ 
  2018 & 0.15 & 39 & 0.34 & 50 & 0.33 & 18 & 0.05 & 19 \\ 
\bottomrule
\end{tabular}
\caption{Unadjusted readmission rates by patient type and number of discharges ($n$) over study period January 2012 to April 2018}
\label{tbl:drg_rates}
\end{table}

\begin{table}[h]
\footnotesize
\centering
\begin{tabular}{c cc cc cc cc cc}
\toprule
& \multicolumn{2}{c}{A} & \multicolumn{2}{c}{B} & \multicolumn{2}{c}{C} & \multicolumn{2}{c}{D} & \multicolumn{2}{c}{E} \\
Year & \% &  n & \% &  n  & \% &  n  & \% &  n & \% &  n \\ 
	\midrule
    2012 & 0.13 & 118 & 0.14 & 91 & 0.09 & 81 & 0.06 & 72 & 0.13 & 111 \\ 
    2013 & 0.11 & 161 & 0.25 & 140 & 0.12 & 128 & 0.09 & 81 & 0.09 & 89 \\ 
    2014 & 0.16 & 189 & 0.19 & 124 & 0.13 & 122 & 0.09 & 97 & 0.15 & 106 \\ 
    2015 & 0.17 & 223 & 0.18 & 144 & 0.18 & 174 & 0.10 & 107 & 0.16 & 131 \\ 
    2016 & 0.17 & 227 & 0.14 & 130 & 0.21 & 198 & 0.16 & 115 & 0.22 & 163 \\ 
    2017 & 0.17 & 182 & 0.21 & 109 & 0.15 & 144 & 0.11 & 100 & 0.18 & 133 \\ 
    2018 & 0.19 & 37 & 0.23 & 13 & 0.25 & 24 & 0.16 & 25 & 0.37 & 27 \\
	\bottomrule
\end{tabular}
\caption{Unadjusted SNF readmission rates and number of discharges ($n$) over study period January 2012 to April 2018}
\label{tbl:snf_rates}
\end{table}

We used logistic regression to estimate readmission rates by SNF and patient type, adjusting for discharge year, HCC, whether or not it was the patient's first hospitalization, LAC, CHF, stroke, dementia, and 3+ elixihauser chronic conditions. These latter variables for which we adjusted estimates were chosen using step-wise variable selection.  Confidence intervals for estimated readmission rates were recovered via bootstrapping.  Our preliminary results from the logistic regression model show that there exist differences between SNFs in caring for different types of patients (Table~\ref{tbltoss}).
\begin{table}[!htbp]
\footnotesize
\centering
\begin{tabular}{cccccc}
  \toprule
SNF & Overall & $\begin{matrix} \rm{Uncomplicated} \\ \rm{Medical}\end{matrix}$ & $\begin{matrix} \rm{Joint} \\ \rm{Surgery}\end{matrix}$ & $\begin{matrix} \rm{Complicated} \\ \rm{Medical}\end{matrix}$ & $\begin{matrix} \rm{Complicated} \\ \rm{Surgery}\end{matrix}$ \\ 
\midrule
A & 15.7 (14.7, 16.9) & 14.3 (13.3, 15.4) &  9.5 ( 8.6, 10.4) & 19.1 (18.0, 20.4) & 19.2 (18.0, 20.4) \\ 
B & 14.5 (13.4, 15.6) & 16.4 (15.3, 17.5) & 12.8 (11.8, 13.9) & 20.1 (18.8, 21.3) & 20.4 (19.0, 21.6) \\ 
C & 14.7 (13.7, 15.8) & 15.6 (14.6, 16.7) & 12.4 (11.4, 13.5) & 20.6 (19.3, 21.8) & 20.2 (19.0, 21.5) \\ 
D & 15.7 (14.5, 16.9) &  9.1 ( 8.2, 10.0) &  8.7 ( 7.9,  9.5) & 11.2 (10.3, 12.1) & 19.6 (18.5, 20.8) \\ 
E & 17.2 (16.0, 18.4) & 20.6 (19.3, 21.8) &  5.8 ( 5.1,  6.5) & 19.0 (17.8, 20.1) & 13.4 (12.4, 14.4) \\ 
\bottomrule
\end{tabular}
\caption{Readmission rates (in \%) for different SNFs by patient type (95\% confidence intervals)}
\label{tbltoss}
\end{table}

\section{Numerical Study} \label{sec:numstudy}

\subsection{Rationale}


We provided conditions under which a simple policy, the myopic policy, is optimal (Proposition~\ref{prop:1}) and under which the optimal policy has a particular threshold structure (Proposition~\ref{prop:3}). Yet, since these conditions may not hold in practice, we noted that when these conditions do not hold, the myopic policy can perform poorly and the optimal policy can be difficult to characterize (Examples~\ref{ex1}--\ref{ex2}). There is a critical need to address these issues, given that it is not easy to compute the optimal policy directly due to the curse of dimensionality associated with the MDP formulation, i.e., the dimension of the state-space is $|\X| = (k+1) \cdot 2^{\ell+1}$ and, at worst, the dimension of the action-space is $\A:=\cup_{x \in \X} A(x) =2^{\ell+1}$. 

Fortunately, our analysis of the structure of the optimal policy also provides key insights into the trade-off between simplicity and performance. That is, a simple myopic policy is optimal provided transfer decisions have no impact on SNF availability, but can perform poorly otherwise. Thus, our goals of this numerical study were to (i) compare the performance of the myopic policy against the optimal policy for a sequence of scenarios in which we increase the degree to which transfer decisions impact SNF availability; and (ii) investigate whether a new heuristic policy, which is an extension of the myopic policy, performs well against the optimal policy in these same scenarios. In short, we seek to determine when the myopic policy may be beneficial to use and, otherwise, when another heuristic policy, which is still relatively easy to compute, may provide a better alternative.  All of our numerical experiments use data from a large tertiary teaching hospital and nearby SNFs to ensure our parameters accurately reflect actual settings.


\subsection{Overview}

\begin{figure}
    \centering
    \includegraphics[width=0.9\textwidth]{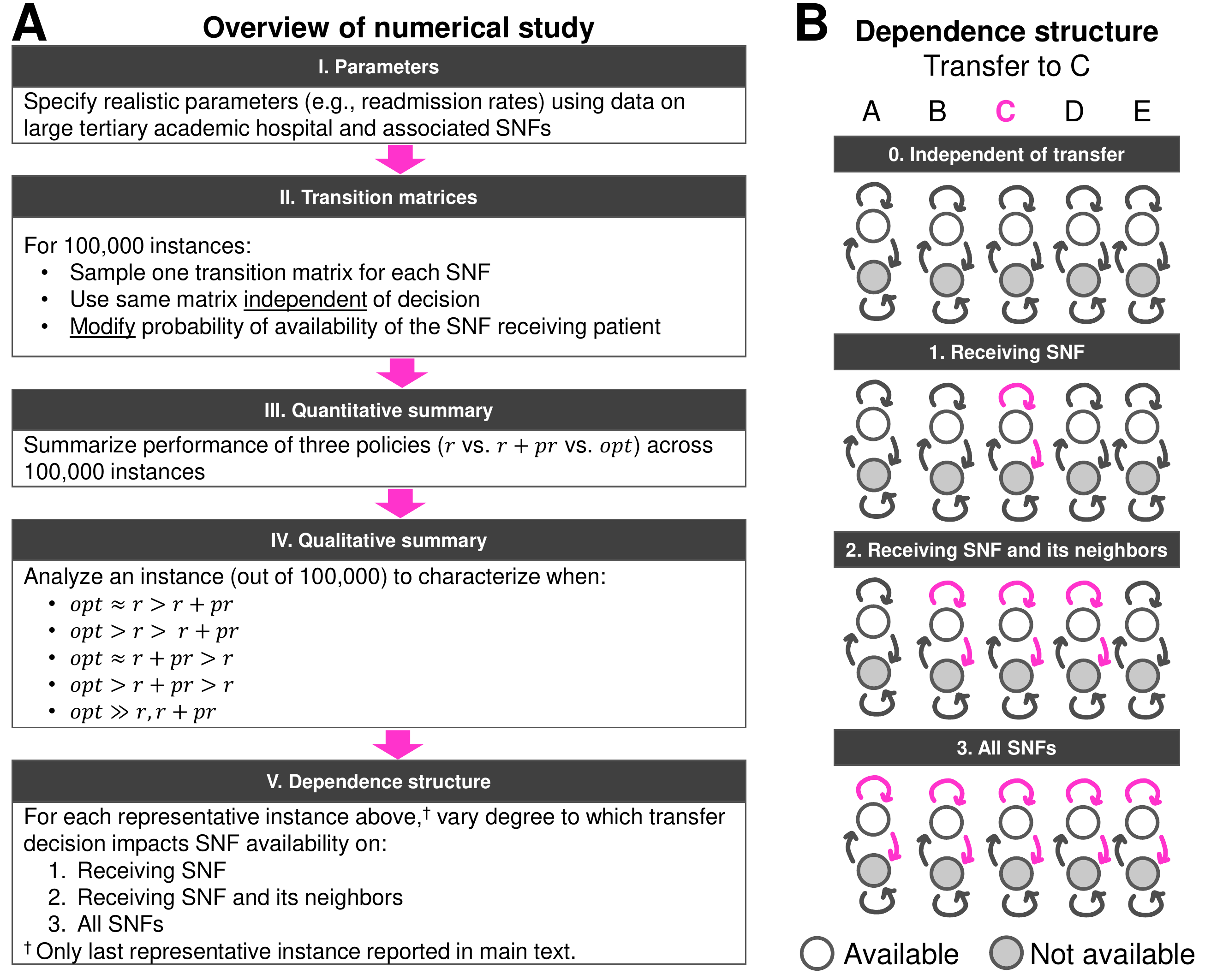}
    \caption{\textbf{A}) The numerical study seeks to evaluate, in five steps, the relative performance of three policies (opt: optimal; $r$: myopic; and $r+pr$: a state-dependent priority rule). \textbf{B}) On the fifth step, policies are examined in three scenarios reflecting difference dependencies between the transfer decision and SNF availability. The probability of SNF availability changes only for the receiving SNF in Scenario 1; the receiving SNF and its neighbors in Scenario 2; and all SNFs in Scenario 3.  }
    \label{fig:overview}
\end{figure}

Our numerical study is divided into five steps that are summarized in Figure~\ref{fig:overview}A. First, we use real-life data as described in Section~\ref{sec:paramest} to estimate the major input parameters for the simulation, such as number of patient types and readmission rates. Second, we describe how we constructed (random) transition probability matrices that allow us to control the degree to which a transfer decision impacts SNF availability. Third, these matrices are randomly sampled in order to quantitatively investigate the relative performance of the three policies (myopic, proposed heuristic, and optimal; as described below) across a wide range of instances. Performance is measured with respect to long-run average readmission rates. Fourth, we then narrow our attention on five different instances to qualitatively investigate when each heuristic performs well and when the optimal policy significantly outperforms the two heuristics. 

The last step is to evaluate how the dependency of the SNF availability on transfer decisions impacts the performance of the two heuristics.  We consider three scenarios to compare policies under different types of dependency structures (Figure~\ref{fig:overview}A). Scenario 1 considers transfer decisions that impact the availability of only one SNF: the SNF receiving the patient.  Scenario 2 considers transfer decisions that impact the availability of the receiving SNF as well as the availability of the `neighboring' SNFs (whenever possible). Scenario 3 considers transfer decisions that impact all SNFs regardless of whether or not they receive a patient. In each instance and scenario, we calculate optimal policies using policy iteration, and compare average costs of the optimal policies with those of two proposed heuristic policies. Note that these three scenarios are evaluated for the experiments in both the third and fourth steps of our numerical study.  

\subsection{Heuristic transfer policies}

Two types of heuristic policies are compared. The first heuristic is the myopic policy. It is referred to as ``$r$" (Definition~\ref{defn:myopic}) and is based only on immediate costs. The second heuristic is a state-dependent priority rule that also depends on the discharge process, immediate and one-step costs, and transition probabilities for SNF availability. It is referred to as ``$r+pr$". These policies are defined as follows:

\medskip

\begin{tabular}{l p{4.8in}}
$r$		        & Prioritize patient types in increasing order of $r_i^a$. \\
$r + pr$ 		& Prioritize patient types in state $x=(i,s_1,\ldots,s_{\ell})$ in increasing order of                   $$r^a_i + \sum_{x'=(i',s'_1,\ldots,s'_n)} \lambda_{i'} \Pi_{j=1}^\ell p^{a,j}_{s_j,s'_j} \min_{a' \in A(x')}\{r^{a'}_{i'}\}.$$
\end{tabular}

\medskip


Policy $r+pr$ attempts to overcome a drawback of policy $r$ (myopic) observed in Examples~\ref{ex1}--\ref{ex2}, Namely, policy $r$ performs poorly when transferring a patient type to a SNF causes that SNF becoming unavailable in the next decision epoch \emph{and} an alternative transfer yields similar immediate costs. By contrast, policy $r+pr$ incorporates information on both immediate \emph{and} future costs, which depend on the discharge process and future SNF availability. It chooses the action that minimizes the sum of the immediate readmission rate and one-step future (expected) readmission rates relative to the patient class and current transfer decision when using the myopic policy at the next stage.  Implementing policy $r+ p r$ requires additional arithmetic operations compared to implementing the myopic heuristic $r$ and decision-makers may have to balance improved performance from implementing the more complicated policy $r+pr$ (via lower long-run average readmission rate) with heuristic complexity.  To this end, we end this section with the following result on the complexity of the two policies, with the proof provided in Section~\ref{appendix:1.5} of the Appendix.
\begin{prop}~\label{prop:complexity}
The algorithm complexity for policy $r$ is $\mathcal{O}\left(M\right)$ and the algorithm complexity for heuristic $r+p r$ is $\mathcal{O}\left(M\cdot 2^l\right)$, where $M = \sum_{x\in\X} |A(x)|$ and $\mathcal{O}(\cdot)$ denotes the Bachmann–Landau notation for Big O \citep{knuth1968art}.
\end{prop}

\subsection{Input parameters}

As described in Section~\ref{sec:paramest}, we consider four different patient types: uncomplicated medical (UM), complicated medical (CM), complicated surgery (CS), and joint surgery (JS) (Table~\ref{tbl1:drg}).  We also consider five different SNFs labeled \textbf{A} through \textbf{E} and enumerated 1 through 5, respectively (Table~\ref{tbl1:snf}). As explained in Section~\ref{sec:paramest}, we chose these SNFs for their proximity to the hospital to mitigate the impact of patient preferences in choosing specific SNFs based on location. We use estimated readmission rates from Table~\ref{tbltoss} for $r^j_i$ and, to begin with, we fixed discharge rates to be the same (i.e.,  $\lambda_i = 0.2, i \in \{\text{UM},\text{CM},\text{CS},\text{JS}\}$) so that each patient is equally likely to be discharged in any period.  We set the discharge rate for each patient type to be equal to more clearly delineate the impact of readmission risks and transition probabilities on the performance of the two heuristic policies.  In each experiment throughout this section, we conduct a sensitivity analysis on discharge rates and report the result of these experiments in the Appendix (see Sections~\ref{appendix:3}-~\ref{appendix:7} of the Appendix). Finally, each patient type can be assigned to any of the five SNFs throughout the study.  

\subsection{Transition probability matrices}

The remaining parameters that need to be specified are the transition probability matrices, which control the  availability of the five SNFs as a function of which SNF receives the patient. Recall that the $s,s'$th entry of a transition probability matrix represents the probability that the corresponding SNF's availability will be $s'$ at the next time period given that it is $s$ in the current period, where $s,s' \in \{0,1\}$, and where $0$ represents the case when SNF is not available and 1 represents the case that the SNF is available. Hence, a total of 25 2-by-2 transition matrices need to be specified for each simulation instance.  Ideally, these transition probabilities could be calculated from data if hospitals collect, over time, how often each SNF is (or is not) available for each patient type after discharge by using standard maximum likelihood estimation approaches (c.f, \cite{billingsley1961statistical}) or extensions thereof. In our case, our partner hospital did not keep records on the outcome from contact with SNFs about capacity availability. We thus chose to generate SNF occupancy rates via transition probabilities.  In generating transition probability matrices, we aimed to capture a range of behaviors observed in our partner hospital. Namely, we observed that while capacity of some SNFs were near-full over time, indicating that the availability depends highly on the transfer decisions, others have low occupancy rates indicating that SNF availability does not depend on allocation decision. We used different matrix structures to capture all of these observations in the data. For the clarity of the presentation, we assumed matrices are stationary, whereas real-life data from our partner hospital show that they may vary over time.      

Based on these considerations, we generated transition matrices in three scenarios (see Figure~\ref{fig:overview}B for a depiction). For each scenario, we generated 5 \emph{baseline} transition probability matrices corresponding to each of the 5 SNFs which are denoted by:
\begin{align*}
    \mathbf{\hat{P}}^{A},\,\, 
    \mathbf{\hat{P}}^{B},\,\,     
    \mathbf{\hat{P}}^{C},\,\, 
    \mathbf{\hat{P}}^{D},\,\,     
    \mathbf{\hat{P}}^{E}.
\end{align*}
With a slight abuse of notation, we also refer to these matrices as $\mathbf{\hat{P}}^{1}$, $\mathbf{\hat{P}}^{2}$, $\mathbf{\hat{P}}^{3}$, $\mathbf{\hat{P}}^{4}$,  $\mathbf{\hat{P}}^{5}$ to simplify latter expressions when we match matrices to (numbered) actions. Matrices were generated by sampling a uniformly distributed random variables $U (0, 1)$ for each diagonal entry of the transition probability matrix and defining the off-diagonal entry so that rows sum to $1$. 

As discussed below, the dependency structure is controlled in each of the three scenarios by the parameters listed in Table~\ref{tab:dependency_parameters}. The associated values we will consider are also listed in the table. 

\begin{table}[ht]
\centering  \footnotesize 
        \begin{tabular}{l l l}
             \toprule
            \textbf{Scenario} & \textbf{Parameter} & \textbf{Value(s)} \\
             \midrule
            \textbf{1} &  $\beta_1$ & $0.1,0.15,\ldots,0.95$ \\
            \midrule
           \multirow{2}{*}{\textbf{2}} & $\beta_2$ & $0.1,0.25,0.5,0.75,0.9$ \\
           & $\gamma_2$ &  $1,2,5,7,10$ \\
           \midrule
	      	\multirow{3}{*}{\textbf{3}}  & $\beta_3$ & $0.25,0.5,0.9$ \\
            &  $\gamma_3$ & $4,6.5,8$ \\
            &  $\delta_3$ & $1,1.75$ \\
 \bottomrule
        \end{tabular}
    \caption{Parameters controlling dependency between transfer decision and SNF availability.} \label{tab:dependency_parameters}
\end{table}

\subsubsection{Scenario 1:} The transfer decision only impacts the SNF to which the transfer decision is made, but not the others. In particular, the baseline matrices are used to govern the availability of the corresponding SNF when a patient is not assigned to this SNF:
\begin{align*}
    \mathbf{P}^{a,j} := \mathbf{\hat{P}}^{j} \qquad \qquad a \neq j. 
\end{align*}
For the remaining transfer decisions, we first set $\mathbf{P}^{a,a} = \mathbf{\hat{P}}^{a}$ and then alter the probability of staying available: $p^{a,a}_{1,1} = \beta_1 \, \hat{p}^a_{1,1}$ where  $\beta_1 \in [0,1]$ represents the dependency of the SNF capacity to the transfer decisions. We also modify off-diagonal entry $p^{a,a}_{1,0} = 1- p^{a,a}_{1,1}$ so that the corresponding row in the matrix sums to 1. 

\subsubsection{Scenario 2:} The transfer decision impacts the SNF to which the patient is transferred to as well as the `neighboring' SNFs.  Considering SNF labels as a sequence \textbf{A}, \textbf{B}, \textbf{C}, \textbf{D}, \textbf{E}, a neighbor is defined as a SNF that follows or precedes the given SNF in the sequence. For example, \textbf{B} is a neighbor of \textbf{A}, and \textbf{C} and \textbf{E} are neighbors of \textbf{D}. As a consequence, a transfer to SNF \textbf{A} impacts the availability of SNF \textbf{B}, whereas a transfer to SNF \textbf{D} impacts the availability of SNF \textbf{C} and \textbf{E}. To construct the 25 transition probability matrices, we proceed like we did in Scenario 1. That is, we again generate 5 baseline transition probability matrices for each SNF, initially set the 25 transition matrices equal to the baseline matrix corresponding to the given SNF, and modify the probability of staying available of the receiving SNF except with $\beta_1$ replaced with $\beta_2/\gamma_2$ (i.e., $p^{a,a}_{1,1} = \frac{\beta_2}{\gamma_2} \, \hat{p}^a_{1,1}$). In addition, we also modify the probability of staying available of the neighboring SNFs according to $p^{a,j}_{1,1} = \beta_2 \, \hat{p}^{j'}_{1,1}$ where $j'$ denotes a neighbor of SNF $j$. The off-diagonal entry is adjusted so that rows sum to 1. 

We let $\beta_2 \in [0,1]$ and $ \gamma_2 \in [1,\infty)$. Akin to scenario 1, $\beta_{2}$ measures the dependency of the SNF availability on the transfer decisions:  all else fixed, higher $\beta_{2}$ values increase the likelihood of the SNF becoming available at the next period.  The parameter, $\gamma_{2}$, also measures dependence but with higher values of $\gamma_{2}$ decreasing the likelihood of the SNF becoming available at the next period.  Parameters $\beta_2$ and $ \gamma_2$ are used to also influence the impact of the transfer decisions on the availability of the receiving SNF and its neighbors. 

\subsubsection{Scenario 3:}  The transfer decision impacts the availability of all SNFs. We proceed as in Scenario 2 to construct the 25 transition probability matrices except we replace $\beta_2$ and $\gamma_2$ with $\beta_3$ and $\gamma_3/\delta_3$. Here, $\beta_3 \in [0,1]$, $\gamma_3 \in [1,\infty)$, and $\delta_3 \in [1,\gamma_3]$. We then modify the probability of SNF availability of the SNFs that are neither the receiving SNF nor the neighbors according to $p^{a,j''}_{1,1} = \beta_3 \,\hat{p}^{j''}_{1,1}$ where $j''$ denotes a non-neighboring SNF of SNF $j$. The off-diagonal entry is adjusted so that rows sum to 1. Parameter $\delta_3$ provides an additional measure of dependence:  increasing its value decreases the likelihood of the SNF being available at the next time period. We bound $\delta_3$ above by $\gamma_3$ to reflect that the transfer decision will impact the availability of the receiving SNF the most, followed by neighboring SNFs, and then the remaining SNFs. 

\subsection{Quantitative analysis} 
\label{sec:Quant-analysis}
\begin{figure}
\centering
\includegraphics[height=2.3in]{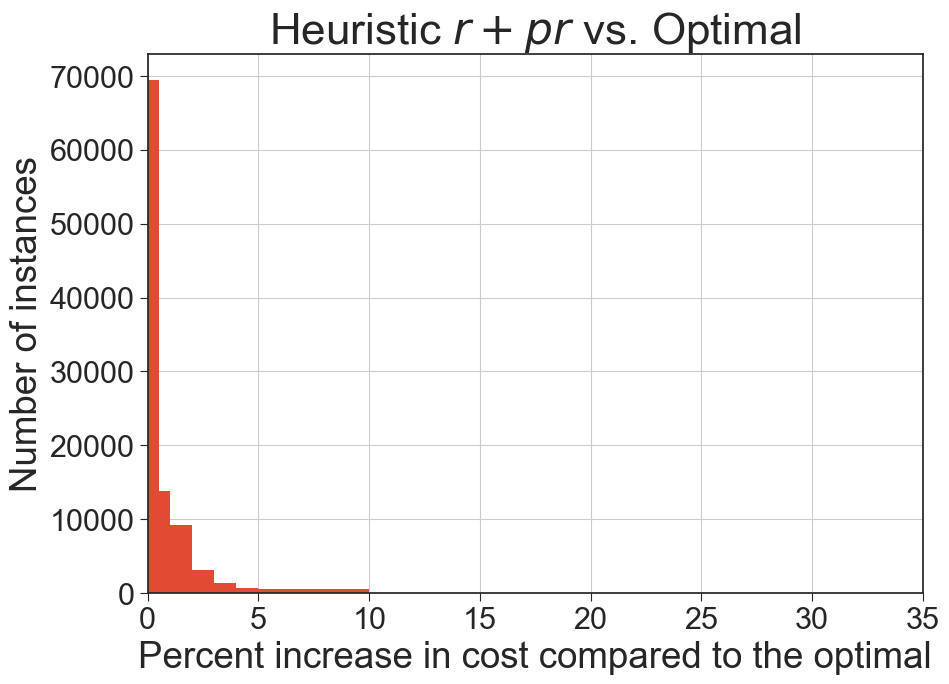}
\includegraphics[height=2.3in]{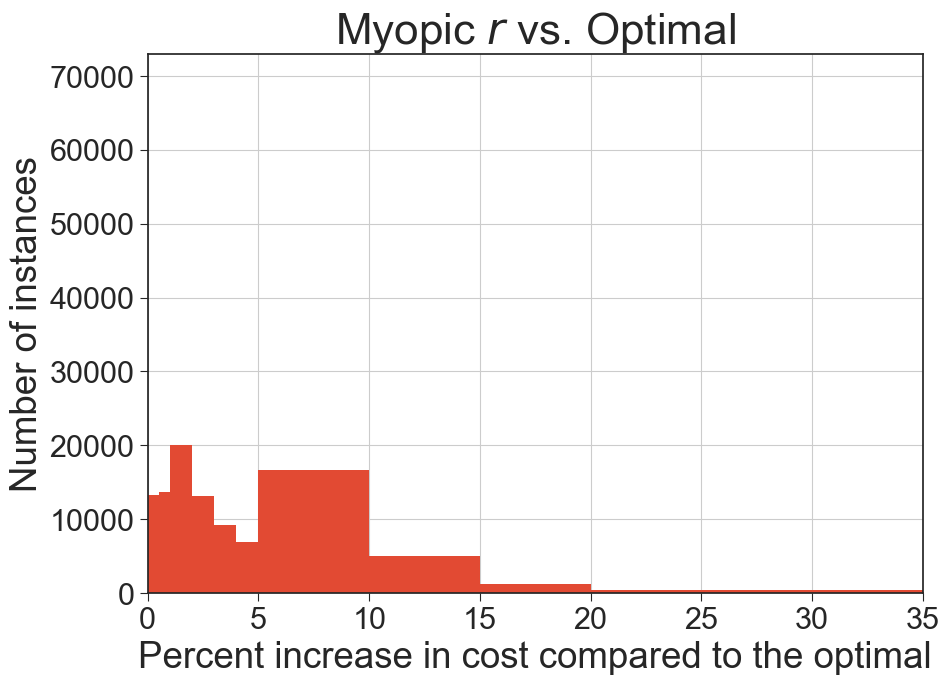} \\
\includegraphics[height=2.3in]{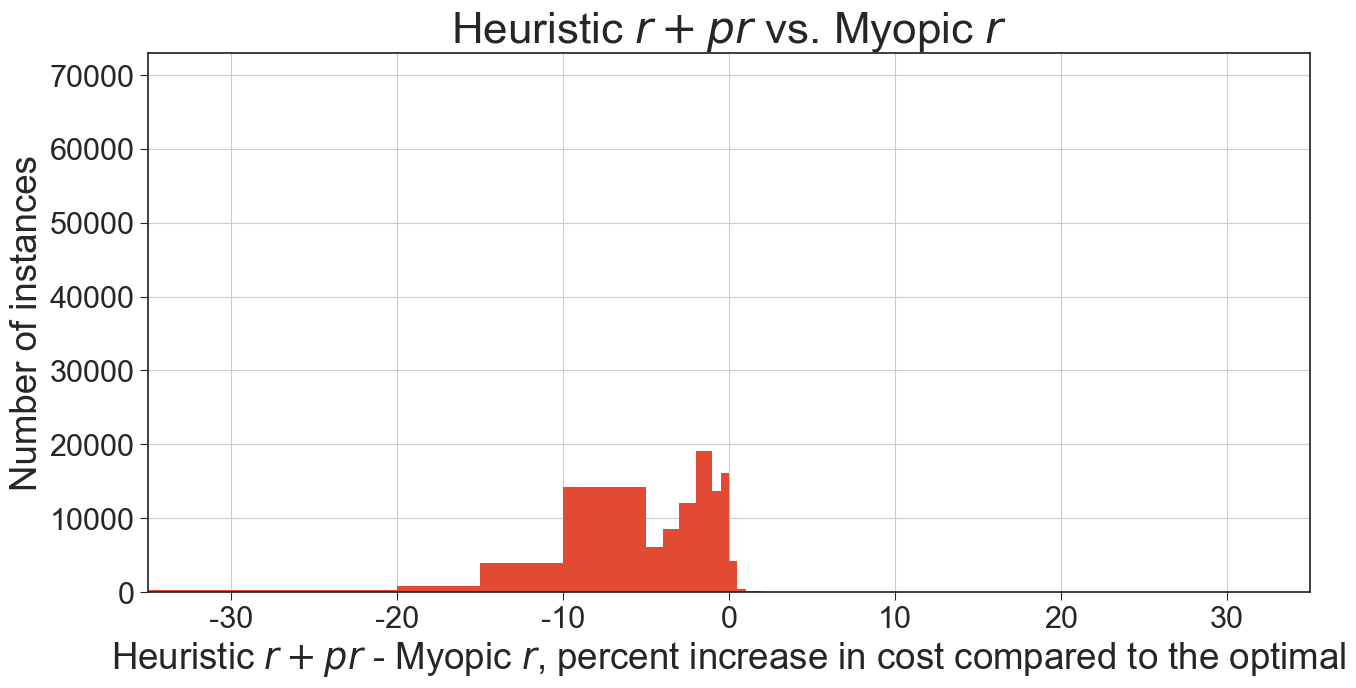}
\caption{Comparison of heuristics $r$ and $p+r$ under Scenario 1 with $\beta_1=0.2$ in terms of percent increase in expected long-term average costs relative to the optimal average costs.} \label{fig:comparisonScen1} 
\end{figure}

We generated 100,000 instances of the 25 transition probability matrices following the procedure described under Scenario 1 in the preceding subsection but with $\beta_1 = 0.2$ only. The results of these experiments are summarized in Figure~\ref{fig:comparisonScen1} where we plot the number of instances against percent increase in cost compared to the optimal long-run average readmission rates when heuristic $r+p r$ is used (top-left panel), percent increase in cost compared to the optimal long-run average readmission rates when myopic policy $r$ is used (top-right panel), and the difference between the percent increase in average readmission rates when heuristic $r+pr$ is used minus the percent increase in average readmission rates when the myopic policy $r$ is used (bottom panel).  Overall, we find that policy $r+pr$ outperformed the myopic policy $r$ in roughly 96\% of the instances.  We also find that the average readmission rates obtained by heuristic $r+p r$ are are within 1\% of those obtained by the optimal policy in more than 80,000 instances. Furthermore, the percent increase in average readmission rates obtained by heuristic $r+p r$ compared to those obtained by the optimal policy never exceeded 10\%.  By contrast, the percent increase in average readmission rates obtained by the myopic policy compared to those obtained by the optimal policy reached 35\% in some cases.

We also generated 10,000 instances of the 25 transition probability matrices following the procedure described under Scenario 2 in the preceding subsection but only using $\beta_2=0.2$ and $\gamma_2=5$, and generated 10,000 instances of the 25 transition probability matrices following the procedure described under Scenario 3 but only with $\beta_3=0.2$, $\gamma_3=5$, and $\delta_3=  1.75$ fixed.  The results of these experiments are summarized in Figures~\ref{fig:comparisonScen2} and~\ref{fig:comparisonScen3} of Section~\ref{appendix:2} of the Appendix, where we observe similar performance of each heuristic policy as in Scenario 1.  For example, when instances are generated following the procedure described under Scenario 2, policy $r+pr$ outperforms the myopic policy $r$ roughly 96\% of the time.  When instances are generated following the procedure described under Scenario 3, policy $r+pr$ outperforms the myopic policy $r$ roughly 88\% of the time.  For completeness, we also evaluated a natural extension of policy $r+p r$ that accounts for SNF availability two periods into the future instead of one. Details are provided in Section~\ref{appendix:2} of the Appendix.

\subsection{Qualitative analysis}

Among the 100,000 instances sampled above, we selected a few to provide insights on when the policies might perform well. We remark that Section~\ref{sec:Quant-analysis} includes many instances that have similar characteristics to the ones considered in this section. The purpose of this section is to highlight common characteristics of the instances included in Section~\ref{sec:Quant-analysis} that lead to superior or inferior performance of each policy.

\subsubsection{When the myopic policy $r$ performs well.} We considered two instances when the myopic policy $r$ outperforms heuristic $r+p r$. In the first instance, the optimal, myopic $r$, and heuristic $r+pr$ policies have average costs (i.e., average readmission rates) of $14.44$, $14.47$, and $15.14$, respectively. Therefore, average costs for the $r$ and $r+pr$ policy are 0.2\% and 4.8\% higher than the average cost of the optimal policy. In the second example, the optimal, myopic $r$, and heuristic $r+pr$ policies have average costs of $14.54$, $14.95$, and $15.78$, respectively. Hence, average costs for the $r$ and $r+pr$ heuristics are 2.8\% and 8.5\% higher than the average cost of the optimal policy. Due to limited space, transition probability matrices for these instances, as well as a detailed investigation into situations when transfer decisions impact SNF availability (i.e., Scenarios 1--3), are reported in Sections~\ref{appendix:3}--\ref{appendix:4} of the Appendix. Policies are also reported by state (there are a total of 384 states) in Tables~\ref{tab:ExampleAppendix4}--~\ref{tab:ExampleAppendix5} of the Appendix. 

These instances highlight that the myopic policy $r$ can still perform better than heuristic $r + p r$, despite only minimizing immediate readmission rates as opposed to immediate readmission rates \emph{and} one-step future expected readmission rates. More specifically, the myopic policy $r$ performs well either when the readmission rate of one SNF for each patient type is significantly lower than those of the other SNFs or when all SNFs are highly likely to be unavailable in one step regardless of the transfer decision. In these cases, little is gained over a myopic policy by planning ahead. If these two conditions hold and, in addition, any SNF with a readmission rate close to the lowest readmission rate is likely to be unavailable in two or more two periods, then myopic policy $r$ can perform better than heuristic policy $r + p r$. Heuristic $r + p r$ only considers that sending a patient to a specific SNF may decrease its availability at the next period, but does not account for the fact that the same SNF can quickly become available after that, while the other SNFs may be unavailable due to their similar long-term availability. By only accounting for the next possible period, heuristic $r + p r$ is unable to capture such longer-term behavior. The two instances illustrate that in such scenarios, it underperforms the myopic policy and may not be close to optimal.

\subsubsection{When the heuristic policy $r + p r$ performs well.}  We considered two instances when heuristic $r+p r$ outperforms the myopic policy $r$. In the first, the optimal, myopic $r$, and heuristic $r+pr$ policies have average costs of $13.1$, $17.75$, and $13.11$, respectively. Therefore, average costs for the $r$ and $r+pr$ policies are 35.5\% and 0.07\% higher than the average cost of the optimal policy. In the second, the optimal, myopic $r$, and heuristic $r+pr$ policies have average costs of $12.49$, $14.64$, and $13.12$, respectively. Hence, average costs for the $r$ and $r+pr$ heuristics are 17.2\% and 5.0\% higher than the average cost of the optimal policy. In each instance, heuristic $r+ p r$ performs well (while the myopic policy $r$ performs poorly) when the SNFs with low readmission rates are highly likely to become unavailable in the long-term while the remaining SNFs are not. In these cases, it is better to use heuristic $r + p r$. For these instances, transition probability matrices, as well as detailed investigation into situations when transfer decisions impact SNF availability, are reported in Sections~\ref{appendix:5}--\ref{appendix:6} of the Appendix. The resulting policies are provided in Tables~\ref{tab:ExampleAppendix6} and~\ref{tab:ExampleAppendix7} of the Appendix. 

\subsubsection{When the optimal policy outperforms heuristics $r$ and $r+ p r$.} \label{section:exopt} We considered one instance when the optimal policy significantly outperforms the myopic policy and heuristic $r+p r$. For this instance, the average cost for the optimal, myopic, and $r+pr$ policies are $14.24$, $17.88$, and $15.69$, respectively.  Thus, average costs for the myopic $r$ and $r+pr$ policies, respectively, are 25.6\% and 10.2\% higher than average costs of the optimal policy. This instance illustrates that when SNF availability status changes quickly (i.e., in one period) from unavailable to available and vice-versa, and when all SNFs have similar long-run availability, it is better to use the optimal policy from our MDP model. Transition probability matrices for this example are given in Section~\ref{appendix:7} of the Appendix. Resulting policies are provided in Table~\ref{tab:ExampleAppendix8} of the Appendix. We proceed to investigate situations when transfer decisions impact SNF availability (Scenarios 1--3) for this instance.

\subsection{Dependency structure}

Here, we want to better characterize the role that dependency between transfers decisions and SNF availability has on heuristic performance. We use the same transition probability matrices from the last instance in Section~\ref{section:exopt}, which represents the case that the optimal policy outperforms heuristics $r$ and $r+ p r$:
\begin{align*}
& \mathbf{\hat{P}}^{A} = \begin{bmatrix}
0.56&0.44\\
0.94&0.06
\end{bmatrix},
\mathbf{\hat{P}}^{B} = \begin{bmatrix}
0.91&0.09\\
0.78&0.22
\end{bmatrix},
\mathbf{\hat{P}}^{C} = \begin{bmatrix}
0.36&0.64\\
0.68&0.32
\end{bmatrix},\\
& \mathbf{\hat{P}}^{D} = \begin{bmatrix}
0.97&0.03\\
0.02&0.98
\end{bmatrix},
\mathbf{\hat{P}}^{E} = \begin{bmatrix}
0.88&0.12\\
0.01&0.99
\end{bmatrix},
\end{align*}
where each entry of the matrix shows the probability that SNF availability changes in each time period. For example, the entry in the first row and first column of $\mathbf{\hat{P}}^{A}$ would denote that SNF \textbf{A} would remain unavailable at the next decision epoch with probability $0.56$. We modify these transition matrices according to Scenarios 1--3 for various parameters listed in Table~\ref{tab:dependency_parameters}. For each scenario, we report the relative gap (in \%) between the expected long-run average costs of each heuristic and the optimal policy, i.e., $\frac{g^{\text{heuristic}}-g^{\text{optimal}}}{g^{\text{optimal}}}\cdot100$. 
   
\subsubsection{Performance of heuristic policies under Scenario 1.}

\begin{figure}[t!]
\centering
\includegraphics[width=0.6\textwidth]{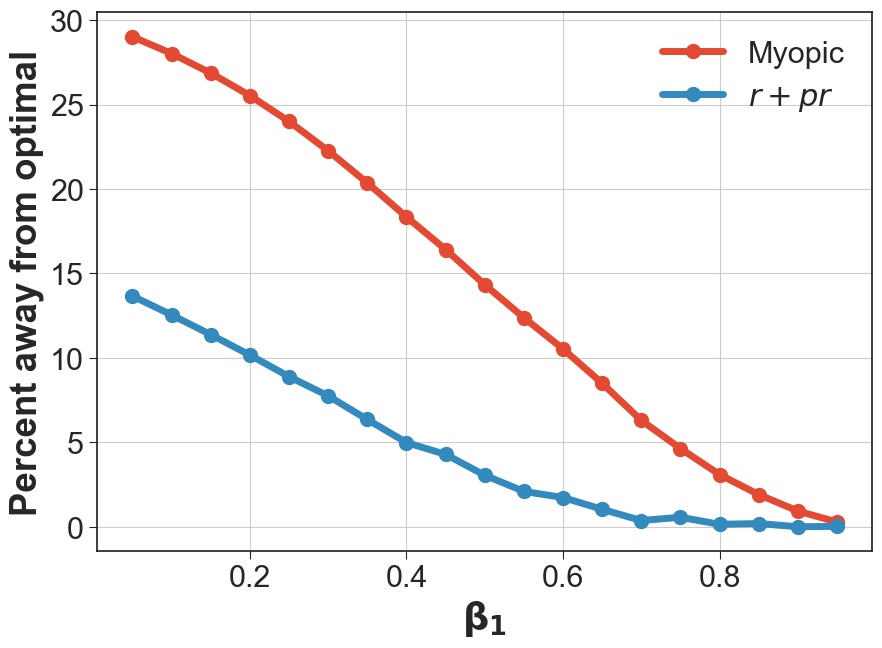}
\caption{Under Scenario 1, the relative gap (in \%) of the expected long-run average costs between each heuristic (myopic $r$ and heuristic $r+pr$) and the optimal policy as a function of $\beta_1$.} \label{fig:scenario1}
\end{figure}

Both heuristics do increasingly worse as we decrease the probability that the receiving SNF stays available via $\beta_1$ (Figure~\ref{fig:scenario1}). At worst ($\beta=0.1$), the myopic policy $r$ led to expected long-term average costs that were 29\% higher than those of the optimal policy; at best, the gap is 0.3\% ($\beta=0.95$). This is expected, considering that a myopic policy is optimal when the transfer decision does not influence transitions (i.e., when $\beta_1=1$). By comparison, heuristic $r+pr$ outperforms the myopic policy for all $\beta_1$ values and degrades more slowly than the myopic policy with decreasing $\beta_1$. At worst, heuristic $r+pr$ led to expected long-term average costs that were 13.7\% ($\beta_1=0.1$) higher than the optimal expected long-term average costs, and, at best, it led to expected long-term average costs that were 0.03\% higher than the optimal expected long-term average costs ($\beta_1=0.95$). 

\subsubsection{Performance of heuristic policies under Scenario 2}
\begin{figure}
\centering
\includegraphics[width=0.9\textwidth]{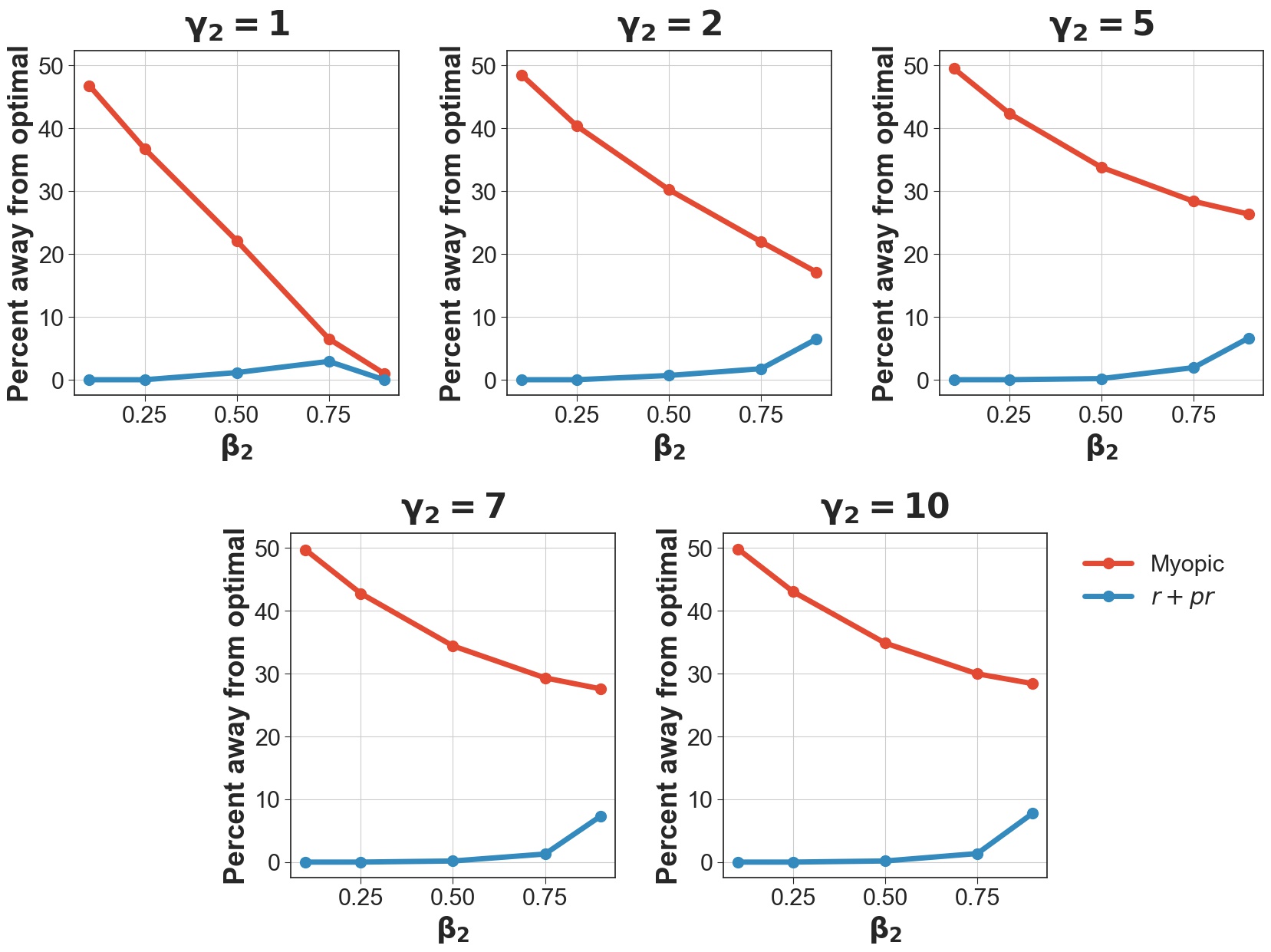}
\caption{Under Scenario 2, the relative gap (in \%) of the expected long-run average costs between each heuristic (myopic $r$ and heuristic $r+pr$) and the optimal policy as a function of $\beta_2$ and $\gamma_2$.} \label{fig:scenario2}
\end{figure}

In contrast to the above scenario, the transfer decision in Scenario 2 not only impacts the SNF to which the patient is sent, but also neighboring SNFs via parameters $\beta_2$ and $\gamma_2$. For this scenario, we find that performance of the myopic policy is extremely sensitive to both $\beta_2$ and $\gamma_2$ (Figure~\ref{fig:scenario2}). Similar to Scenario 1, the myopic policy performs worse with decreasing availability of the receiving SNF (i.e., as $\gamma_2$ increases). The myopic policy also performs worse with decreasing availability of the receiving SNF and its neighbors (i.e., as $\beta_2$ decreases). Expected long-term average costs were, at worst, 50\% higher for the myopic policy than the optimal expected long-term average cost ($\beta_1=0.1$, $\gamma_2=10$); and at best, within 1\% ($\beta_2=0.9$, $\gamma_2=1$) of the optimal cost.

Heuristic $r+pr$, however, outperformed the myopic policy in all scenarios and was relatively less sensitive to $\beta_2$ and $\gamma_2$ compared to the myopic policy. In particular, the expected long-term average costs for heuristic $r+pr$ were within 2\% of the optimal cost in 20 out of the 25 instances.  At worst, the expected long-term average cost for heuristic $r+pr$ was 7.7\% higher than the optimal expected long-term average cost.  Surprisingly, increasing the dependency of the transfer decision only on the receiving SNFs by increasing $\gamma_2$ does not greatly affect the performance of heuristic $r+pr$ except when the dependency of the transfer decision on neighboring SNFs is low (i.e. $\beta_2=0.9$). 

\subsubsection{Performance of heuristic policies under Scenario 3.}

\begin{figure}
\centering
\includegraphics[width=0.9\textwidth]{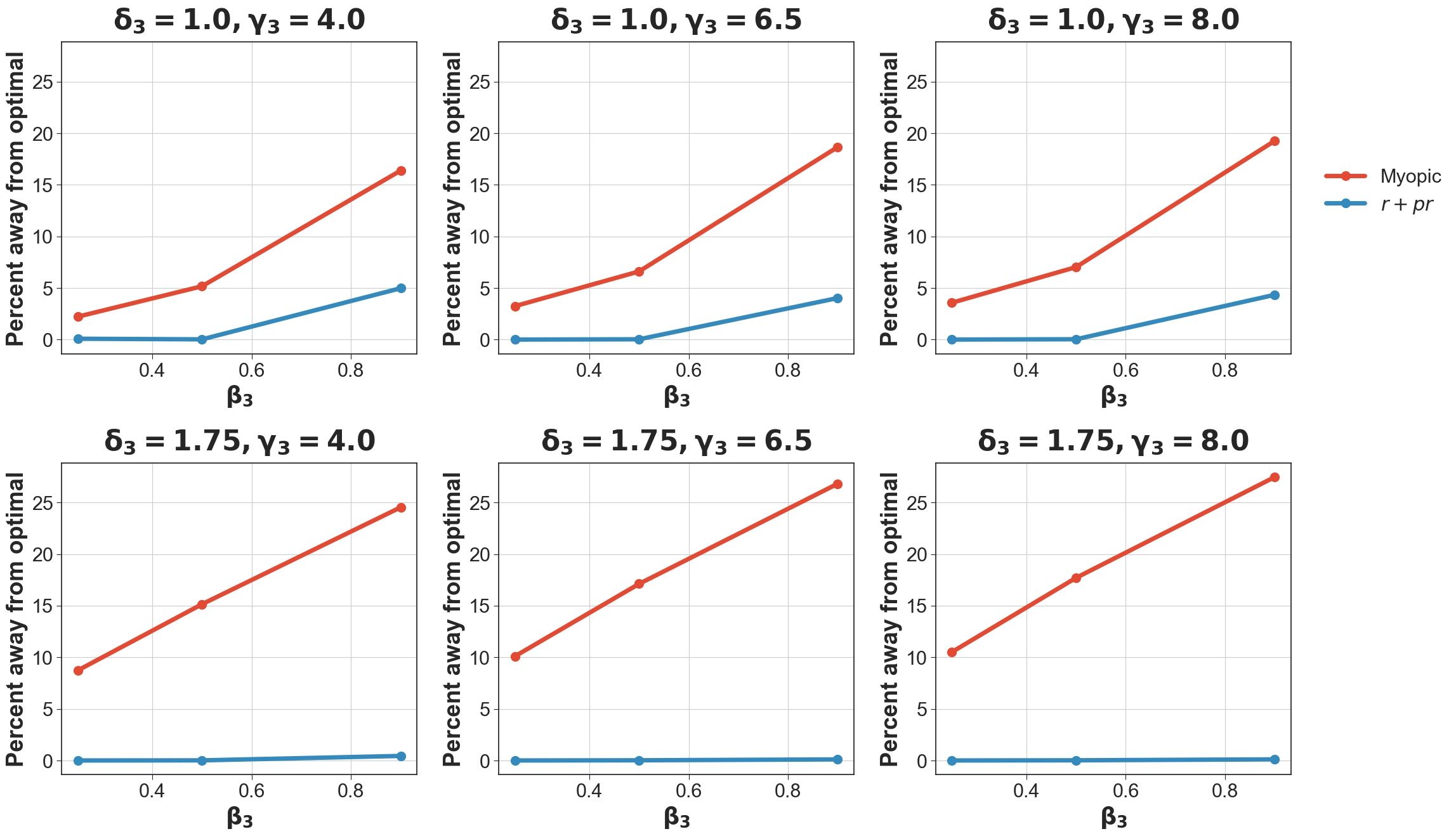}
\caption{Under Scenario 3, the relative gap (in \%) of the expected long-run average costs between each heuristic (myopic $r$ and heuristic $r+pr$) and the optimal policy as a function of $\beta_3$, $\gamma_3$, and $\delta_3$.} \label{fig:scenario3}
\end{figure}

Our last scenario considers when the transfer decision impacts the availability of all SNFs via parameters $\beta_3$, $\gamma_3$, and $\delta_3$. The myopic policy performs slightly better in this scenario compared to prior scenarios, but still worse than heuristic $r + pr$ (Figure~\ref{fig:scenario3}). Its expected long-term average cost is 27.5\% higher than the optimal expected long-term cost ($\beta_3=0.9$, $\gamma_3=8.0$, $\delta_3=1.75$) and within 2.3\% ($\beta_3=0.25$, $\gamma_3=4$, $\delta_3=1$) of the optimal value.  Interestingly, decreasing the dependency of the availability of \emph{all} SNFs on the transfer decision by increasing $\beta_3$ has the surprising effect of worsening performance of the myopic policy. This suggests that it is not broadly how transfer decisions impact SNF availability, but specifically how transfer decisions impact SNF availability on the receiving SNF and its neighbors relative to other SNFs (as was explored in Scenario 1 and 2). To that point, we also observe that decreasing the availability of only the receiving SNF by increasing $\delta_3$ leads to a worsening performance of the myopic policy. 

Meanwhile, expected long-term average costs for heuristic $r+pr$ were within 0.5\% of the optimal in all but three scenarios; at worst, its expected long-term average cost was 4.98\% higher than the optimal expected long-term cost. As before, heuristic $r+pr$ is less sensitive to the changes in the dependency structure than the myopic policy, helping to protect against poor performance. Yet, at the same time, the dependency structure is still important, as one can change the performance of heuristic $r+pr$ by increasing $\beta_3$ while holding $\delta_3$ fixed at $1$.

\section{Conclusion} \label{sec:conclusion}

In this paper, we proposed and analyzed transfer policies for patient care transitions between hospital after a patient discharge and SNFs.  To this end, we used an MDP formulation that minimizes discounted and long-run average readmission rates which were estimated using electronic health records from a large tertiary academic hospital.  The Patient Protection and Affordable Care Act prompted scrutinization of hospital quality of care by holding hospitals accountable for patient outcomes after discharge.  This was partly done through financial penalties on hospitals with excess readmission rates as part of Medicare's Hospital Readmissions Reduction Program.  As a result, hospitals have focused on post-acute care to reduce readmissions and SNFs represent one of the most common setting for post-acute care in the US, but still result in high rates of readmissions.  The present work is based on the premise that one potential way to improve readmissions is to send patients to the most appropriate SNF after discharge.

For this purpose, we estimated readmission rates by SNF and patient types using observational data from a large tertiary teaching hospital and nearby SNFs.  We then developed, analyzed, and solved a sequential decision-making model of patient transfer decisions that minimizes readmissions using the aforementioned estimated readmission rates.  Our model accounts for patient discharge patterns and SNF capacity availability.  We provide conditions for when an easy-to-use myopic policy, which assigns discharged patients to the SNF corresponding to the lowest readmission rate that is available, is optimal. We also show when an optimal policy has a threshold-like structure. Other than these cases, structural properties of the optimal policy are difficult to characterize.  Further, the model suffers from the curse of state-space (and action-space) dimensionality, making approximation of optimal policies difficult. Using estimated readmission rates from real-life data, we compare the performance of the myopic policy and a proposed policy that minimizes the sum of immediate readmission rates and one-step expected readmission rates under the the myopic policy.  In contrast to the myopic policy, this policy depends on the discharge process, readmission rates, and future SNF availability.  We evaluate when the myopic policy may be beneficial to use and, otherwise, when the proposed transfer heuristic provides a better alternative.  Other than these cases, we contend that using MDPs for guiding these transfer decisions may help reduce readmissions.

Our analysis produced several useful insights about transfer decisions.  First, when the transfer decision has little impact on SNF availability, the myopic policy should be used.  Second, when either one SNF is significantly better for a patient than all others or all SNFs are likely to be unavailable regardless of the transfer decision, then the myopic policy $r$ performs well. If, in addition, any SNF with a readmission rate close to the lowest readmission rate is likely to be unavailable in two or more two steps, then myopic policy $r$ can perform better than heuristic policy $r + p r$.  Third, when SNFs with higher readmission rates are much more likely to be available in the long-run than the remaining SNFs, heuristic $r+pr$ performs well (while the myopic policy does not).  Fourth, when SNF availability status changes quickly (i.e., in one period), and when all SNFs have similar long-run availability, it is better to use the optimal policy.  Lastly, we have also shown that as SNF availability becomes increasingly dependent on the transfer decision, the myopic policy's performance worsens whereas the proposed heuristic $r+pr$ generally performs well (at worst, the performance gap between average readmission rate and the optimal readmission rate never exceeds 10\%). Therefore, heuristic $r+pr$ provides a viable option for implementation in practice, especially when the number of patient types and the number of SNFs are not very large. 

Considering that these insights are sensitive to modeling assumptions, it is important to examine whether alternative assumptions or modeling approaches converge on the same conclusions in an effort to strengthen the validity of these conclusions. These guiding principles should be interpreted with caution as one does with all optimization approaches, before making recommendations for continuing, adapting, and/or adopting new strategies for transferring patients to SNFs.  For example, we used observational data for estimating readmission rates but data on how a transfer decision impacts SNF availability, or SNF availability in general, was not available.  Additionally, the description of the state space included whether a particular SNF was available or not.  It is possible  that alternative descriptions of the state space are more relevant, such as one that includes a measure of congestion in each SNF (e.g., the number of each patient type currently receiving care in each SNF).  We did not so here because congestion is not usually available to the social worker when making the SNF recommendation to the patient, and importantly, because congestion is often not the main reason for a SNF being available for a particular patient type and other factors, such as patient type and insurance coverage play an important role.  
Finally, given that estimates of readmission rates were obtained from observational data, transfer policies based on robust MDP formulations may be used. 

In summary, our work is the first to consider transfer decisions between hospital and SNFs using a sequential decision-making model with costs estimated from observational data. Finding overall reductions in readmission rates from transfer decisions using our proposed approach suggests that including measures of readmission combined with the probabilistic nature of SNF availability can improve care and reduce costs.  Our results do not support a change in current clinical practice, especially in light of the considerations and potential extensions mentioned above; rather they contribute to a new body of literature suggesting that embedding more realistic models may yield benefits, and point to the necessity of further assessment to identify specific criteria to guide the adoption of such models. Lastly, our approach and analysis can be used to optimize transfer decisions to other rehabilitation centers or patient homes.




\bibliographystyle{unsrtnat}

\bibliography{references}





\appendix
\section{Appendix}

\subsection{Proofs of Propositions~\ref{prop:dcoe}-~\ref{prop:acoe}}
\label{appendix:0}
\textbf{Proof of Proposition~\ref{prop:dcoe}}
Since $|\X|<\infty$ and $|A(x)|<\infty$ for each $x \in \X$, the proof is an immediate consequence of Theorem 6.2.10 of~\cite{puterman94}.

\noindent
\textbf{Proof of Proposition~\ref{prop:acoe}}
The conditions $0<\lambda_i<1$ for all $i$ and $0 <p^{a,j}_{s,s'} < 1$ for all $a$, $j \geq 1$, $s$, and $s'$ imply that, for any stationary deterministic policy, all states of the form $(i,1,s_1,\ldots,s_{\ell})$ are reachable from any other state and that these states communicate.  This yields that there is at most one recurrent class.  Since $r^a_i$ are uniformly bounded in $a$ and $i$ and $|\X|<\infty$ and $|A(x)|<\infty$ for each $x \in \X$, the result follows from Theorem 8.4.3 of~\cite{puterman94}.

\subsection{Optimal Policy and Myopic Policy for Examples~\ref{ex1} and \ref{ex2}} \label{appendix:0.5}

\begin{table}[htbp]
\centering
\begin{tabular}{l c c c c c c }
\toprule
\multicolumn{1}{c}{\textbf{$i$}} & \multicolumn{1}{c}{\textbf{$s_1$}} & \multicolumn{1}{c}{\textbf{$s_2$}} & \multicolumn{1}{c}{\textbf{Optimal Policy}} & \multicolumn{1}{c}{\textbf{Myopic Policy}} \\
    \midrule    1 & 0 & 0 & 0 & 0\\
    1 & 0 & 1 & 2 & 2 \\
    1 & 1 & 0 & 1 & 1 \\
    1 & 1 & 1 & 1 & 1 \\
    2 & 0 & 0 & 0 & 0 \\
    2 & 0 & 1 & 2 & 2 \\
    2 & 1 & 0 & 1 & 1 \\
    2 & 1 & 1 & 1 & 2 \\
\bottomrule
\end{tabular}
\caption{Actions chosen in each state by the optimal policy and the myopic policy for instance described in Example~\ref{ex1}.}
\label{tab:ex1}
\end{table}

\begin{table}[htbp]
\centering
\begin{tabular}{l c c c c c c }
\toprule
\multicolumn{1}{c}{\textbf{$i$}} & \multicolumn{1}{c}{\textbf{$s_1$}} & \multicolumn{1}{c}{\textbf{$s_2$}} & \multicolumn{1}{c}{\textbf{Optimal Policy}} & \multicolumn{1}{c}{\textbf{Myopic Policy}} \\
    \midrule    1 & 0 & 0 & 0 & 0 \\
    1 & 0 & 1 & 2 & 2 \\
    1 & 1 & 0 & 1 & 1 \\
    1 & 1 & 1 & 2 & 1 \\
    2 & 0 & 0 & 0 & 0 \\
    2 & 0 & 1 & 2 & 2 \\
    2 & 1 & 0 & 1 & 1 \\
    2 & 1 & 1 & 2 & 2 \\
\bottomrule
\end{tabular}
\caption{Actions chosen in each state by the optimal policy and the myopic policy for the instance described in Example~\ref{ex2}.}
\label{tab:ex2}
\end{table}

\subsection{Proofs of Propositions~\ref{prop:1} and~\ref{prop:3}} \label{appendix:1}

\noindent
\textbf{Proof of Proposition~\ref{prop:1}.}  

\noindent
By assumption,
\begin{align*}
 \sum_{x'=(i',1,s'_1,\ldots,s'_{\ell}) \in \X} \lambda_{i'} \Pi_{j=1}^\ell p^{a,j}_{s_j,s'_j} v_{\alpha}(x') = \sum_{x'=(i',1,s'_1,\ldots,s'_{\ell}) \in \X} \lambda_{i'} \Pi_{j=1}^\ell p^{j}_{s_j,s'_j} v_{\alpha}(x')
\end{align*}
for all $a$ and for all $s'_1$, $s_2$', $\ldots$, and $s'_{\ell}$.  As a result, we can simplify the optimality equations
\begin{align*}
 v_{\alpha}(x) & = \min_{a \in A(x) } \{ r^a_i + \alpha \sum_{x'=(i',1,s'_1,\ldots,s'_{\ell}) \in \X} \lambda_{i'} \Pi_{j=1}^\ell p^{a,j}_{s_j,s'_j} v_{\alpha}(x') \} 
 \end{align*}
 to
 \begin{align*}
 v_{\alpha}(x) & = \min_{a \in A(x) } \{ r^a_i + \alpha \sum_{x'=(i',1,s'_1,\ldots,s'_{\ell}) \in \X} \lambda_{i'} \Pi_{j=1}^\ell p^{j}_{s_j,s'_j} v_{\alpha}(x') \} 
\end{align*}
for all $x=(i,s_0,s_1,\ldots,s_{\ell})$.  This last expression implies that immediate actions only impact the immediate costs incurred via $r^a_i$ and not the future evolution of SNF availability.  As a result, the myopic policy is optimal.

\noindent
\textbf{Proof of Proposition~\ref{prop:3}.}

\noindent
For the discounted expected cost problem, observe that the inequality 
\begin{align}
 & \alpha \sum_{x''=(i'',s''_1,\ldots,s''_{\ell}) \in \X} \lambda_{i'} \Pi_{j=1}^\ell p^{a^*,j}_{s'_j,s''_j} v_{\alpha}(x'') - \min_{a \in A(x')\setminus\{a^*\} } \{ r^a_i + \alpha \sum_{x''=(i'',s''_1,\ldots,s''_{\ell}) \in \X} \lambda_{i''} \Pi_{j=1}^\ell p^{a,j}_{s'_j,s''_j} v_{\alpha}(x'') \} \nonumber \\
 & \quad - \alpha \sum_{x''=(i'',s''_1,\ldots,s''_{\ell}) \in \X} \lambda_{i''} \Pi_{j=1}^\ell p^{a^*,j}_{s_j,s''_j} v_{\alpha}(x'') + \min_{a \in A(x)\setminus\{a^*\} } \{ r^a_i + \alpha \sum_{x''=(i'',s''_1,\ldots,s''_{\ell}) \in \X} \lambda_{i''} \Pi_{j=1}^\ell p^{a,j}_{s_j,s''_j} v_{\alpha}(x'') \} \nonumber\\
 & \qquad \leq 0 \label{ineq1:prop3}
 \end{align}
implies that
if 
$$r^{a^*}_i + \alpha \sum_{x''=(i'',s''_1,\ldots,s''_{\ell}) \in \X} \lambda_{i''} \Pi_{j=1}^\ell p^{a^*,j}_{s_j,s''_j} v_{\alpha}(x'') = \min_{a \in A(x)} \{ r^a_i + \alpha \sum_{x''=(i'',s''_1,\ldots,s''_{\ell}) \in \X} \lambda_{i''} \Pi_{j=1}^\ell p^{a,j}_{s_j,s''_j} v_{\alpha}(x'') \} ,$$ then 
$$r^{a^*}_i + \alpha \sum_{x''=(i'',s''_1,\ldots,s''_{\ell}) \in \X} \lambda_{i''} \Pi_{j=1}^\ell p^{a^*,j}_{s'_j,s''_j} v_{\alpha}(x'') = \min_{a \in A(x')} \{ r^a_i + \alpha \sum_{x''=(i'',s''_1,\ldots,s''_{\ell}) \in \X} \lambda_{i''} \Pi_{j=1}^\ell p^{a^*,j}_{s'_j,s''_j} v_{\alpha}(x'') \}.$$  
Intuitively, this implies that if it is optimal to allocate patient type $i$ to SNF $a^*$ in state $x=(i,s_0,s_1,s_2,\ldots,s_{\ell})$, then it is also optimal to allocate patient type $i$ to SNF $a^*$ in state $x'=(i,s_0,s_1',s_2',\ldots,s_n')$.  We prove~\eqref{ineq1:prop3} by induction assuming that $s_{a^*}=s'_{a^*}$ and $s_j \geq s'_j$ for $j \neq a^*$.

Let $v^0_{\alpha}(x) = 0$ for all $x \in S$ and recursively define $v^{n+1}_{\alpha}(x) = T v^n_{\alpha}(x)$ for all $x$ and $n=1,2,\ldots$.  Note that the inequality trivially holds for $n=0$.  Assume inequality~\eqref{ineq1:prop3} holds for $n \geq 0$ and consider it for $n+1$.  Let  $$a'=\argmin_{a \in A(x')\setminus\{a^*\} } \{ r^a_i + \alpha \sum_{x''=(i'',s''_1,\ldots,s''_{\ell}) \in \X} \lambda_{i''} \Pi_{j=1}^\ell p^{a,j}_{s'_j,s''_j} v^n_{\alpha}(x'') \}.$$ 
Fix $x''\in S$ and consider the terms in~\ref{ineq1:prop3} corresponding to state $x''$ with $v_{\alpha}$ replaced by $v^{n+1}_{\alpha}$ and note that they are bounded above by
\begin{multline*}
 \alpha \lambda_{i''} \Pi_{j=1}^\ell p^{a^*,j}_{s_j',s''_j} v^{n+1}_{\alpha}(x'') - \alpha \lambda_{i''} \Pi_{j=1}^\ell p^{a',j}_{s_j',s''_j} v_{\alpha}(x'') 
 - \alpha \lambda_{i''} \Pi_{j=1}^\ell p^{a^*,j}_{s_j,s''_j} v_{\alpha}(x'') + \alpha \lambda_{i''} \Pi_{j=1}^\ell p^{a',j}_{s_j,s''_j} v_{\alpha}(x''),
\end{multline*}
which can be rewritten as 
\begin{align*}
\alpha \lambda_{i''} v^{n}_{\alpha}(x'') \Big[
\Pi_{j=1}^\ell p^{a^*,j}_{s_j',s''_j} - \Pi_{j=1}^\ell p^{a',j}_{s_j',s''_j}
- \Pi_{j=1}^\ell p^{a^*,j}_{s_j,s''_j} + \Pi_{j=1}^\ell p^{a',j}_{s_j,s''_j} \Big].  
\end{align*}
The inductive hypothesis and assumption~\ref{ass:prop3} implies that this last expression is bounded above by $0$.  This yields the result for the discounted expected cost case.  The proof of the long-run average cost case is similar and is omitted for brevity. 

\subsection{Proof of Proposition~\ref{prop:complexity}}\label{appendix:1.5}
For the myopic policy $r$, for every possible state $x\in \X$, the algorithm selects the SNF in $A(x)$ with the lowest reward, which is done by checking all possible $|A(x)|$ cases. Therefore, the algorithm takes $M=\sum_{x\in \mathbb{X}}|A(x)|$ steps in total to construct the policy.

For the $r+pr$ policy, for every possible state $x\in \mathbb{X}$, the algorithm chooses the action $a\in A(x)$ with the lowest score given by
$$\gamma(i, a) := r^a_i + \sum_{x'=(i',s'_1,\ldots,s'_n)} \lambda_{i'} \Pi_{j=1}^\ell p^{a,j}_{s_j,s'_j} \min_{a' \in A(x')}\{r^{a'}_{i'}\}.$$
To compute the score, the algorithm checks every state $x'$, and every action $a'\in A(x')$, and multiplies it by a weight that has $l$ terms. Then, $\gamma(i, a)$ takes $\mathcal{O}(M\cdot l)$ steps to compute, and that has to be done $M$ times. Therefore, the policy $r+pr$ has complexity $\mathcal{O}(M^2\cdot l)$.

\subsection{Additional details for the `quantitative analysis' section} \label{appendix:2}

\begin{figure}
\centering
\includegraphics[height=2.3in]{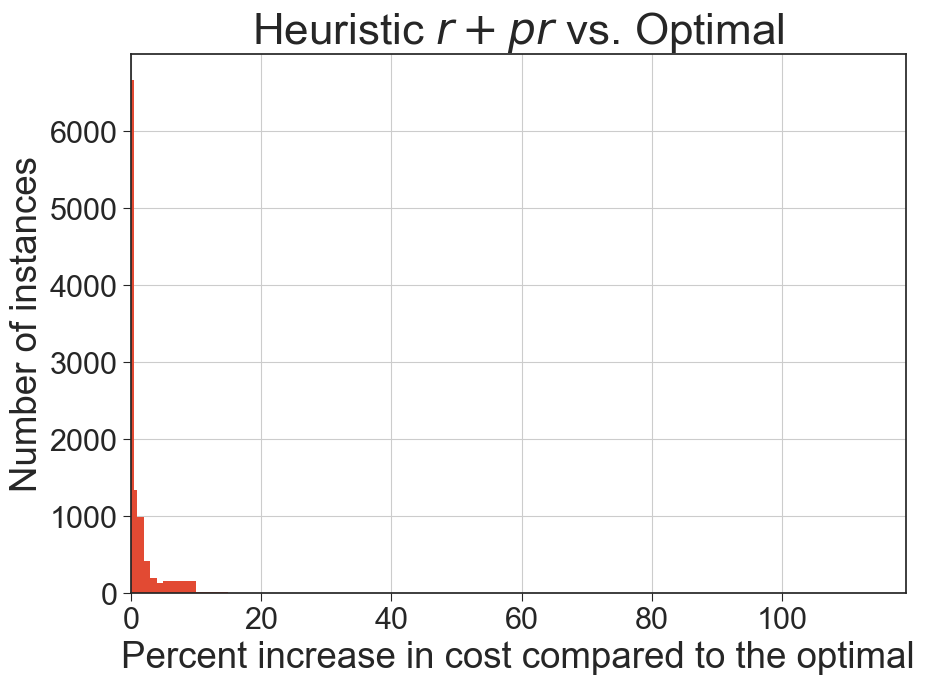}
\includegraphics[height=2.3in]{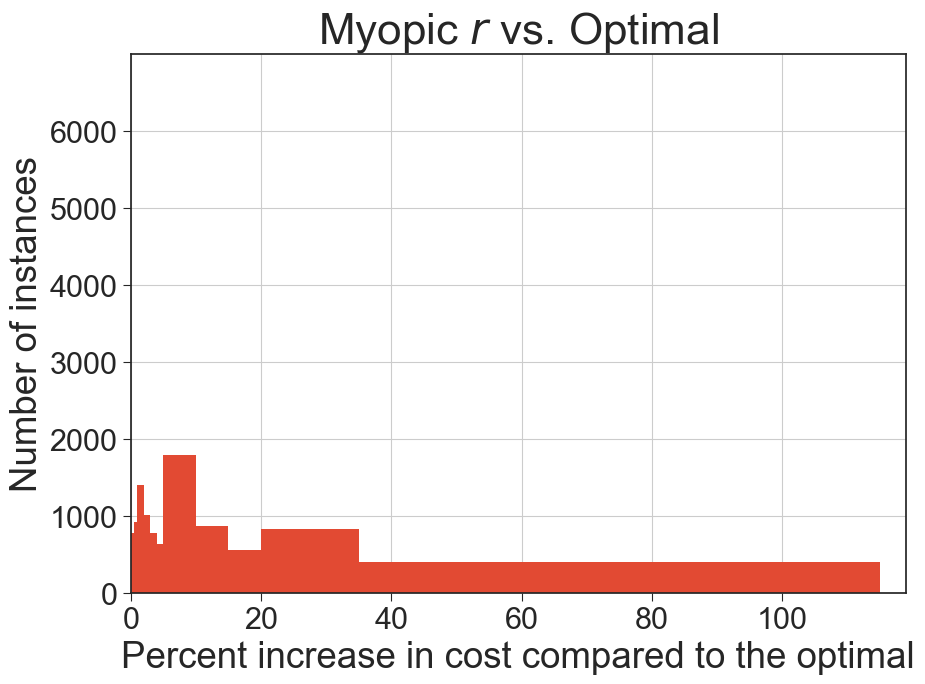} \\
\includegraphics[height=2.3in]{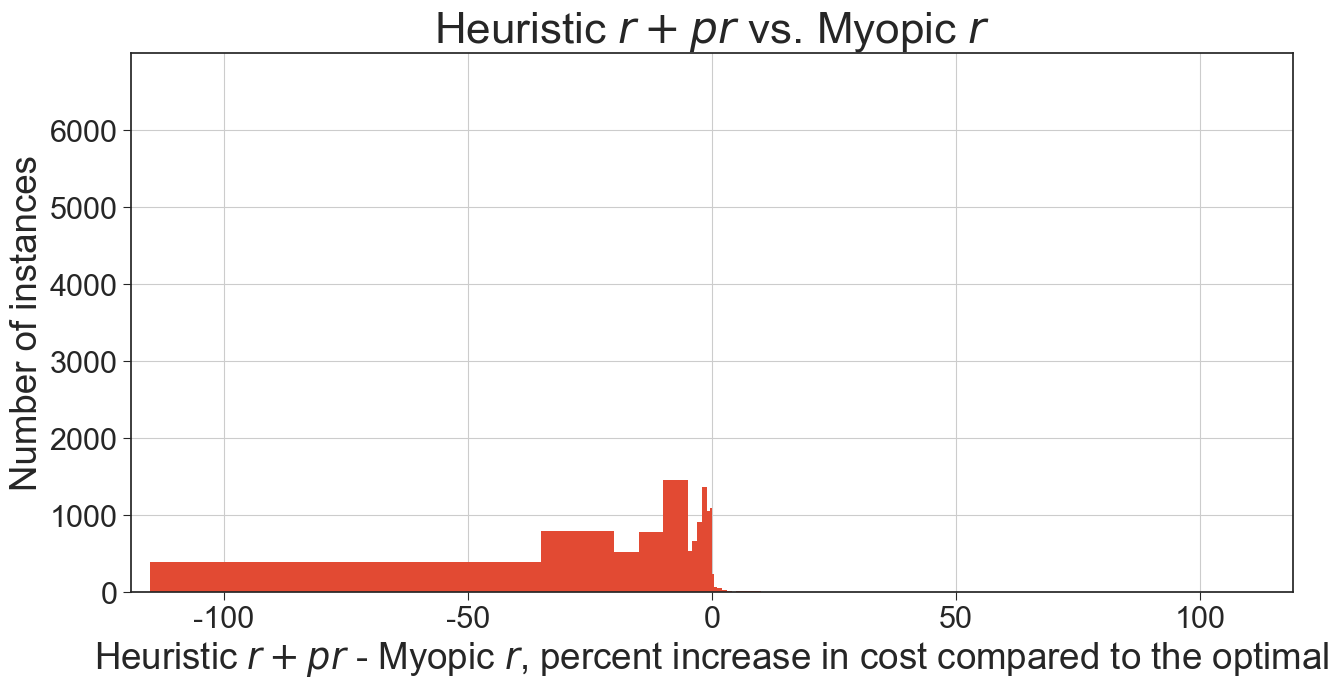}
\caption{Comparison of heuristics $r$ and $p+r$ under Scenario 2 with $\beta_2=0.2$,$\gamma_2=5$ in terms of percent increase in expected long-term average costs relative to optimal expected long-term average cost.} \label{fig:comparisonScen2}
\end{figure}

\begin{figure}
\centering
\includegraphics[height=2.3in]{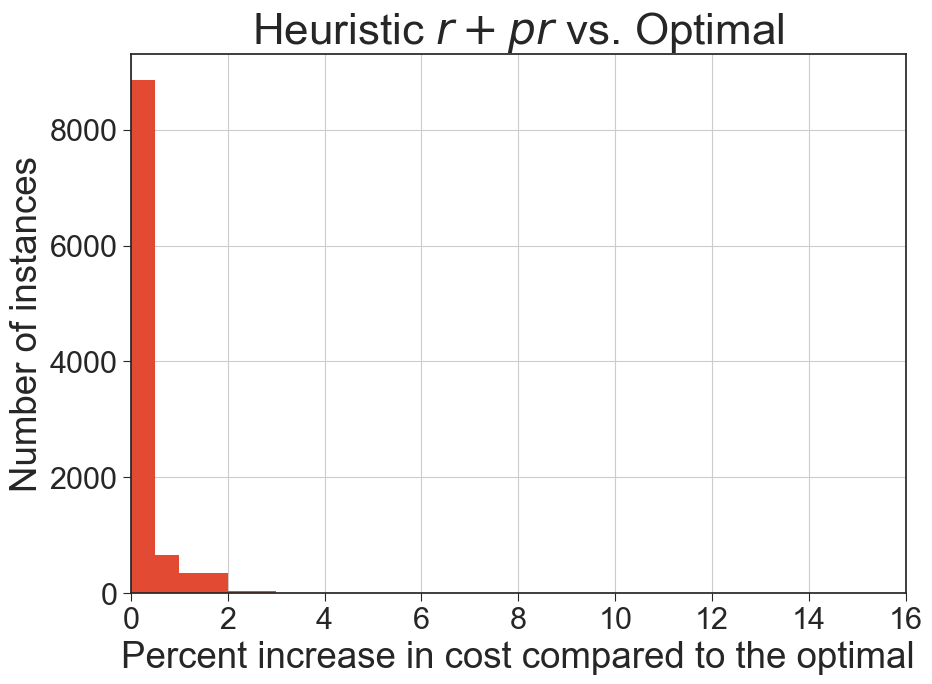}
\includegraphics[height=2.3in]{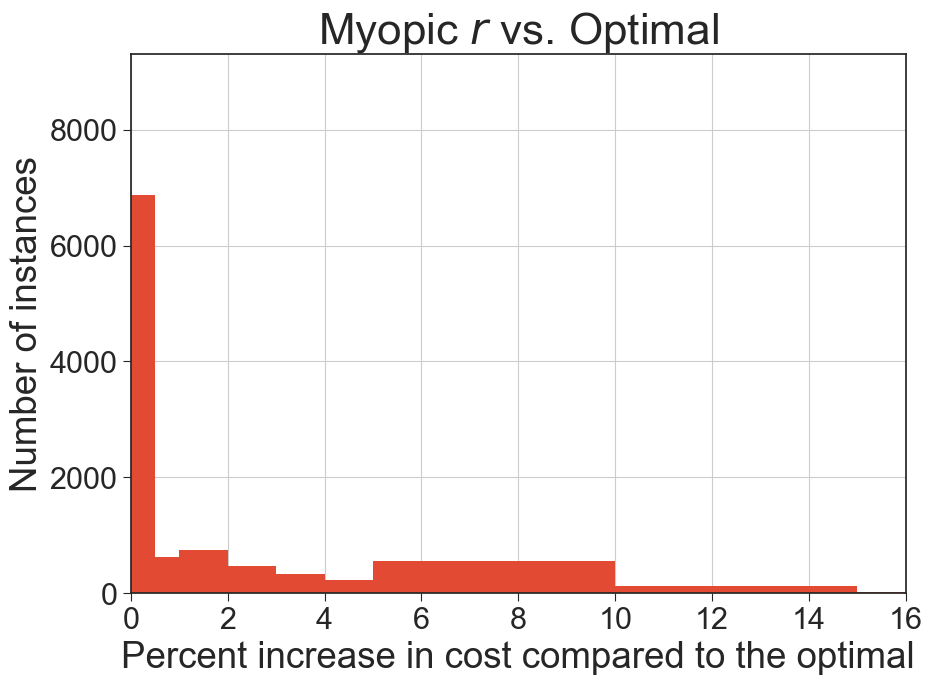} \\
\includegraphics[height=2.3in]{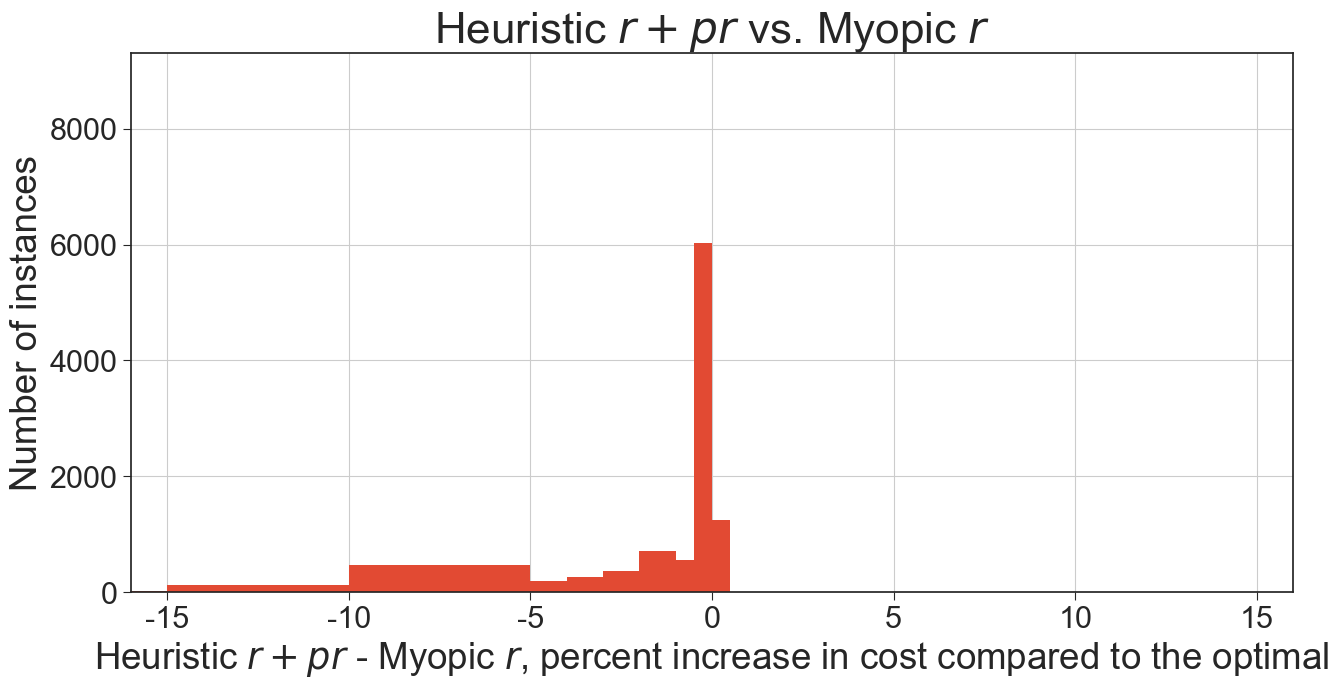}
\caption{Comparison of heuristics $r$ and $p+r$ under Scenario 2 with $\beta_3=0.2$, $ \gamma_3=5$ in terms of percent increase in expected long-term average costs relative to optimal policy.} \label{fig:comparisonScen3}
\end{figure}

\subsubsection{An extension to heuristic $r+pr$}

\begin{figure}
\centering
\includegraphics[height=2.3in]{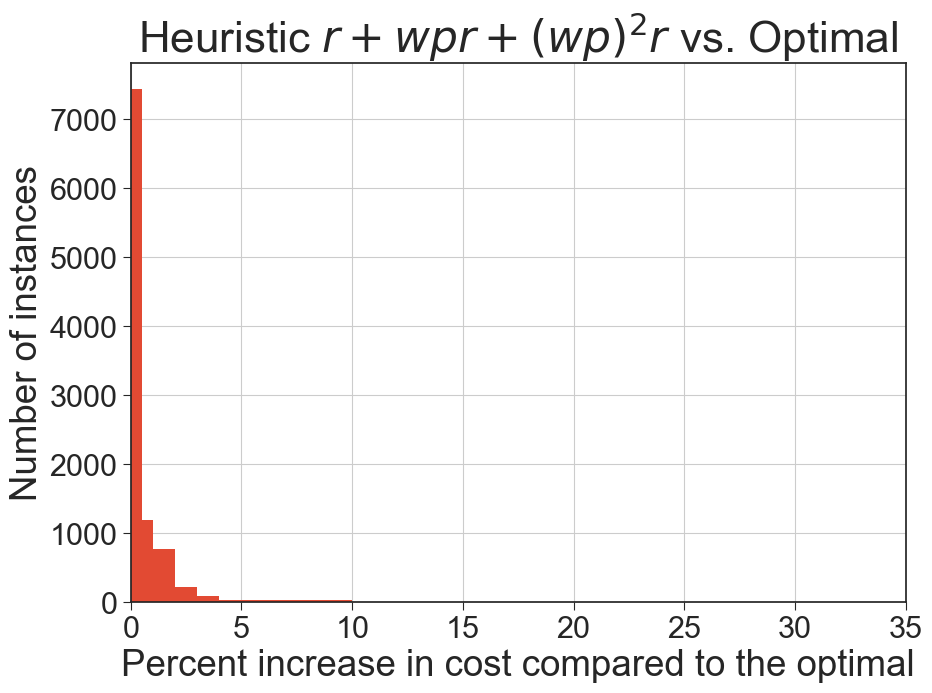}\\
\includegraphics[height=2.3in]{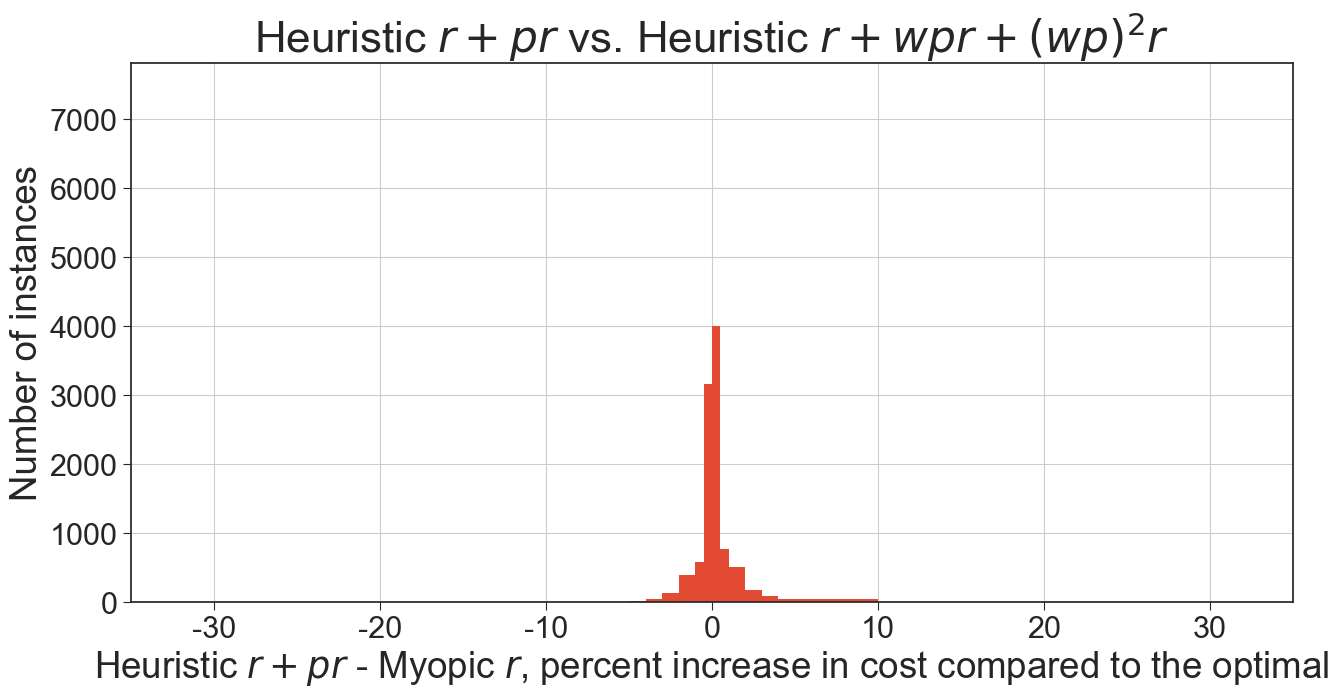} 
\caption{Comparison of heuristics $r+ wpr + w^2p^2r$ and $p+r$ under Scenario 1 with $\beta_1=0.2$in terms of percent increase in expected long-term average costs relative to optimal policy.} \label{fig:comparisonHeur3}
\end{figure}

We proposed and evaluated a natural extension of policy $r+p r$ that accounts for SNF availability two periods into the future instead of one, which we will label as $r + w p r + w^2 p^2 r$ with $w$ denoting a number between $0$ and $1$.  To evaluate this policy, we used the 100,000 instances of the 25 transition probability matrices generated by following the procedure described under Scenario 1 presented in Section~\ref{sec:numstudy} with $\beta_1 = 0.2$ only and tested weights $w=0.5, 0.75, 0.9,1$. The results of these experiments for $w=0.75$ only are summarized in Figure~\ref{fig:comparisonHeur3} where we plot the number of instances against percent increase in average readmission rates obtained by heuristic $r+ w p r + w^2 p^2 r$ compared to those obtained by the optimal policy (top panel) and the number of instances against the difference between percent increase in average readmission rates obtained by heuristic $r+pr$ compared to those obtained by the optimal policy minus percent increase in average readmission rates obtained by heuristic $r + w p + w^2 P^2 r$ compared to those obtained by the optimal policy (bottom panel).  Overall, we find that policy $r+ w p r + w^2 P^2 r$ outperformed heuristic $r + pr$ in roughly 55\% of the instances.

\subsection{Additional details for the `qualitative analysis' and `dependency structure' sections}

\subsubsection{First instance in which myopic policy $r$ performs well} \label{appendix:3}

For this instance, we started with baseline transition matrices:

\begin{align*}
& \mathbf{\hat{P}}^{A} = \begin{bmatrix}
0.04&0.96\\
0.11&0.89
\end{bmatrix},
\mathbf{\hat{P}}^{B} = \begin{bmatrix}
0.54&0.46\\
0.18&0.82
\end{bmatrix},
\mathbf{\hat{P}}^{C}= \begin{bmatrix}
0.93&0.07\\
0.35&0.65
\end{bmatrix},\\
& \mathbf{\hat{P}}^{D} = \begin{bmatrix}
0.88&0.12\\
0.49&0.51
\end{bmatrix},
\mathbf{\hat{P}}^{E} = \begin{bmatrix}
0.97&0.03\\
0.41&0.59
\end{bmatrix}.
\end{align*}

\medskip

From these baseline matrices, we generated the final 25 2-by-2 transition matrices under Scenario 1 with $\beta_1=0.2$:

\footnotesize
\begin{align*}
\mathbf{P}^{1,1} = \begin{bmatrix}
0.04&0.96\\
0.82&0.18
\end{bmatrix}, &\mathbf{P}^{2,1} = \begin{bmatrix}
0.04&0.96\\
0.11&0.89
\end{bmatrix}, &\mathbf{P}^{3,1} = \begin{bmatrix}
0.04&0.96\\
0.11&0.89
\end{bmatrix}, &\mathbf{P}^{4,1} = \begin{bmatrix}
0.04&0.96\\
0.11&0.89
\end{bmatrix}, &\mathbf{P}^{5,1} = \begin{bmatrix}
0.04&0.96\\
0.11&0.89
\end{bmatrix},\\
\mathbf{P}^{1,2} = \begin{bmatrix}
0.54&0.46\\
0.18&0.82
\end{bmatrix}, & \mathbf{P}^{2,2} = \begin{bmatrix}
0.54&0.46\\
0.84&0.16
\end{bmatrix}, & \mathbf{P}^{3,2} = \begin{bmatrix}
0.54&0.46\\
0.18&0.82
\end{bmatrix}, & \mathbf{P}^{4,2} = \begin{bmatrix}
0.54&0.46\\
0.18&0.82
\end{bmatrix}, & \mathbf{P}^{5,2} = \begin{bmatrix}
0.54&0.46\\
0.18&0.82
\end{bmatrix},\\
\mathbf{P}^{1,3} = \begin{bmatrix}
0.93&0.07\\
0.35&0.65
\end{bmatrix}, & \mathbf{P}^{2,3} = \begin{bmatrix}
0.93&0.07\\
0.35&0.65
\end{bmatrix}, & \mathbf{P}^{3,3} = \begin{bmatrix}
0.93&0.07\\
0.87&0.13
\end{bmatrix}, & \mathbf{P}^{4,3} = \begin{bmatrix}
0.93&0.07\\
0.35&0.65
\end{bmatrix}, & \mathbf{P}^{5,3} = \begin{bmatrix}
0.93&0.07\\
0.35&0.65
\end{bmatrix},\\
\mathbf{P}^{1,4} = \begin{bmatrix}
0.88&0.12\\
0.49&0.51
\end{bmatrix}, & \mathbf{P}^{2,4} = \begin{bmatrix}
0.88&0.12\\
0.49&0.51
\end{bmatrix}, & \mathbf{P}^{3,4} = \begin{bmatrix}
0.88&0.12\\
0.49&0.51
\end{bmatrix}, & \mathbf{P}^{4,4} = \begin{bmatrix}
0.88&0.12\\
0.9&0.1
\end{bmatrix}, & \mathbf{P}^{5,4} = \begin{bmatrix}
0.88&0.12\\
0.49&0.51
\end{bmatrix},\\
\mathbf{P}^{1,5} = \begin{bmatrix}
0.97&0.03\\
0.41&0.59
\end{bmatrix}, & \mathbf{P}^{2,5} = \begin{bmatrix}
0.97&0.03\\
0.41&0.59
\end{bmatrix}, & \mathbf{P}^{3,5} = \begin{bmatrix}
0.97&0.03\\
0.41&0.59
\end{bmatrix}, & \mathbf{P}^{4,5} = \begin{bmatrix}
0.97&0.03\\
0.41&0.59
\end{bmatrix}, & \mathbf{P}^{5,5} = \begin{bmatrix}
0.97&0.03\\
0.88&0.12
\end{bmatrix}.
\end{align*}

\medskip

\normalsize  Note, the matrices above are identical to the corresponding baseline transition matrices except for entry $p^{a,a}_{1,1}$ which is obtained by multiplying the corresponding entry in $\mathbf{\hat{P}}^a$ by $\beta_1$. Thus, since $\hat{p}^A_{1,1}=0.89$, then $p^{1,1}_{1,1}=0.89\times0.2=0.18$. Based on these transition matrices, we obtained each of the three policies, which are described by action and state in Table~\ref{tab:ExampleAppendix4}.

We also investigated the influence of discharge rates on our results using the same transition probabilities. We obtained qualitatively similar long-run average costs when changing discharge rates to $\lambda_{UM}=0.4$, $\lambda_{CM}=0.3$, $\lambda_{CS}=0.2$, and $\lambda_{JS}=0.075$ and when the discharge rates are $\lambda_{UM}=0.025$, $\lambda_{CM}=0.075$, $\lambda_{CS}=0.2$, and $\lambda_{JS}=0.3$.  In the former case, the optimal average cost is approximately $7.91$, the average cost of the myopic policy is approximately $7.91$, and the average cost of heuristic $r + p r$ is approximately $8.09$.  It follows that the average costs for the myopic and $r+pr$ heuristics, respectively, are 0\% and 2.28\% higher than the optimal average cost.  In the latter case, the optimal average cost is approximately $22.36$, the average cost of the myopic policy is approximately $22.41$, and the average cost of heuristic $r + p r$ is approximately $23.6$.  It follows that the average costs of the myopic and $r+pr$ heuristics, respectively, are 0.22\% and 5.55\% higher than the optimal average cost.

\begin{figure}
\centering
\includegraphics[width=0.35\textwidth]{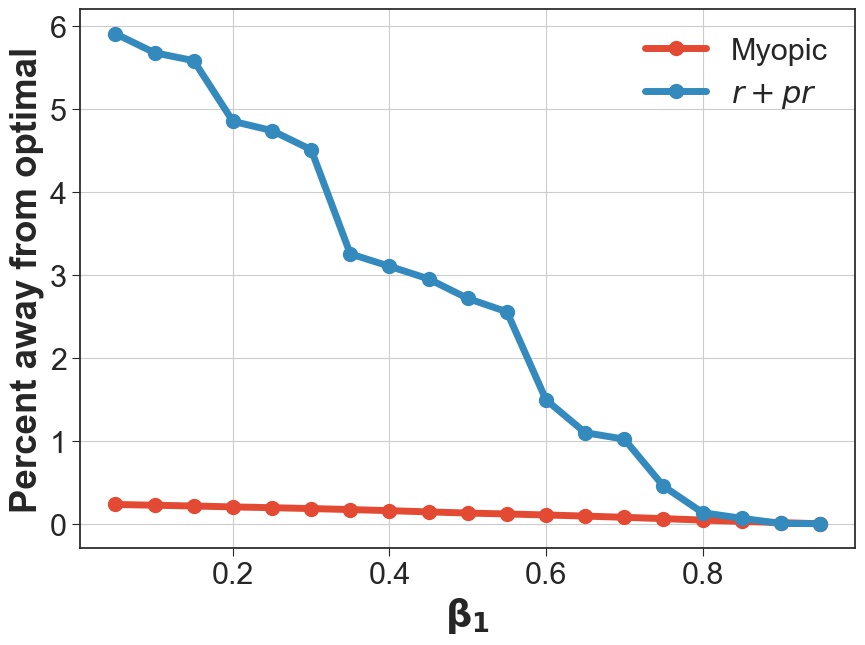} \\
\includegraphics[width=0.65\textwidth]{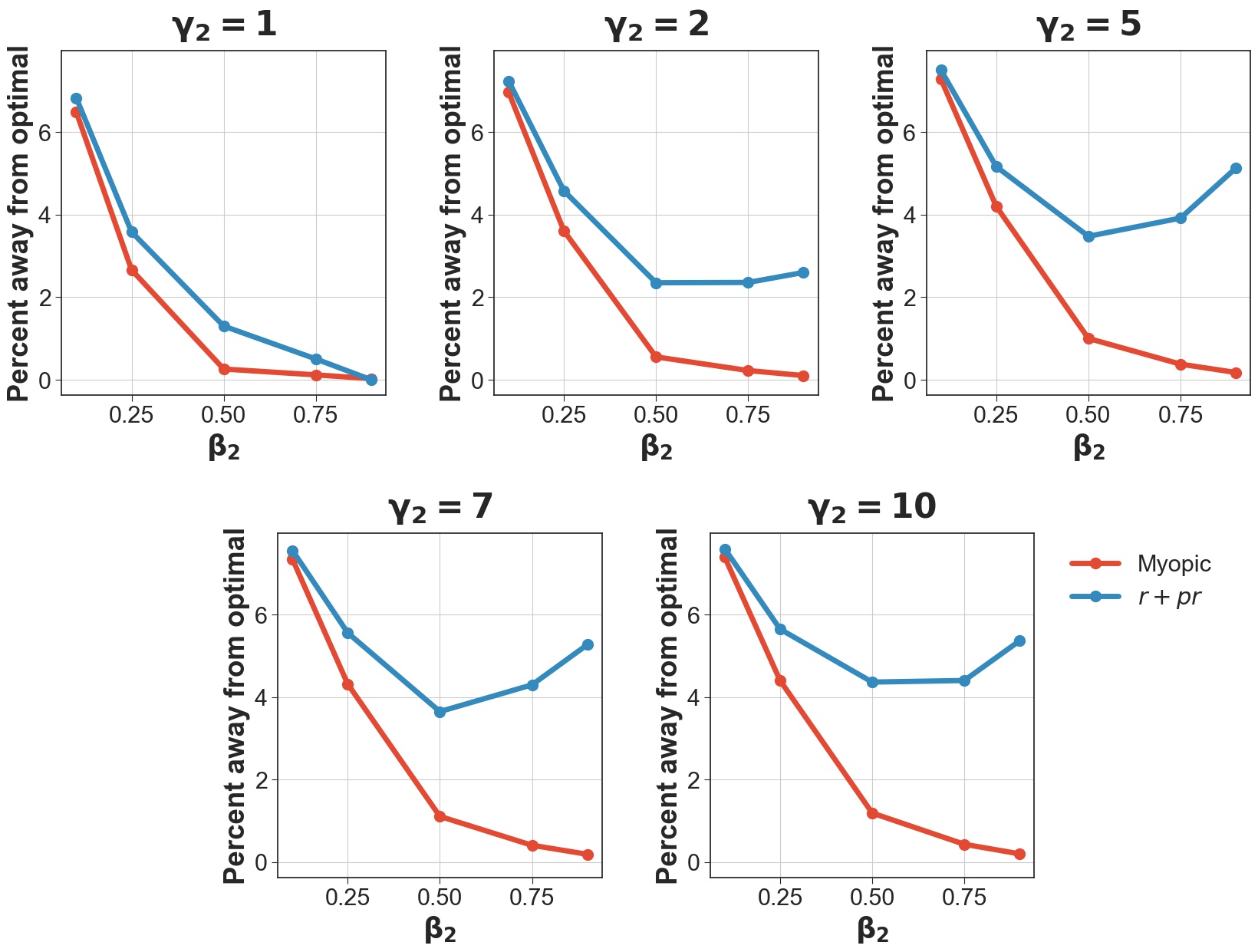} \\
\includegraphics[width=0.8\textwidth]{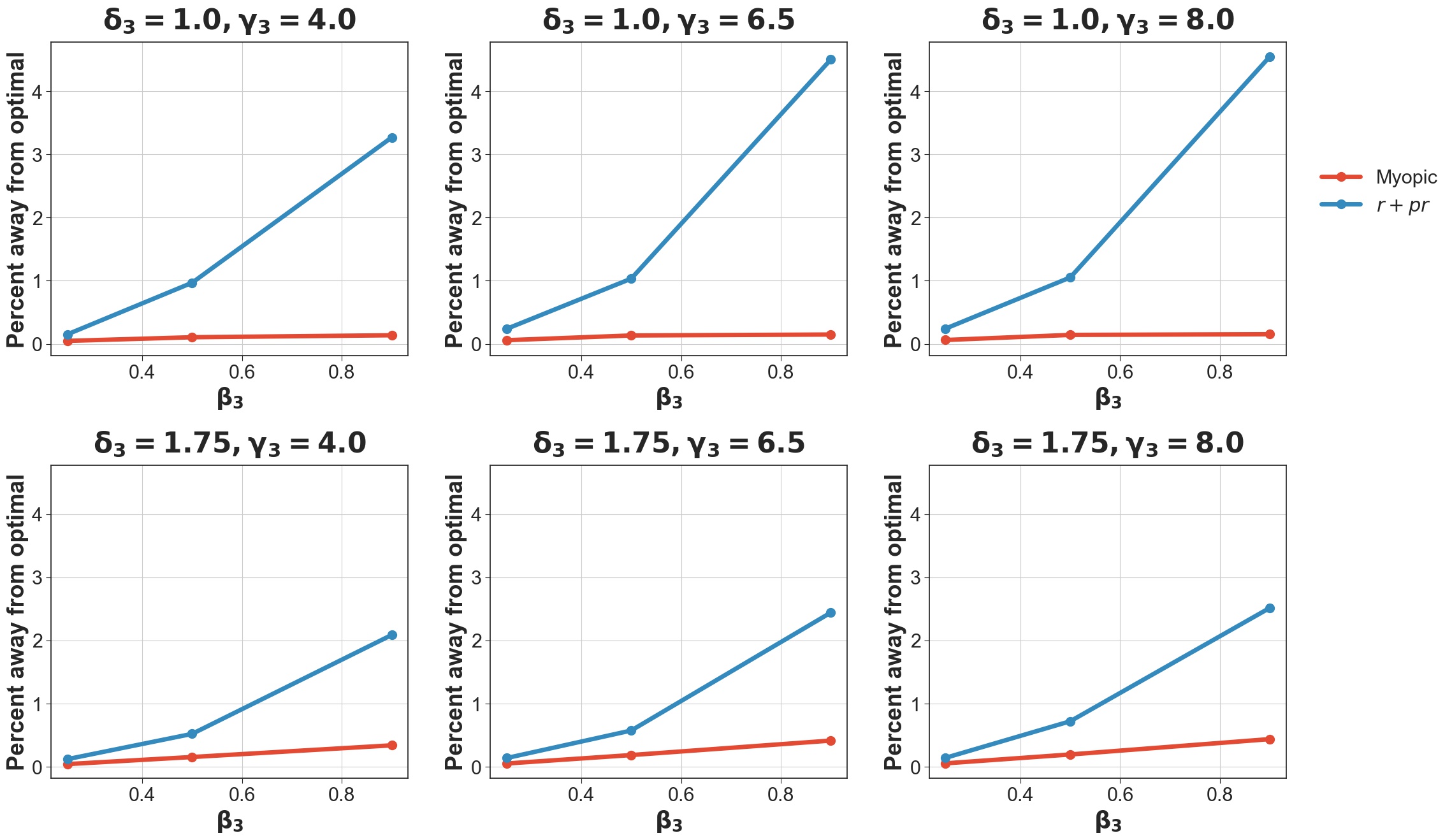}
\caption{For first instance considered when myopic performs well, the relative gap (in \%) of the expected long-run average costs between each heuristic (myopic $r$ and heuristic $r+pr$) under Scenarios 1--3.} \label{fig:instance1_myopic}
\end{figure}

We also used the same baseline transition matrices in this instance to investigate Scenarios 1--3 for all the parameters listed in Table~\ref{tab:dependency_parameters}. Results from this parameter study are summarized in Figure~\ref{fig:instance1_myopic}.

\subsubsection{Second instance in which myopic policy $r$ performs well} \label{appendix:4}

For this instance, we started with baseline transition matrices:

\begin{align*}
& \mathbf{\hat{P}}^{A} = \begin{bmatrix}
0.86&0.14\\
0.3&0.7
\end{bmatrix},
\mathbf{\hat{P}}^{B} = \begin{bmatrix}
0.93&0.07\\
0.05&0.95
\end{bmatrix},
\mathbf{\hat{P}}^{C} = \begin{bmatrix}
0.16&0.84\\
0.01&0.99
\end{bmatrix},\\
& \mathbf{\hat{P}}^{D} = \begin{bmatrix}
0.92&0.08\\
0.75&0.25
\end{bmatrix},
\mathbf{\hat{P}}^{E} = \begin{bmatrix}
0.8&0.2\\
0.97&0.03
\end{bmatrix}.
\end{align*}

\medskip

From these baseline matrices, we generated the final 25 2-by-2 transition matrices under Scenario 1 with $\beta_1=0.2$:

\footnotesize
\begin{align*}
\mathbf{P}^{1,1} = \begin{bmatrix}
0.86&0.14\\
0.86&0.14
\end{bmatrix}, &\mathbf{P}^{2,1} = \begin{bmatrix}
0.86&0.14\\
0.3&0.7
\end{bmatrix}, &\mathbf{P}^{3,1} = \begin{bmatrix}
0.86&0.14\\
0.3&0.7
\end{bmatrix}, &\mathbf{P}^{4,1} = \begin{bmatrix}
0.86&0.14\\
0.3&0.7
\end{bmatrix}, &\mathbf{P}^{5,1} = \begin{bmatrix}
0.86&0.14\\
0.3&0.7
\end{bmatrix},\\
\mathbf{P}^{1,2} = \begin{bmatrix}
0.93&0.07\\
0.05&0.95
\end{bmatrix}, &\mathbf{P}^{2,2} = \begin{bmatrix}
0.93&0.07\\
0.81&0.19
\end{bmatrix}, &\mathbf{P}^{3,2} = \begin{bmatrix}
0.93&0.07\\
0.05&0.95
\end{bmatrix}, &\mathbf{P}^{4,2} = \begin{bmatrix}
0.93&0.07\\
0.05&0.95
\end{bmatrix}, &\mathbf{P}^{5,2} = \begin{bmatrix}
0.93&0.07\\
0.05&0.95
\end{bmatrix},\\
\mathbf{P}^{1,3} = \begin{bmatrix}
0.16&0.84\\
0.01&0.99
\end{bmatrix}, &\mathbf{P}^{2,3} = \begin{bmatrix}
0.16&0.84\\
0.01&0.99
\end{bmatrix}, &\mathbf{P}^{3,3} = \begin{bmatrix}
0.16&0.84\\
0.8&0.2
\end{bmatrix}, &\mathbf{P}^{4,3} = \begin{bmatrix}
0.16&0.84\\
0.01&0.99
\end{bmatrix}, &\mathbf{P}^{5,3} = \begin{bmatrix}
0.16&0.84\\
0.01&0.99
\end{bmatrix},\\
\mathbf{P}^{1,4} = \begin{bmatrix}
0.92&0.08\\
0.75&0.25
\end{bmatrix}, &\mathbf{P}^{2,4} = \begin{bmatrix}
0.92&0.08\\
0.75&0.25
\end{bmatrix}, &\mathbf{P}^{3,4} = \begin{bmatrix}
0.92&0.08\\
0.75&0.25
\end{bmatrix}, &\mathbf{P}^{4,4} = \begin{bmatrix}
0.92&0.08\\
0.95&0.05
\end{bmatrix}, &\mathbf{P}^{5,4} = \begin{bmatrix}
0.92&0.08\\
0.75&0.25
\end{bmatrix},\\
\mathbf{P}^{1,5} = \begin{bmatrix}
0.8&0.2\\
0.97&0.03
\end{bmatrix}, &\mathbf{P}^{2,5} = \begin{bmatrix}
0.8&0.2\\
0.97&0.03
\end{bmatrix}, &\mathbf{P}^{3,5} = \begin{bmatrix}
0.8&0.2\\
0.97&0.03
\end{bmatrix}, &\mathbf{P}^{4,5} = \begin{bmatrix}
0.8&0.2\\
0.97&0.03
\end{bmatrix}, &\mathbf{P}^{5,5} = \begin{bmatrix}
0.8&0.2\\
0.99&0.01
\end{bmatrix}.
\end{align*}

\medskip

\normalsize  Note that these matrices above are identical to the corresponding baseline transition matrices except for entry $p^{a,a}_{1,1}$ which is obtained by multiplying the corresponding entry in $\mathbf{\hat{P}}^a$ by $\beta_1$. Thus, since $\hat{p}^A_{1,1}=0.14$, then $p^{1,1}_{1,}=0.7\times0.2=0.14$. Based on these transition matrices, we obtained each of the three policies, which are described by action and state in Table~\ref{tab:ExampleAppendix5}.


We also investigated the influence of discharge rates on our results using the same transition probabilities. We recovered qualitatively similar long-run average costs when changing discharge rates to $\lambda_{UM}=0.4$, $\lambda_{CM}=0.3$, $\lambda_{CS}=0.2$, and $\lambda_{JS}=0.075$ and when the discharge rates are $\lambda_{UM}=0.025$, $\lambda_{CM}=0.075$, $\lambda_{CS}=0.2$, and $\lambda_{JS}=0.3$.  In the former case, the optimal average cost is approximately $7.96$, the average cost of the myopic policy is approximately $8.07$, and the average cost of heuristic $r + p r$ is approximately $8.47$.  It follows that the average costs of the myopic and $r+pr$ heuristics, respectively, are 1.38\% and 6.41\% higher than the optimal average cost.  In the latter case, the optimal average cost is approximately $22.71$, the average cost of the myopic policy is approximately $23.43$, and the average cost of heuristic $r + p r$ is approximately $24.64$.  It follows that the average costs of the myopic and $r+pr$ heuristics, respectively, are 3.17\% and 8.5\% higher than the optimal average cost.

\begin{figure}
\centering
\includegraphics[width=0.3\textwidth]{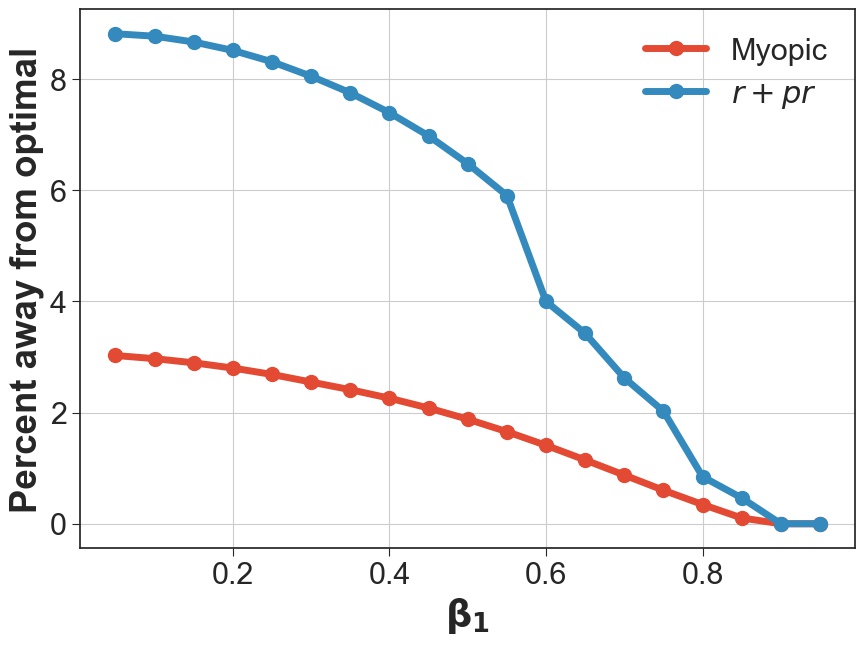} \\
\includegraphics[width=0.6\textwidth]{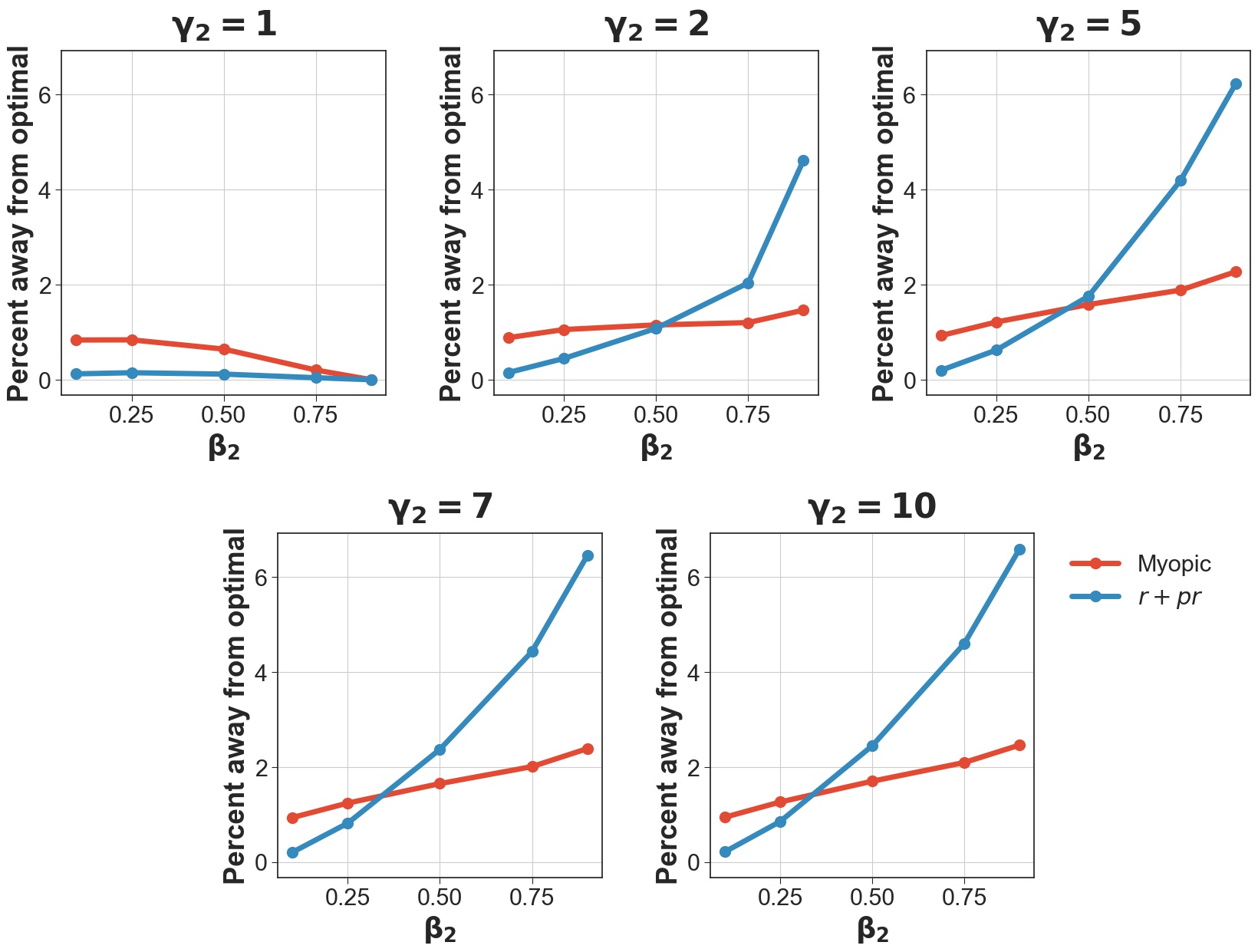} \\
\includegraphics[width=0.5\textwidth]{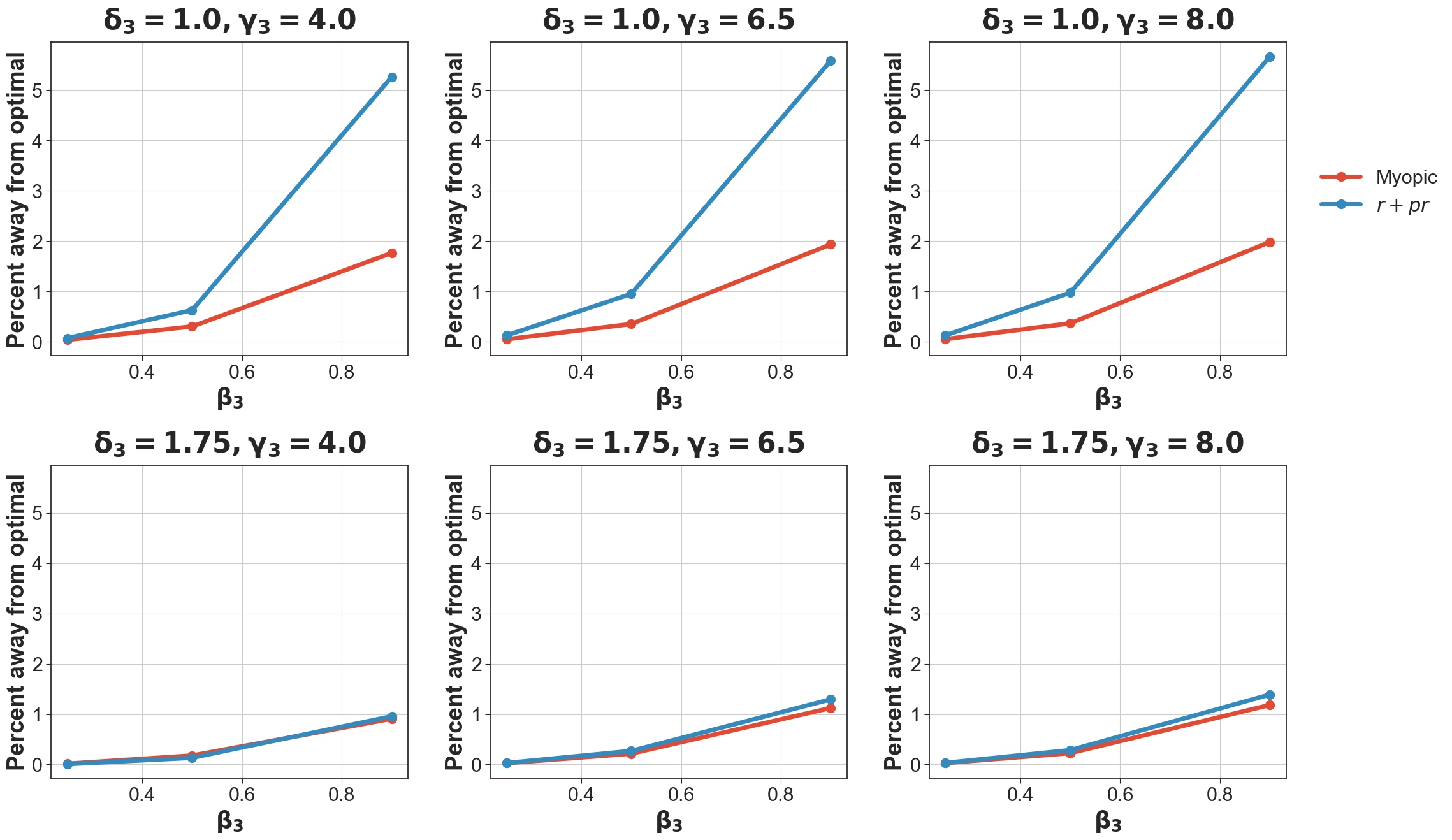}
\caption{For second instance considered when myopic performs well, the relative gap (in \%) of the expected long-run average costs between each heuristic (myopic $r$ and heuristic $r+pr$) under Scenarios 1--3.} \label{fig:instance2_myopic}
\end{figure}

We also used the same baseline transition matrices in this instance to investigate Scenarios 1--3 for all the parameters listed in Table~\ref{tab:dependency_parameters}. Results from this parameter study are summarized in Figure~\ref{fig:instance2_myopic}.

\subsubsection{First instance in which policy $r + p r$ performs well} \label{appendix:5}

For this instance, we started with baseline transition matrices:

\begin{align*}
& \mathbf{\hat{P}}^{A} = \begin{bmatrix}
0.01&0.99\\
0.97&0.03
\end{bmatrix},
\mathbf{\hat{P}}^{B} = \begin{bmatrix}
0.67&0.33\\
0.97&0.03
\end{bmatrix},
\mathbf{\hat{P}}^{C} = \begin{bmatrix}
0.83&0.17\\
0.84&0.16
\end{bmatrix},\\
& \mathbf{\hat{P}}^{D} = \begin{bmatrix}
0.11&0.89\\
0.13&0.87
\end{bmatrix},
\mathbf{\hat{P}}^{E} = \begin{bmatrix}
0.72&0.28\\
0.97&0.03
\end{bmatrix}.
\end{align*}

\medskip
From these baseline matrices, we generated the final 25 2-by-2 transition matrices under Scenario 1 with $\beta_1=0.2$:

\footnotesize 
\begin{align*}
\mathbf{P}^{1,1} = \begin{bmatrix}
0.01&0.99\\
0.99&0.01
\end{bmatrix}, &\mathbf{P}^{2,1} = \begin{bmatrix}
0.01&0.99\\
0.97&0.03
\end{bmatrix}, &\mathbf{P}^{3,1} = \begin{bmatrix}
0.01&0.99\\
0.97&0.03
\end{bmatrix}, &\mathbf{P}^{4,1} = \begin{bmatrix}
0.01&0.99\\
0.97&0.03
\end{bmatrix}, &\mathbf{P}^{5,1} = \begin{bmatrix}
0.01&0.99\\
0.97&0.03
\end{bmatrix},\\
\mathbf{P}^{1,2} = \begin{bmatrix}
0.67&0.33\\
0.97&0.03
\end{bmatrix}, &\mathbf{P}^{2,2} = \begin{bmatrix}
0.67&0.33\\
0.99&0.01
\end{bmatrix}, &\mathbf{P}^{3,2} = \begin{bmatrix}
0.67&0.33\\
0.97&0.03
\end{bmatrix}, &\mathbf{P}^{4,2} = \begin{bmatrix}
0.67&0.33\\
0.97&0.03
\end{bmatrix}, &\mathbf{P}^{5,2} = \begin{bmatrix}
0.67&0.33\\
0.97&0.03
\end{bmatrix},\\
\mathbf{P}^{1,3} = \begin{bmatrix}
0.83&0.17\\
0.84&0.16
\end{bmatrix}, &\mathbf{P}^{2,3} = \begin{bmatrix}
0.83&0.17\\
0.84&0.16
\end{bmatrix}, &\mathbf{P}^{3,3} = \begin{bmatrix}
0.83&0.17\\
0.97&0.03
\end{bmatrix}, &\mathbf{P}^{4,3} = \begin{bmatrix}
0.83&0.17\\
0.84&0.16
\end{bmatrix}, &\mathbf{P}^{5,3} = \begin{bmatrix}
0.83&0.17\\
0.84&0.16
\end{bmatrix},\\
\mathbf{P}^{1,4} = \begin{bmatrix}
0.11&0.89\\
0.13&0.87
\end{bmatrix}, &\mathbf{P}^{2,4} = \begin{bmatrix}
0.11&0.89\\
0.13&0.87
\end{bmatrix}, &\mathbf{P}^{3,4} = \begin{bmatrix}
0.11&0.89\\
0.13&0.87
\end{bmatrix}, &\mathbf{P}^{4,4} = \begin{bmatrix}
0.11&0.89\\
0.83&0.17
\end{bmatrix}, &\mathbf{P}^{5,4} = \begin{bmatrix}
0.11&0.89\\
0.13&0.87
\end{bmatrix},\\
\mathbf{P}^{1,5} = \begin{bmatrix}
0.72&0.28\\
0.97&0.03
\end{bmatrix}, &\mathbf{P}^{2,5} = \begin{bmatrix}
0.72&0.28\\
0.97&0.03
\end{bmatrix}, &\mathbf{P}^{3,5} = \begin{bmatrix}
0.72&0.28\\
0.97&0.03
\end{bmatrix}, &\mathbf{P}^{4,5} = \begin{bmatrix}
0.72&0.28\\
0.97&0.03
\end{bmatrix}, &\mathbf{P}^{5,5} = \begin{bmatrix}
0.72&0.28\\
0.99&0.01
\end{bmatrix}.
\end{align*}

\normalsize  Note, these matrices above are identical to the corresponding baseline transition matrices except for entry $p^{a,a}_{1,1}$ which we obtained by multiplying the corresponding entry in $\mathbf{\hat{P}}^a$ by $\beta_1$. Thus, since $\hat{p}^A_{1,1}=0.03$, then $p^{1,1}_{1,1}=0.03\times0.2=0.01$. Based on these transition matrices, we recovered each of the three policies, which are described by action and state in Table~\ref{tab:ExampleAppendix6}.

We also investigated the influence of discharge rates on our results using the same transition probabilities. We recovered qualitatively similar long-run average costs when changing discharge rates to $\lambda_{UM}=0.4$, $\lambda_{CM}=0.3$, $\lambda_{CS}=0.2$, and $\lambda_{JS}=0.075$ and when the discharge rates are $\lambda_{UM}=0.025$, $\lambda_{CM}=0.075$, $\lambda_{CS}=0.2$, and $\lambda_{JS}=0.3$.  In the former case, the optimal average cost is approximately $8.7$, the average cost of the myopic policy is approximately $13.37$, and the average cost of heuristic $r + p r$ is approximately $8.7$.  It follows that the average costs of the myopic and $r+pr$ heuristics, respectively, are 53.68\% and 0\% higher than the optimal average cost.  In the latter case, the optimal average cost is approximately $17.08$, the average cost of the myopic policy is approximately $20.31$, and the average cost of heuristic $r + p r$ is approximately $17.12$.  It follows that the average costs of the myopic and $r+pr$ heuristics, respectively, are 18.9\% and 0.23\% higher than the optimal average cost.

\begin{figure}
\centering
\includegraphics[width=0.3\textwidth]{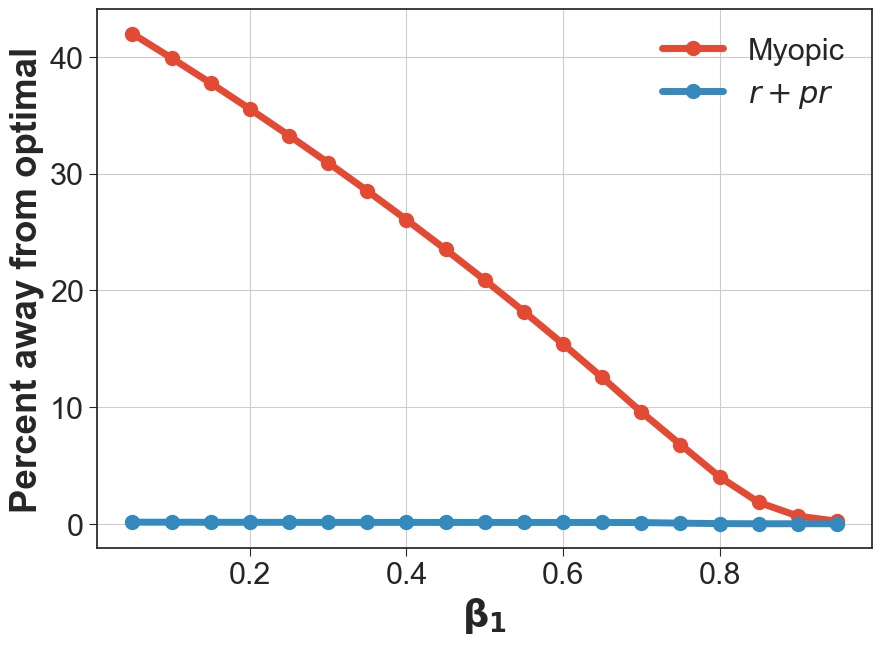} \\
\includegraphics[width=0.6\textwidth]{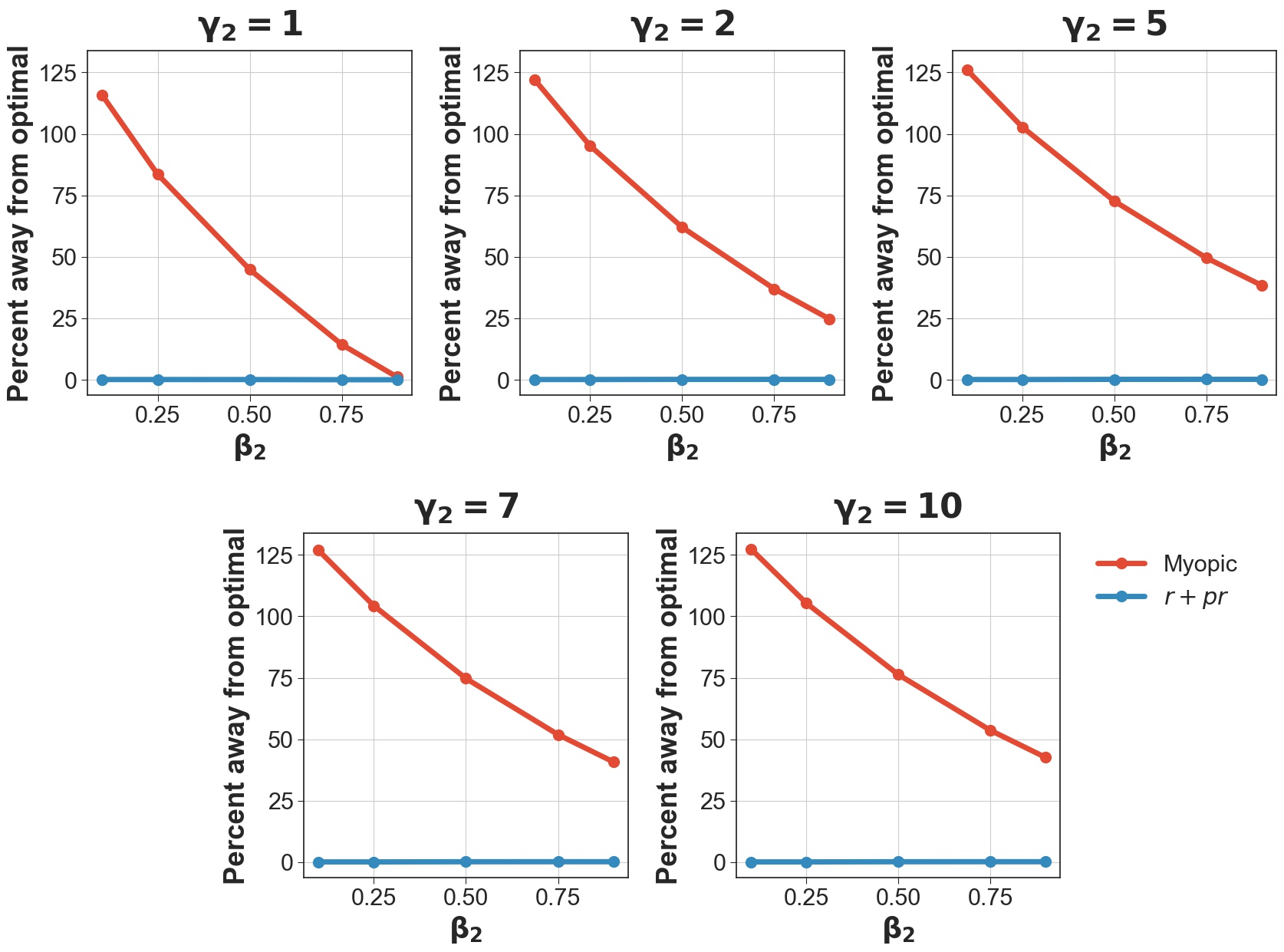} \\
\includegraphics[width=0.5\textwidth]{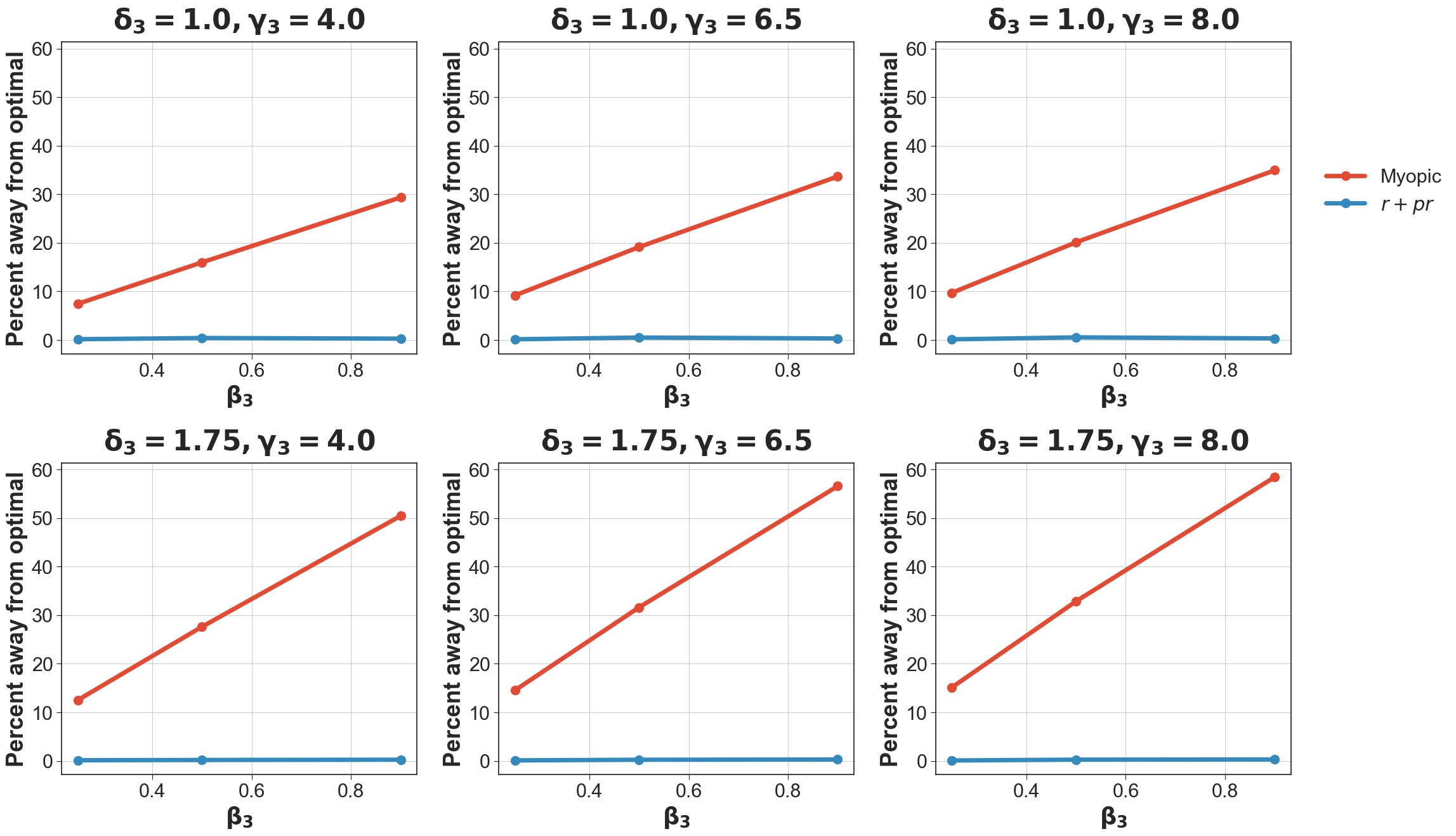}
\caption{For first instance considered when heuristic $r+p r$ performs well, the relative gap (in \%) of the expected long-run average costs between each heuristic (myopic $r$ and heuristic $r+pr$) under Scenarios 1--3.} \label{fig:instance3_myopic}
\end{figure}

We also used the same baseline transition matrices in this instance to investigate Scenarios 1--3 for all the parameters listed in Table~\ref{tab:dependency_parameters}. Results from this parameter study are summarized in Figure~\ref{fig:instance3_myopic}.

\subsubsection{Second instance considered when policy $r + p r$ performs well} \label{appendix:6}

For this instance, we started with baseline transition matrices:

\begin{align*}
& \mathbf{\hat{P}}^{A} = \begin{bmatrix}
0.41&0.59\\
0.69&0.31
\end{bmatrix},
\mathbf{\hat{P}}^{B} = \begin{bmatrix}
0.05&0.95\\
0.64&0.36
\end{bmatrix},
\mathbf{\hat{P}}^{C} = \begin{bmatrix}
0.91&0.09\\
0.66&0.34
\end{bmatrix},\\
& \mathbf{\hat{P}}^{D} = \begin{bmatrix}
0.16&0.84\\
0.89&0.11
\end{bmatrix},
\mathbf{\hat{P}}^{E} = \begin{bmatrix}
0.97&0.03\\
0.01&0.99
\end{bmatrix}.
\end{align*}

\medskip

From these baseline matrices, we generated the final 25 2-by-2 transition matrices under Scenario 1 with $\beta_1=0.2$:

\footnotesize
\begin{align*}
\mathbf{P}^{1,1} = \begin{bmatrix}
0.41&0.59\\
0.94&0.06
\end{bmatrix}, &\mathbf{P}^{2,1} = \begin{bmatrix}
0.41&0.59\\
0.69&0.31
\end{bmatrix}, &\mathbf{P}^{3,1} = \begin{bmatrix}
0.41&0.59\\
0.69&0.31
\end{bmatrix}, &\mathbf{P}^{4,1} = \begin{bmatrix}
0.41&0.59\\
0.69&0.31
\end{bmatrix}, &\mathbf{P}^{5,1} = \begin{bmatrix}
0.41&0.59\\
0.69&0.31
\end{bmatrix},\\
\mathbf{P}^{1,2} = \begin{bmatrix}
0.05&0.95\\
0.64&0.36
\end{bmatrix}, &\mathbf{P}^{2,2} = \begin{bmatrix}
0.05&0.95\\
0.93&0.07
\end{bmatrix}, &\mathbf{P}^{3,2} = \begin{bmatrix}
0.05&0.95\\
0.64&0.36
\end{bmatrix}, &\mathbf{P}^{4,2} = \begin{bmatrix}
0.05&0.95\\
0.64&0.36
\end{bmatrix}, &\mathbf{P}^{5,2} = \begin{bmatrix}
0.05&0.95\\
0.64&0.36
\end{bmatrix},\\
\mathbf{P}^{1,3} = \begin{bmatrix}
0.91&0.09\\
0.66&0.34
\end{bmatrix}, &\mathbf{P}^{2,3} = \begin{bmatrix}
0.91&0.09\\
0.66&0.34
\end{bmatrix}, &\mathbf{P}^{3,3} = \begin{bmatrix}
0.91&0.09\\
0.93&0.07
\end{bmatrix}, &\mathbf{P}^{4,3} = \begin{bmatrix}
0.91&0.09\\
0.66&0.34
\end{bmatrix}, &\mathbf{P}^{5,3} = \begin{bmatrix}
0.91&0.09\\
0.66&0.34
\end{bmatrix},\\
\mathbf{P}^{1,4} = \begin{bmatrix}
0.16&0.84\\
0.89&0.11
\end{bmatrix}, &\mathbf{P}^{2,4} = \begin{bmatrix}
0.16&0.84\\
0.89&0.11
\end{bmatrix}, &\mathbf{P}^{3,4} = \begin{bmatrix}
0.16&0.84\\
0.89&0.11
\end{bmatrix}, &\mathbf{P}^{4,4} = \begin{bmatrix}
0.16&0.84\\
0.98&0.02
\end{bmatrix}, &\mathbf{P}^{5,4} = \begin{bmatrix}
0.16&0.84\\
0.89&0.11
\end{bmatrix},\\
\mathbf{P}^{1,5} = \begin{bmatrix}
0.97&0.03\\
0.01&0.99
\end{bmatrix}, &\mathbf{P}^{2,5} = \begin{bmatrix}
0.97&0.03\\
0.01&0.99
\end{bmatrix}, &\mathbf{P}^{3,5} = \begin{bmatrix}
0.97&0.03\\
0.01&0.99
\end{bmatrix}, &\mathbf{P}^{4,5} = \begin{bmatrix}
0.97&0.03\\
0.01&0.99
\end{bmatrix}, &\mathbf{P}^{5,5} = \begin{bmatrix}
0.97&0.03\\
0.8&0.2
\end{bmatrix}.
\end{align*}

\medskip

\normalsize Note that these matrices above are identical to the corresponding baseline transition matrices except for entry $p^{a,a}_{1,1}$ which we obtained by multiplying the corresponding entry in $\mathbf{\hat{P}}^a$ by $\beta_1$. Thus, since $\hat{p}^A_{1,1}=0.31$, then $p^{1,1}_{1,}=0.31\times0.2=0.06$. Based on these transition matrices, we recovered each of the three policies, which are described by action and state in Table~\ref{tab:ExampleAppendix7}.

We also investigated the influence of discharge rates on our results using the same transition probabilities. We recovered qualitatively similar long-run average costs when changing discharge rates to $\lambda_{UM}=0.4$, $\lambda_{CM}=0.3$, $\lambda_{CS}=0.2$, and $\lambda_{JS}=0.075$ and when the discharge rates are $\lambda_{UM}=0.025$, $\lambda_{CM}=0.075$, $\lambda_{CS}=0.2$, and $\lambda_{JS}=0.3$.  In the former case, the optimal average cost is approximately $7.45$, the average cost of the myopic policy is approximately $8.44$, and the average cost of heuristic $r + p r$ is approximately $7.53$.  It follows that the average costs of the myopic and $r+pr$ heuristics, respectively, are 13.2\% and 1.01\% higher than the optimal average cost.  In the latter case, the optimal average cost is approximately $17.9$, the average cost of the myopic policy is approximately $20.98$, and the average cost of heuristic $r + p r$ is approximately $19.5$.  It follows that the average costs of the myopic and $r+pr$ heuristics, respectively, are 17.21\% and 8.94\% higher than the optimal average cost.

\begin{figure}
\centering
\includegraphics[width=0.3\textwidth]{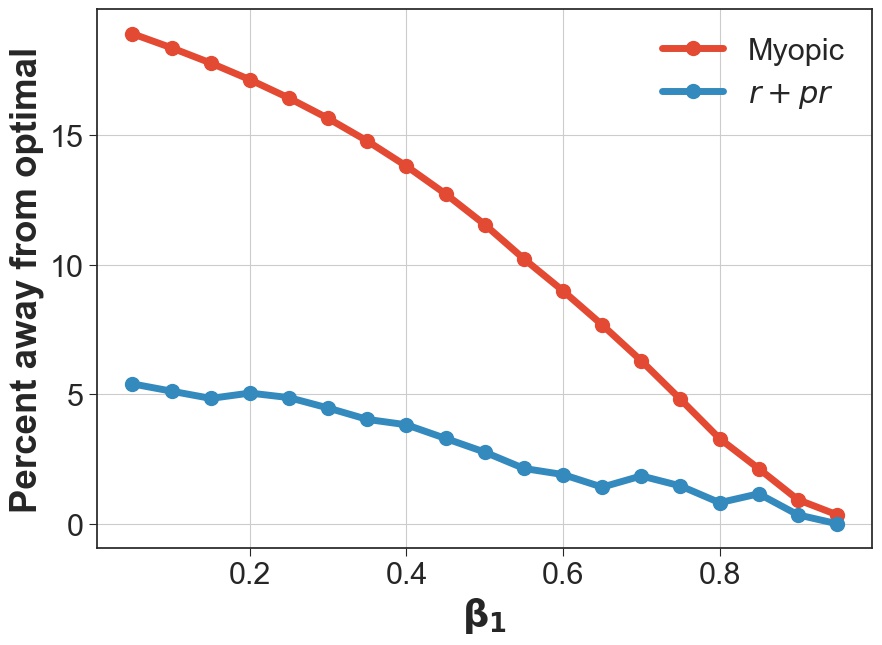} \\
\includegraphics[width=0.6\textwidth]{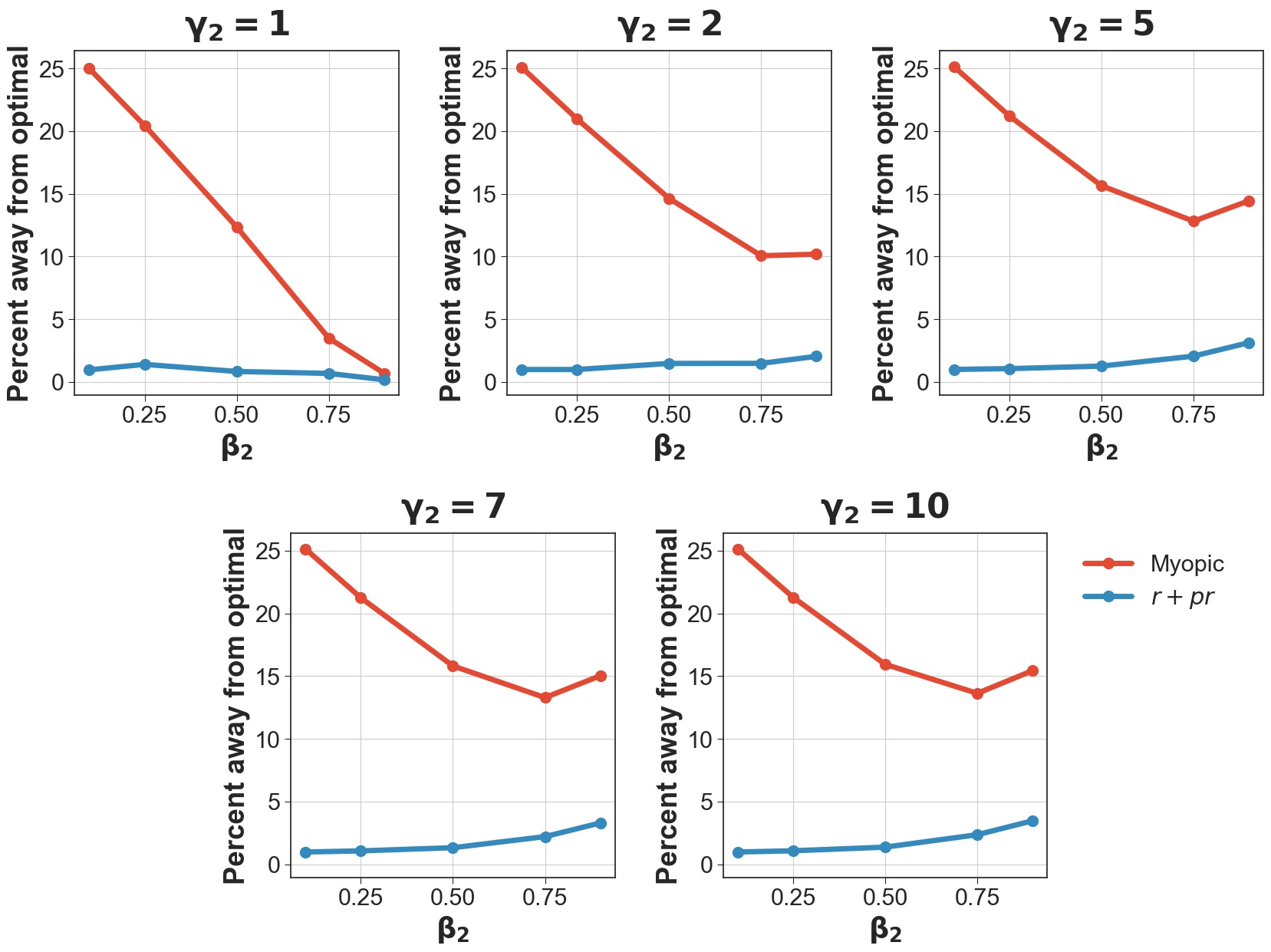} \\
\includegraphics[width=0.5\textwidth]{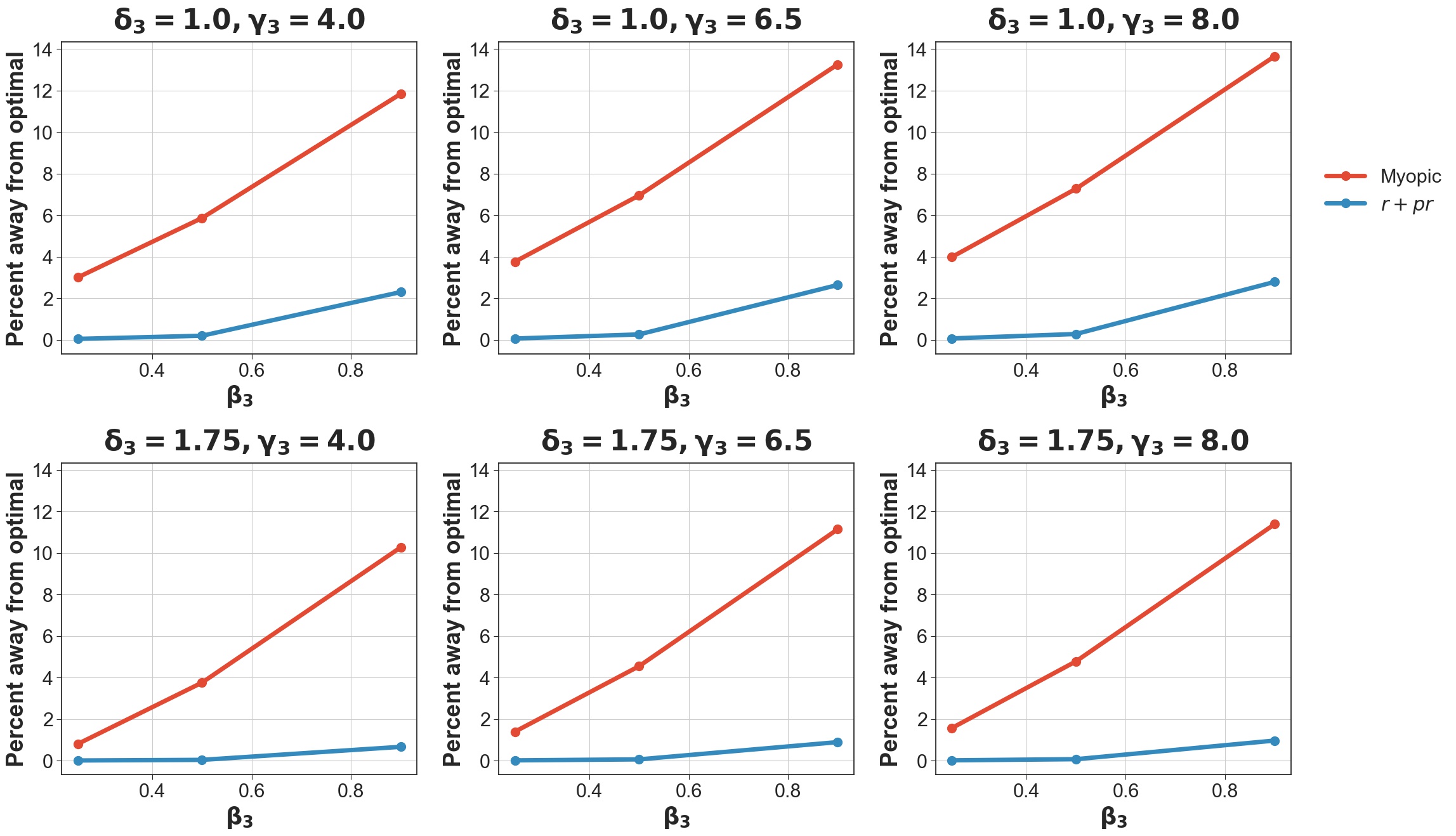}
\caption{For second instance considered when heuristic $r+p r$ performs well, the relative gap (in \%) of the expected long-run average costs between each heuristic (myopic $r$ and heuristic $r+pr$) under Scenarios 1--3.} \label{fig:instance4_myopic}
\end{figure}

We also used the same baseline transition matrices in this instance to investigate Scenarios 1--3 for all of the parameters listed in Table~\ref{tab:dependency_parameters}. Results from this experiment are summarized in Figure~\ref{fig:instance4_myopic}.

\subsubsection{Instance in which the optimal policy outperforms the myopic policy $r$ and policy $r + p r$} \label{appendix:7}

For this instance, we started with baseline transition matrices provided in Section~\ref{sec:numstudy} of the main text.  From these baseline matrices, we generated the final 25 2-by-2 transition matrices under Scenario 1 with $\beta_1=0.2$:

\footnotesize
\begin{align*}
\mathbf{P}^{1,1} = \begin{bmatrix}
0.56&0.44\\
0.99&0.01
\end{bmatrix}, & \mathbf{P}^{2,1} = \begin{bmatrix}
0.56&0.44\\
0.94&0.06
\end{bmatrix}, & \mathbf{P}^{3,1} = \begin{bmatrix}
0.56&0.44\\
0.94&0.06
\end{bmatrix}, & \mathbf{P}^{4,1} = \begin{bmatrix}
0.56&0.44\\
0.94&0.06
\end{bmatrix}, & \mathbf{P}^{5,1} = \begin{bmatrix}
0.56&0.44\\
0.94&0.06
\end{bmatrix},\\
\mathbf{P}^{1,2} = \begin{bmatrix}
0.91&0.09\\
0.78&0.22
\end{bmatrix}, & \mathbf{P}^{2,2} = \begin{bmatrix}
0.91&0.09\\
0.96&0.04
\end{bmatrix}, & \mathbf{P}^{3,2} = \begin{bmatrix}
0.91&0.09\\
0.78&0.22
\end{bmatrix}, & \mathbf{P}^{4,2} = \begin{bmatrix}
0.91&0.09\\
0.78&0.22
\end{bmatrix}, & \mathbf{P}^{5,2} = \begin{bmatrix}
0.91&0.09\\
0.78&0.22
\end{bmatrix},\\
\mathbf{P}^{1,3} = \begin{bmatrix}
0.36&0.64\\
0.68&0.32
\end{bmatrix}, &\mathbf{P}^{2,3} = \begin{bmatrix}
0.36&0.64\\
0.68&0.32
\end{bmatrix}, &\mathbf{P}^{3,3} = \begin{bmatrix}
0.36&0.64\\
0.94&0.06
\end{bmatrix}, &\mathbf{P}^{4,3} = \begin{bmatrix}
0.36&0.64\\
0.68&0.32
\end{bmatrix}, &\mathbf{P}^{5,3} = \begin{bmatrix}
0.36&0.64\\
0.68&0.32
\end{bmatrix},\\
\mathbf{P}^{1,4} = \begin{bmatrix}
0.97&0.03\\
0.02&0.98
\end{bmatrix}, &\mathbf{P}^{2,4} = \begin{bmatrix}
0.97&0.03\\
0.02&0.98
\end{bmatrix}, &\mathbf{P}^{3,4} = \begin{bmatrix}
0.97&0.03\\
0.02&0.98
\end{bmatrix}, &\mathbf{P}^{4,4} = \begin{bmatrix}
0.97&0.03\\
0.8&0.2
\end{bmatrix}, &\mathbf{P}^{5,4} = \begin{bmatrix}
0.97&0.03\\
0.02&0.98
\end{bmatrix},\\
\mathbf{P}^{1,5} = \begin{bmatrix}
0.88&0.12\\
0.01&0.99
\end{bmatrix}, &\mathbf{P}^{2,5} = \begin{bmatrix}
0.88&0.12\\
0.01&0.99
\end{bmatrix}, &\mathbf{P}^{3,5} = \begin{bmatrix}
0.88&0.12\\
0.01&0.99
\end{bmatrix}, &\mathbf{P}^{4,5} = \begin{bmatrix}
0.88&0.12\\
0.01&0.99
\end{bmatrix}, &\mathbf{P}^{5,5} = \begin{bmatrix}
0.88&0.12\\
0.8&0.2
\end{bmatrix}.
\end{align*}
\normalsize

\medskip

Note, these matrices above are identical to the corresponding baseline transition matrices except for entry $p^{a,a}_{1,1}$ which recovered by multiplying the corresponding entry in $\mathbf{\hat{P}}^a$ by $\beta_1$. Thus, since $\hat{p}^A_{1,1}=0.06$, then $p^{1,1}_{1,}=0.06\times0.2=0.01$. Based on these transition matrices, we recovered each of the three policies, which are described by action and state in Table~\ref{tab:ExampleAppendix8}.

We also investigated the influence of discharge rates on our results using the same transition probabilities. We recovered qualitatively similar long-run average costs when changing discharge rates to $\lambda_{UM}=0.4$, $\lambda_{CM}=0.3$, $\lambda_{CS}=0.2$, and $\lambda_{JS}=0.075$ and when the discharge rates are $\lambda_{UM}=0.025$, $\lambda_{CM}=0.075$, $\lambda_{CS}=0.2$, and $\lambda_{JS}=0.3$.  In the former case, the optimal average cost is approximately $8.26$, the average cost of the myopic policy is approximately $9.49$, and the average cost of heuristic $r + p r$ is approximately $8.73$.  It follows that the average costs of the myopic and $r+pr$ heuristics, respectively, are 14.9\% and 5.69\% higher than the optimal average cost.  In the latter case, the optimal average cost is approximately $21.81$, the average cost of the myopic policy is approximately $27.6$, and the average cost of heuristic $r + p r$ is approximately $23.98$.  It follows that the average costs of the myopic and $r+pr$ heuristics, respectively, are 26.5\% and 9.95\% higher than the optimal average cost.

\subsection{Policies accompanying qualitative analysis section}

\begin{center}
\begin{longtable}{|l| l | l | l | l | l | l | l | l |}
\caption{For first instance considered when myopic policy $r$ performs well, actions chosen in each state by the myopic, $r+pr$, and optimal policies.} \label{tab:ExampleAppendix4} \\

\hline \multicolumn{1}{|c|}{\textbf{$i$}} & \multicolumn{1}{c|}{\textbf{$s_1$}} & \multicolumn{1}{c|}{\textbf{$s_2$}} & 
\multicolumn{1}{c|}{\textbf{$s_3$}} & \multicolumn{1}{c|}{\textbf{$s_4$}} & 
\multicolumn{1}{c|}{\textbf{$s_5$}} & \multicolumn{1}{c|}{\textbf{$\mathbf{r}$}} & 
\multicolumn{1}{c|}{\textbf{$\mathbf{r+pr}$}} & \multicolumn{1}{c|}{\textbf{Optimal}}  \\ \hline 
\endfirsthead

\multicolumn{9}{c}%
{{\bfseries \tablename\ \thetable{} -- continued from previous page}} \\
\hline 
\multicolumn{1}{|c|}{\textbf{$i$}} & \multicolumn{1}{c|}{\textbf{$s_1$}} & \multicolumn{1}{c|}{\textbf{$s_2$}} & 
\multicolumn{1}{c|}{\textbf{$s_3$}} & \multicolumn{1}{c|}{\textbf{$s_4$}} & 
\multicolumn{1}{c|}{\textbf{$s_5$}} & \multicolumn{1}{c|}{\textbf{$\mathbf{r}$}} & 
\multicolumn{1}{c|}{\textbf{$\mathbf{r+pr}$}} & \multicolumn{1}{c|}{\textbf{Optimal}}  \\ \hline
\endhead

\hline \multicolumn{9}{|r|}{{Continued on next page}} \\ \hline
\endfoot

\hline \hline
\endlastfoot

1&0&0&0&0&0&0&0&0\\
1&0&0&0&0&1&5&5&5\\
1&0&0&0&1&0&4&4&4\\
1&0&0&0&1&1&4&4&4\\
1&0&0&1&0&0&3&3&3\\
1&0&0&1&0&1&3&3&3\\
1&0&0&1&1&0&4&4&4\\
1&0&0&1&1&1&4&4&4\\
1&0&1&0&0&0&2&2&2\\
1&0&1&0&0&1&2&2&2\\
1&0&1&0&1&0&4&4&4\\
1&0&1&0&1&1&4&4&4\\
1&0&1&1&0&0&3&3&3\\
1&0&1&1&0&1&3&3&3\\
1&0&1&1&1&0&4&4&4\\
1&0&1&1&1&1&4&4&4\\
1&1&0&0&0&0&1&1&1\\
1&1&0&0&0&1&1&5&1\\
1&1&0&0&1&0&4&4&4\\
1&1&0&0&1&1&4&4&4\\
1&1&0&1&0&0&1&3&1\\
1&1&0&1&0&1&1&3&1\\
1&1&0&1&1&0&4&4&4\\
1&1&0&1&1&1&4&4&4\\
1&1&1&0&0&0&1&2&1\\
1&1&1&0&0&1&1&1&1\\
1&1&1&0&1&0&4&4&4\\
1&1&1&0&1&1&4&4&4\\
1&1&1&1&0&0&1&3&1\\
1&1&1&1&0&1&1&1&1\\
1&1&1&1&1&0&4&4&4\\
1&1&1&1&1&1&4&4&4\\
2&0&0&0&0&0&0&0&0\\
2&0&0&0&0&1&5&5&5\\
2&0&0&0&1&0&4&4&4\\
2&0&0&0&1&1&5&5&5\\
2&0&0&1&0&0&3&3&3\\
2&0&0&1&0&1&5&5&5\\
2&0&0&1&1&0&4&4&4\\
2&0&0&1&1&1&5&5&5\\
2&0&1&0&0&0&2&2&2\\
2&0&1&0&0&1&5&5&5\\
2&0&1&0&1&0&4&4&4\\
2&0&1&0&1&1&5&5&5\\
2&0&1&1&0&0&3&3&3\\
2&0&1&1&0&1&5&5&5\\
2&0&1&1&1&0&4&4&4\\
2&0&1&1&1&1&5&5&5\\
2&1&0&0&0&0&1&1&1\\
2&1&0&0&0&1&5&5&5\\
2&1&0&0&1&0&4&4&4\\
2&1&0&0&1&1&5&5&5\\
2&1&0&1&0&0&1&3&1\\
2&1&0&1&0&1&5&5&5\\
2&1&0&1&1&0&4&4&4\\
2&1&0&1&1&1&5&5&5\\
2&1&1&0&0&0&1&2&1\\
2&1&1&0&0&1&5&5&5\\
2&1&1&0&1&0&4&4&1\\
2&1&1&0&1&1&5&5&5\\
2&1&1&1&0&0&1&1&1\\
2&1&1&1&0&1&5&5&5\\
2&1&1&1&1&0&4&4&1\\
2&1&1&1&1&1&5&5&5\\
3&0&0&0&0&0&0&0&0\\
3&0&0&0&0&1&5&5&5\\
3&0&0&0&1&0&4&4&4\\
3&0&0&0&1&1&4&4&4\\
3&0&0&1&0&0&3&3&3\\
3&0&0&1&0&1&5&5&5\\
3&0&0&1&1&0&4&4&4\\
3&0&0&1&1&1&4&4&4\\
3&0&1&0&0&0&2&2&2\\
3&0&1&0&0&1&5&5&5\\
3&0&1&0&1&0&4&4&4\\
3&0&1&0&1&1&4&4&4\\
3&0&1&1&0&0&2&2&3\\
3&0&1&1&0&1&5&5&5\\
3&0&1&1&1&0&4&4&4\\
3&0&1&1&1&1&4&4&4\\
3&1&0&0&0&0&1&1&1\\
3&1&0&0&0&1&5&5&5\\
3&1&0&0&1&0&4&4&4\\
3&1&0&0&1&1&4&4&4\\
3&1&0&1&0&0&1&3&1\\
3&1&0&1&0&1&5&5&1\\
3&1&0&1&1&0&4&4&4\\
3&1&0&1&1&1&4&4&4\\
3&1&1&0&0&0&1&2&1\\
3&1&1&0&0&1&5&5&1\\
3&1&1&0&1&0&4&4&4\\
3&1&1&0&1&1&4&4&4\\
3&1&1&1&0&0&1&3&1\\
3&1&1&1&0&1&5&5&1\\
3&1&1&1&1&0&4&4&4\\
3&1&1&1&1&1&4&4&4\\
4&0&0&0&0&0&0&0&0\\
4&0&0&0&0&1&5&5&5\\
4&0&0&0&1&0&4&4&4\\
4&0&0&0&1&1&5&5&5\\
4&0&0&1&0&0&3&3&3\\
4&0&0&1&0&1&5&5&5\\
4&0&0&1&1&0&4&3&4\\
4&0&0&1&1&1&5&5&5\\
4&0&1&0&0&0&2&2&2\\
4&0&1&0&0&1&5&5&5\\
4&0&1&0&1&0&4&4&4\\
4&0&1&0&1&1&5&5&5\\
4&0&1&1&0&0&3&3&3\\
4&0&1&1&0&1&5&5&5\\
4&0&1&1&1&0&4&3&4\\
4&0&1&1&1&1&5&5&5\\
4&1&0&0&0&0&1&1&1\\
4&1&0&0&0&1&5&5&5\\
4&1&0&0&1&0&1&4&4\\
4&1&0&0&1&1&5&5&5\\
4&1&0&1&0&0&1&3&1\\
4&1&0&1&0&1&5&5&5\\
4&1&0&1&1&0&1&3&1\\
4&1&0&1&1&1&5&5&5\\
4&1&1&0&0&0&1&2&1\\
4&1&1&0&0&1&5&5&5\\
4&1&1&0&1&0&1&4&1\\
4&1&1&0&1&1&5&5&5\\
4&1&1&1&0&0&1&3&1\\
4&1&1&1&0&1&5&5&5\\
4&1&1&1&1&0&1&3&1\\
4&1&1&1&1&1&5&5&5\\
\bottomrule
\end{longtable}
\end{center}

\begin{center}
\begin{longtable}{|l | l | l | l | l | l | l | l | l |}
\caption{Second instance in which  myopic policy $r$ performs well, actions chosen in each state by the myopic, $r+pr$, and optimal policies.} \label{tab:ExampleAppendix5} \\

\hline \multicolumn{1}{|c|}{\textbf{$i$}} & \multicolumn{1}{c|}{\textbf{$s_1$}} & \multicolumn{1}{c|}{\textbf{$s_2$}} & 
\multicolumn{1}{c|}{\textbf{$s_3$}} & \multicolumn{1}{c|}{\textbf{$s_4$}} & 
\multicolumn{1}{c|}{\textbf{$s_5$}} & \multicolumn{1}{c|}{\textbf{$\mathbf{r}$}} & 
\multicolumn{1}{c|}{\textbf{$\mathbf{r+pr}$}} & \multicolumn{1}{c|}{\textbf{Optimal}}  \\ \hline 
\endfirsthead

\multicolumn{9}{c}%
{{\bfseries \tablename\ \thetable{} -- continued from previous page}} \\
\hline 
\multicolumn{1}{|c|}{\textbf{$i$}} & \multicolumn{1}{c|}{\textbf{$s_1$}} & \multicolumn{1}{c|}{\textbf{$s_2$}} & 
\multicolumn{1}{c|}{\textbf{$s_3$}} & \multicolumn{1}{c|}{\textbf{$s_4$}} & 
\multicolumn{1}{c|}{\textbf{$s_5$}} & \multicolumn{1}{c|}{\textbf{$\mathbf{r}$}} & 
\multicolumn{1}{c|}{\textbf{$\mathbf{r+pr}$}} & \multicolumn{1}{c|}{\textbf{Optimal}}  \\ \hline
\endhead

\hline \multicolumn{9}{|r|}{{Continued on next page}} \\ \hline
\endfoot

\hline \hline
\endlastfoot

1&0&0&0&0&0&0&0&0\\
1&0&0&0&0&1&5&5&5\\
1&0&0&0&1&0&4&4&4\\
1&0&0&0&1&1&4&4&4\\
1&0&0&1&0&0&3&3&3\\
1&0&0&1&0&1&3&5&5\\
1&0&0&1&1&0&4&4&4\\
1&0&0&1&1&1&4&4&4\\
1&0&1&0&0&0&2&2&2\\
1&0&1&0&0&1&2&5&5\\
1&0&1&0&1&0&4&4&4\\
1&0&1&0&1&1&4&4&4\\
1&0&1&1&0&0&3&2&3\\
1&0&1&1&0&1&3&2&5\\
1&0&1&1&1&0&4&4&4\\
1&0&1&1&1&1&4&4&4\\
1&1&0&0&0&0&1&1&1\\
1&1&0&0&0&1&1&1&5\\
1&1&0&0&1&0&4&4&4\\
1&1&0&0&1&1&4&4&4\\
1&1&0&1&0&0&1&1&1\\
1&1&0&1&0&1&1&1&5\\
1&1&0&1&1&0&4&4&4\\
1&1&0&1&1&1&4&4&4\\
1&1&1&0&0&0&1&1&1\\
1&1&1&0&0&1&1&1&1\\
1&1&1&0&1&0&4&4&4\\
1&1&1&0&1&1&4&4&4\\
1&1&1&1&0&0&1&1&1\\
1&1&1&1&0&1&1&1&1\\
1&1&1&1&1&0&4&4&4\\
1&1&1&1&1&1&4&4&4\\
2&0&0&0&0&0&0&0&0\\
2&0&0&0&0&1&5&5&5\\
2&0&0&0&1&0&4&4&4\\
2&0&0&0&1&1&5&5&5\\
2&0&0&1&0&0&3&3&3\\
2&0&0&1&0&1&5&5&5\\
2&0&0&1&1&0&4&4&4\\
2&0&0&1&1&1&5&5&5\\
2&0&1&0&0&0&2&2&2\\
2&0&1&0&0&1&5&5&5\\
2&0&1&0&1&0&4&4&4\\
2&0&1&0&1&1&5&5&5\\
2&0&1&1&0&0&3&2&3\\
2&0&1&1&0&1&5&5&5\\
2&0&1&1&1&0&4&4&4\\
2&0&1&1&1&1&5&5&5\\
2&1&0&0&0&0&1&1&1\\
2&1&0&0&0&1&5&5&5\\
2&1&0&0&1&0&4&4&4\\
2&1&0&0&1&1&5&5&5\\
2&1&0&1&0&0&1&1&1\\
2&1&0&1&0&1&5&5&5\\
2&1&0&1&1&0&4&4&4\\
2&1&0&1&1&1&5&5&5\\
2&1&1&0&0&0&1&1&1\\
2&1&1&0&0&1&5&5&5\\
2&1&1&0&1&0&4&4&4\\
2&1&1&0&1&1&5&5&5\\
2&1&1&1&0&0&1&1&1\\
2&1&1&1&0&1&5&5&5\\
2&1&1&1&1&0&4&4&4\\
2&1&1&1&1&1&5&5&5\\
3&0&0&0&0&0&0&0&0\\
3&0&0&0&0&1&5&5&5\\
3&0&0&0&1&0&4&4&4\\
3&0&0&0&1&1&4&4&4\\
3&0&0&1&0&0&3&3&3\\
3&0&0&1&0&1&5&5&5\\
3&0&0&1&1&0&4&4&4\\
3&0&0&1&1&1&4&4&4\\
3&0&1&0&0&0&2&2&2\\
3&0&1&0&0&1&5&5&5\\
3&0&1&0&1&0&4&4&4\\
3&0&1&0&1&1&4&4&4\\
3&0&1&1&0&0&2&2&3\\
3&0&1&1&0&1&5&5&5\\
3&0&1&1&1&0&4&4&4\\
3&0&1&1&1&1&4&4&4\\
3&1&0&0&0&0&1&1&1\\
3&1&0&0&0&1&5&5&5\\
3&1&0&0&1&0&4&4&4\\
3&1&0&0&1&1&4&4&4\\
3&1&0&1&0&0&1&1&1\\
3&1&0&1&0&1&5&5&5\\
3&1&0&1&1&0&4&4&4\\
3&1&0&1&1&1&4&4&4\\
3&1&1&0&0&0&1&1&1\\
3&1&1&0&0&1&5&5&5\\
3&1&1&0&1&0&4&4&4\\
3&1&1&0&1&1&4&4&4\\
3&1&1&1&0&0&1&1&1\\
3&1&1&1&0&1&5&5&5\\
3&1&1&1&1&0&4&4&4\\
3&1&1&1&1&1&4&4&4\\
4&0&0&0&0&0&0&0&0\\
4&0&0&0&0&1&5&5&5\\
4&0&0&0&1&0&4&4&4\\
4&0&0&0&1&1&5&5&5\\
4&0&0&1&0&0&3&3&3\\
4&0&0&1&0&1&5&5&5\\
4&0&0&1&1&0&4&4&4\\
4&0&0&1&1&1&5&5&5\\
4&0&1&0&0&0&2&2&2\\
4&0&1&0&0&1&5&5&5\\
4&0&1&0&1&0&4&4&4\\
4&0&1&0&1&1&5&5&5\\
4&0&1&1&0&0&3&2&3\\
4&0&1&1&0&1&5&5&5\\
4&0&1&1&1&0&4&4&4\\
4&0&1&1&1&1&5&5&5\\
4&1&0&0&0&0&1&1&1\\
4&1&0&0&0&1&5&5&5\\
4&1&0&0&1&0&1&4&4\\
4&1&0&0&1&1&5&5&5\\
4&1&0&1&0&0&1&1&1\\
4&1&0&1&0&1&5&5&5\\
4&1&0&1&1&0&1&1&4\\
4&1&0&1&1&1&5&5&5\\
4&1&1&0&0&0&1&1&1\\
4&1&1&0&0&1&5&5&5\\
4&1&1&0&1&0&1&1&4\\
4&1&1&0&1&1&5&5&5\\
4&1&1&1&0&0&1&1&1\\
4&1&1&1&0&1&5&5&5\\
4&1&1&1&1&0&1&1&4\\
4&1&1&1&1&1&5&5&5
\end{longtable}
\end{center}

\begin{center}
\begin{longtable}{|l| l | l | l | l | l | l | l | l |}
\caption{First instance in which heuristic $r+pr$  performs well, actions chosen in each state by the myopic, $r+pr$, and optimal policies.} \label{tab:ExampleAppendix6} \\

\hline \multicolumn{1}{|c|}{\textbf{$i$}} & \multicolumn{1}{c|}{\textbf{$s_1$}} & \multicolumn{1}{c|}{\textbf{$s_2$}} & 
\multicolumn{1}{c|}{\textbf{$s_3$}} & \multicolumn{1}{c|}{\textbf{$s_4$}} & 
\multicolumn{1}{c|}{\textbf{$s_5$}} & \multicolumn{1}{c|}{\textbf{$\mathbf{r}$}} & 
\multicolumn{1}{c|}{\textbf{$\mathbf{r+pr}$}} & \multicolumn{1}{c|}{\textbf{Optimal}}  \\ \hline 
\endfirsthead

\multicolumn{9}{c}%
{{\bfseries \tablename\ \thetable{} -- continued from previous page}} \\
\hline 
\multicolumn{1}{|c|}{\textbf{$i$}} & \multicolumn{1}{c|}{\textbf{$s_1$}} & \multicolumn{1}{c|}{\textbf{$s_2$}} & 
\multicolumn{1}{c|}{\textbf{$s_3$}} & \multicolumn{1}{c|}{\textbf{$s_4$}} & 
\multicolumn{1}{c|}{\textbf{$s_5$}} & \multicolumn{1}{c|}{\textbf{$\mathbf{r}$}} & 
\multicolumn{1}{c|}{\textbf{$\mathbf{r+pr}$}} & \multicolumn{1}{c|}{\textbf{Optimal}}  \\ \hline
\endhead    

\hline \multicolumn{9}{|r|}{{Continued on next page}} \\ \hline
\endfoot

\hline \hline
\endlastfoot

1&0&0&0&0&0&0&0&0\\
1&0&0&0&0&1&5&5&5\\
1&0&0&0&1&0&4&4&4\\
1&0&0&0&1&1&4&4&4\\
1&0&0&1&0&0&3&3&3\\
1&0&0&1&0&1&3&3&3\\
1&0&0&1&1&0&4&4&4\\
1&0&0&1&1&1&4&4&4\\
1&0&1&0&0&0&2&2&2\\
1&0&1&0&0&1&2&2&2\\
1&0&1&0&1&0&4&4&4\\
1&0&1&0&1&1&4&4&4\\
1&0&1&1&0&0&3&3&3\\
1&0&1&1&0&1&3&3&3\\
1&0&1&1&1&0&4&4&4\\
1&0&1&1&1&1&4&4&4\\
1&1&0&0&0&0&1&1&1\\
1&1&0&0&0&1&1&1&1\\
1&1&0&0&1&0&4&1&1\\
1&1&0&0&1&1&4&1&1\\
1&1&0&1&0&0&1&1&1\\
1&1&0&1&0&1&1&1&1\\
1&1&0&1&1&0&4&1&1\\
1&1&0&1&1&1&4&1&1\\
1&1&1&0&0&0&1&1&1\\
1&1&1&0&0&1&1&1&1\\
1&1&1&0&1&0&4&1&1\\
1&1&1&0&1&1&4&1&1\\
1&1&1&1&0&0&1&1&1\\
1&1&1&1&0&1&1&1&1\\
1&1&1&1&1&0&4&1&1\\
1&1&1&1&1&1&4&1&1\\
2&0&0&0&0&0&0&0&0\\
2&0&0&0&0&1&5&5&5\\
2&0&0&0&1&0&4&4&4\\
2&0&0&0&1&1&5&5&5\\
2&0&0&1&0&0&3&3&3\\
2&0&0&1&0&1&5&5&5\\
2&0&0&1&1&0&4&4&4\\
2&0&0&1&1&1&5&5&5\\
2&0&1&0&0&0&2&2&2\\
2&0&1&0&0&1&5&5&5\\
2&0&1&0&1&0&4&4&4\\
2&0&1&0&1&1&5&5&5\\
2&0&1&1&0&0&3&3&3\\
2&0&1&1&0&1&5&5&5\\
2&0&1&1&1&0&4&4&4\\
2&0&1&1&1&1&5&5&5\\
2&1&0&0&0&0&1&1&1\\
2&1&0&0&0&1&5&5&5\\
2&1&0&0&1&0&4&1&1\\
2&1&0&0&1&1&5&5&5\\
2&1&0&1&0&0&1&1&1\\
2&1&0&1&0&1&5&5&5\\
2&1&0&1&1&0&4&1&1\\
2&1&0&1&1&1&5&5&5\\
2&1&1&0&0&0&1&1&1\\
2&1&1&0&0&1&5&5&5\\
2&1&1&0&1&0&4&1&1\\
2&1&1&0&1&1&5&5&5\\
2&1&1&1&0&0&1&1&1\\
2&1&1&1&0&1&5&5&5\\
2&1&1&1&1&0&4&1&1\\
2&1&1&1&1&1&5&5&5\\
3&0&0&0&0&0&0&0&0\\
3&0&0&0&0&1&5&5&5\\
3&0&0&0&1&0&4&4&4\\
3&0&0&0&1&1&4&4&4\\
3&0&0&1&0&0&3&3&3\\
3&0&0&1&0&1&5&5&5\\
3&0&0&1&1&0&4&4&4\\
3&0&0&1&1&1&4&4&4\\
3&0&1&0&0&0&2&2&2\\
3&0&1&0&0&1&5&5&5\\
3&0&1&0&1&0&4&4&4\\
3&0&1&0&1&1&4&4&4\\
3&0&1&1&0&0&2&2&2\\
3&0&1&1&0&1&5&5&5\\
3&0&1&1&1&0&4&4&4\\
3&0&1&1&1&1&4&4&4\\
3&1&0&0&0&0&1&1&1\\
3&1&0&0&0&1&5&5&1\\
3&1&0&0&1&0&4&1&1\\
3&1&0&0&1&1&4&5&1\\
3&1&0&1&0&0&1&1&1\\
3&1&0&1&0&1&5&5&1\\
3&1&0&1&1&0&4&1&1\\
3&1&0&1&1&1&4&5&1\\
3&1&1&0&0&0&1&1&1\\
3&1&1&0&0&1&5&5&1\\
3&1&1&0&1&0&4&1&1\\
3&1&1&0&1&1&4&5&1\\
3&1&1&1&0&0&1&1&1\\
3&1&1&1&0&1&5&5&1\\
3&1&1&1&1&0&4&1&1\\
3&1&1&1&1&1&4&5&1\\
4&0&0&0&0&0&0&0&0\\
4&0&0&0&0&1&5&5&5\\
4&0&0&0&1&0&4&4&4\\
4&0&0&0&1&1&5&5&5\\
4&0&0&1&0&0&3&3&3\\
4&0&0&1&0&1&5&5&5\\
4&0&0&1&1&0&4&3&4\\
4&0&0&1&1&1&5&5&5\\
4&0&1&0&0&0&2&2&2\\
4&0&1&0&0&1&5&5&5\\
4&0&1&0&1&0&4&2&4\\
4&0&1&0&1&1&5&5&5\\
4&0&1&1&0&0&3&3&3\\
4&0&1&1&0&1&5&5&5\\
4&0&1&1&1&0&4&3&4\\
4&0&1&1&1&1&5&5&5\\
4&1&0&0&0&0&1&1&1\\
4&1&0&0&0&1&5&5&5\\
4&1&0&0&1&0&1&1&1\\
4&1&0&0&1&1&5&5&5\\
4&1&0&1&0&0&1&1&1\\
4&1&0&1&0&1&5&5&5\\
4&1&0&1&1&0&1&1&1\\
4&1&0&1&1&1&5&5&5\\
4&1&1&0&0&0&1&1&1\\
4&1&1&0&0&1&5&5&5\\
4&1&1&0&1&0&1&1&1\\
4&1&1&0&1&1&5&5&5\\
4&1&1&1&0&0&1&1&1\\
4&1&1&1&0&1&5&5&5\\
4&1&1&1&1&0&1&1&1\\
4&1&1&1&1&1&5&5&5
\end{longtable}
\end{center}

\begin{center}
\begin{longtable}{|l| l | l | l | l | l | l | l | l |}
\caption{Second instance in which heuristic $r+pr$ performs well, actions chosen in each state by the myopic, $r+pr$, and optimal policies.} \label{tab:ExampleAppendix7} \\

\hline \multicolumn{1}{|c|}{\textbf{$i$}} & \multicolumn{1}{c|}{\textbf{$s_1$}} & \multicolumn{1}{c|}{\textbf{$s_2$}} & 
\multicolumn{1}{c|}{\textbf{$s_3$}} & \multicolumn{1}{c|}{\textbf{$s_4$}} & 
\multicolumn{1}{c|}{\textbf{$s_5$}} & \multicolumn{1}{c|}{\textbf{$\mathbf{r}$}} & 
\multicolumn{1}{c|}{\textbf{$\mathbf{r+pr}$}} & \multicolumn{1}{c|}{\textbf{Optimal}}  \\ \hline 
\endfirsthead

\multicolumn{9}{c}%
{{\bfseries \tablename\ \thetable{} -- continued from previous page}} \\
\hline 
\multicolumn{1}{|c|}{\textbf{$i$}} & \multicolumn{1}{c|}{\textbf{$s_1$}} & \multicolumn{1}{c|}{\textbf{$s_2$}} & 
\multicolumn{1}{c|}{\textbf{$s_3$}} & \multicolumn{1}{c|}{\textbf{$s_4$}} & 
\multicolumn{1}{c|}{\textbf{$s_5$}} & \multicolumn{1}{c|}{\textbf{$\mathbf{r}$}} & 
\multicolumn{1}{c|}{\textbf{$\mathbf{r+pr}$}} & \multicolumn{1}{c|}{\textbf{Optimal}}  \\ \hline
\endhead

\hline \multicolumn{9}{|r|}{{Continued on next page}} \\ \hline
\endfoot

\hline \hline
\endlastfoot

1&0&0&0&0&0&0&0&0\\
1&0&0&0&0&1&5&5&5\\
1&0&0&0&1&0&4&4&4\\
1&0&0&0&1&1&4&4&4\\
1&0&0&1&0&0&3&3&3\\
1&0&0&1&0&1&3&3&3\\
1&0&0&1&1&0&4&4&4\\
1&0&0&1&1&1&4&4&4\\
1&0&1&0&0&0&2&2&2\\
1&0&1&0&0&1&2&2&2\\
1&0&1&0&1&0&4&4&4\\
1&0&1&0&1&1&4&4&4\\
1&0&1&1&0&0&3&3&2\\
1&0&1&1&0&1&3&3&3\\
1&0&1&1&1&0&4&4&4\\
1&0&1&1&1&1&4&4&4\\
1&1&0&0&0&0&1&1&1\\
1&1&0&0&0&1&1&1&1\\
1&1&0&0&1&0&4&4&4\\
1&1&0&0&1&1&4&4&4\\
1&1&0&1&0&0&1&1&1\\
1&1&0&1&0&1&1&1&1\\
1&1&0&1&1&0&4&4&4\\
1&1&0&1&1&1&4&4&4\\
1&1&1&0&0&0&1&1&1\\
1&1&1&0&0&1&1&1&1\\
1&1&1&0&1&0&4&4&4\\
1&1&1&0&1&1&4&4&4\\
1&1&1&1&0&0&1&1&1\\
1&1&1&1&0&1&1&1&1\\
1&1&1&1&1&0&4&4&4\\
1&1&1&1&1&1&4&4&4\\
2&0&0&0&0&0&0&0&0\\
2&0&0&0&0&1&5&5&5\\
2&0&0&0&1&0&4&4&4\\
2&0&0&0&1&1&5&5&4\\
2&0&0&1&0&0&3&3&3\\
2&0&0&1&0&1&5&5&3\\
2&0&0&1&1&0&4&4&4\\
2&0&0&1&1&1&5&5&4\\
2&0&1&0&0&0&2&2&2\\
2&0&1&0&0&1&5&5&2\\
2&0&1&0&1&0&4&4&4\\
2&0&1&0&1&1&5&4&4\\
2&0&1&1&0&0&3&3&2\\
2&0&1&1&0&1&5&5&2\\
2&0&1&1&1&0&4&4&4\\
2&0&1&1&1&1&5&4&4\\
2&1&0&0&0&0&1&1&1\\
2&1&0&0&0&1&5&5&1\\
2&1&0&0&1&0&4&4&4\\
2&1&0&0&1&1&5&4&4\\
2&1&0&1&0&0&1&1&1\\
2&1&0&1&0&1&5&5&1\\
2&1&0&1&1&0&4&4&4\\
2&1&0&1&1&1&5&5&4\\
2&1&1&0&0&0&1&1&1\\
2&1&1&0&0&1&5&1&1\\
2&1&1&0&1&0&4&4&4\\
2&1&1&0&1&1&5&4&4\\
2&1&1&1&0&0&1&1&1\\
2&1&1&1&0&1&5&1&1\\
2&1&1&1&1&0&4&4&4\\
2&1&1&1&1&1&5&4&4\\
3&0&0&0&0&0&0&0&0\\
3&0&0&0&0&1&5&5&5\\
3&0&0&0&1&0&4&4&4\\
3&0&0&0&1&1&4&4&4\\
3&0&0&1&0&0&3&3&3\\
3&0&0&1&0&1&5&5&3\\
3&0&0&1&1&0&4&4&4\\
3&0&0&1&1&1&4&4&4\\
3&0&1&0&0&0&2&2&2\\
3&0&1&0&0&1&5&2&2\\
3&0&1&0&1&0&4&4&4\\
3&0&1&0&1&1&4&4&4\\
3&0&1&1&0&0&2&2&2\\
3&0&1&1&0&1&5&2&2\\
3&0&1&1&1&0&4&4&4\\
3&0&1&1&1&1&4&4&4\\
3&1&0&0&0&0&1&1&1\\
3&1&0&0&0&1&5&1&1\\
3&1&0&0&1&0&4&4&4\\
3&1&0&0&1&1&4&4&4\\
3&1&0&1&0&0&1&1&1\\
3&1&0&1&0&1&5&1&1\\
3&1&0&1&1&0&4&4&4\\
3&1&0&1&1&1&4&4&4\\
3&1&1&0&0&0&1&1&2\\
3&1&1&0&0&1&5&1&1\\
3&1&1&0&1&0&4&4&4\\
3&1&1&0&1&1&4&4&4\\
3&1&1&1&0&0&1&1&2\\
3&1&1&1&0&1&5&1&1\\
3&1&1&1&1&0&4&4&4\\
3&1&1&1&1&1&4&4&4\\
4&0&0&0&0&0&0&0&0\\
4&0&0&0&0&1&5&5&5\\
4&0&0&0&1&0&4&4&4\\
4&0&0&0&1&1&5&5&4\\
4&0&0&1&0&0&3&3&3\\
4&0&0&1&0&1&5&5&3\\
4&0&0&1&1&0&4&4&4\\
4&0&0&1&1&1&5&5&4\\
4&0&1&0&0&0&2&2&2\\
4&0&1&0&0&1&5&5&2\\
4&0&1&0&1&0&4&4&4\\
4&0&1&0&1&1&5&4&4\\
4&0&1&1&0&0&3&3&2\\
4&0&1&1&0&1&5&5&2\\
4&0&1&1&1&0&4&4&4\\
4&0&1&1&1&1&5&4&4\\
4&1&0&0&0&0&1&1&1\\
4&1&0&0&0&1&5&5&1\\
4&1&0&0&1&0&1&4&4\\
4&1&0&0&1&1&5&5&1\\
4&1&0&1&0&0&1&1&1\\
4&1&0&1&0&1&5&5&1\\
4&1&0&1&1&0&1&1&4\\
4&1&0&1&1&1&5&5&1\\
4&1&1&0&0&0&1&1&2\\
4&1&1&0&0&1&5&5&1\\
4&1&1&0&1&0&1&4&4\\
4&1&1&0&1&1&5&1&4\\
4&1&1&1&0&0&1&1&2\\
4&1&1&1&0&1&5&5&1\\
4&1&1&1&1&0&1&4&4\\
4&1&1&1&1&1&5&1&4
\end{longtable}
\end{center}

\begin{center}
\begin{longtable}{|l| l | l | l | l | l | l | l | l |}
\caption{Instance in which the optimal policy outperforms heuristics $r$ and $r+pr$, actions chosen in each state by the myopic, $r+pr$, and optimal policies.} \label{tab:ExampleAppendix8} \\

\hline \multicolumn{1}{|c|}{\textbf{$i$}} & \multicolumn{1}{c|}{\textbf{$s_1$}} & \multicolumn{1}{c|}{\textbf{$s_2$}} & 
\multicolumn{1}{c|}{\textbf{$s_3$}} & \multicolumn{1}{c|}{\textbf{$s_4$}} & 
\multicolumn{1}{c|}{\textbf{$s_5$}} & \multicolumn{1}{c|}{\textbf{$\mathbf{r}$}} & 
\multicolumn{1}{c|}{\textbf{$\mathbf{r+pr}$}} & \multicolumn{1}{c|}{\textbf{Optimal}}  \\ \hline 
\endfirsthead

\multicolumn{9}{c}%
{{\bfseries \tablename\ \thetable{} -- continued from previous page}} \\
\hline 
\multicolumn{1}{|c|}{\textbf{$i$}} & \multicolumn{1}{c|}{\textbf{$s_1$}} & \multicolumn{1}{c|}{\textbf{$s_2$}} & 
\multicolumn{1}{c|}{\textbf{$s_3$}} & \multicolumn{1}{c|}{\textbf{$s_4$}} & 
\multicolumn{1}{c|}{\textbf{$s_5$}} & \multicolumn{1}{c|}{\textbf{$\mathbf{r}$}} & 
\multicolumn{1}{c|}{\textbf{$\mathbf{r+pr}$}} & \multicolumn{1}{c|}{\textbf{Optimal}}  \\ \hline
\endhead

\hline \multicolumn{9}{|r|}{{Continued on next page}} \\ \hline
\endfoot

\hline \hline
\endlastfoot

1&0&0&0&0&0&0&0&0\\
1&0&0&0&0&1&5&5&5\\
1&0&0&0&1&0&4&4&4\\
1&0&0&0&1&1&4&4&4\\
1&0&0&1&0&0&3&3&3\\
1&0&0&1&0&1&3&3&3\\
1&0&0&1&1&0&4&3&3\\
1&0&0&1&1&1&4&4&3\\
1&0&1&0&0&0&2&2&2\\
1&0&1&0&0&1&2&2&2\\
1&0&1&0&1&0&4&2&2\\
1&0&1&0&1&1&4&4&2\\
1&0&1&1&0&0&3&2&3\\
1&0&1&1&0&1&3&3&3\\
1&0&1&1&1&0&4&3&3\\
1&0&1&1&1&1&4&4&3\\
1&1&0&0&0&0&1&1&1\\
1&1&0&0&0&1&1&1&1\\
1&1&0&0&1&0&4&1&1\\
1&1&0&0&1&1&4&4&1\\
1&1&0&1&0&0&1&1&1\\
1&1&0&1&0&1&1&1&1\\
1&1&0&1&1&0&4&1&1\\
1&1&0&1&1&1&4&4&1\\
1&1&1&0&0&0&1&1&1\\
1&1&1&0&0&1&1&1&1\\
1&1&1&0&1&0&4&1&1\\
1&1&1&0&1&1&4&4&1\\
1&1&1&1&0&0&1&1&1\\
1&1&1&1&0&1&1&1&1\\
1&1&1&1&1&0&4&1&1\\
1&1&1&1&1&1&4&4&1\\
2&0&0&0&0&0&0&0&0\\
2&0&0&0&0&1&5&5&5\\
2&0&0&0&1&0&4&4&4\\
2&0&0&0&1&1&5&5&5\\
2&0&0&1&0&0&3&3&3\\
2&0&0&1&0&1&5&3&3\\
2&0&0&1&1&0&4&3&3\\
2&0&0&1&1&1&5&5&3\\
2&0&1&0&0&0&2&2&2\\
2&0&1&0&0&1&5&2&2\\
2&0&1&0&1&0&4&2&2\\
2&0&1&0&1&1&5&5&2\\
2&0&1&1&0&0&3&2&3\\
2&0&1&1&0&1&5&3&3\\
2&0&1&1&1&0&4&3&3\\
2&0&1&1&1&1&5&5&3\\
2&1&0&0&0&0&1&1&1\\
2&1&0&0&0&1&5&1&1\\
2&1&0&0&1&0&4&1&1\\
2&1&0&0&1&1&5&5&1\\
2&1&0&1&0&0&1&1&1\\
2&1&0&1&0&1&5&1&1\\
2&1&0&1&1&0&4&1&1\\
2&1&0&1&1&1&5&5&1\\
2&1&1&0&0&0&1&1&1\\
2&1&1&0&0&1&5&1&1\\
2&1&1&0&1&0&4&1&1\\
2&1&1&0&1&1&5&5&1\\
2&1&1&1&0&0&1&1&1\\
2&1&1&1&0&1&5&1&1\\
2&1&1&1&1&0&4&1&1\\
2&1&1&1&1&1&5&5&1\\
3&0&0&0&0&0&0&0&0\\
3&0&0&0&0&1&5&5&5\\
3&0&0&0&1&0&4&4&4\\
3&0&0&0&1&1&4&4&4\\
3&0&0&1&0&0&3&3&3\\
3&0&0&1&0&1&5&3&3\\
3&0&0&1&1&0&4&3&3\\
3&0&0&1&1&1&4&4&3\\
3&0&1&0&0&0&2&2&2\\
3&0&1&0&0&1&5&2&2\\
3&0&1&0&1&0&4&2&2\\
3&0&1&0&1&1&4&4&2\\
3&0&1&1&0&0&2&2&2\\
3&0&1&1&0&1&5&2&2\\
3&0&1&1&1&0&4&2&2\\
3&0&1&1&1&1&4&4&2\\
3&1&0&0&0&0&1&1&1\\
3&1&0&0&0&1&5&1&1\\
3&1&0&0&1&0&4&1&1\\
3&1&0&0&1&1&4&4&1\\
3&1&0&1&0&0&1&1&1\\
3&1&0&1&0&1&5&1&1\\
3&1&0&1&1&0&4&1&1\\
3&1&0&1&1&1&4&4&1\\
3&1&1&0&0&0&1&1&1\\
3&1&1&0&0&1&5&1&1\\
3&1&1&0&1&0&4&1&1\\
3&1&1&0&1&1&4&4&1\\
3&1&1&1&0&0&1&1&1\\
3&1&1&1&0&1&5&1&1\\
3&1&1&1&1&0&4&1&1\\
3&1&1&1&1&1&4&4&1\\
4&0&0&0&0&0&0&0&0\\
4&0&0&0&0&1&5&5&5\\
4&0&0&0&1&0&4&4&4\\
4&0&0&0&1&1&5&5&5\\
4&0&0&1&0&0&3&3&3\\
4&0&0&1&0&1&5&3&3\\
4&0&0&1&1&0&4&3&3\\
4&0&0&1&1&1&5&5&3\\
4&0&1&0&0&0&2&2&2\\
4&0&1&0&0&1&5&2&2\\
4&0&1&0&1&0&4&2&2\\
4&0&1&0&1&1&5&5&2\\
4&0&1&1&0&0&3&2&3\\
4&0&1&1&0&1&5&3&3\\
4&0&1&1&1&0&4&3&3\\
4&0&1&1&1&1&5&5&3\\
4&1&0&0&0&0&1&1&1\\
4&1&0&0&0&1&5&1&1\\
4&1&0&0&1&0&1&1&1\\
4&1&0&0&1&1&5&5&1\\
4&1&0&1&0&0&1&1&1\\
4&1&0&1&0&1&5&1&1\\
4&1&0&1&1&0&1&1&1\\
4&1&0&1&1&1&5&5&1\\
4&1&1&0&0&0&1&1&1\\
4&1&1&0&0&1&5&1&1\\
4&1&1&0&1&0&1&1&1\\
4&1&1&0&1&1&5&5&1\\
4&1&1&1&0&0&1&1&1\\
4&1&1&1&0&1&5&1&1\\
4&1&1&1&1&0&1&1&1\\
4&1&1&1&1&1&5&5&1
\end{longtable}
\end{center}

\end{document}